\documentclass[preprint,3p,sort&compress,final,times]{elsarticle}

\usepackage{amssymb}
\usepackage{amsthm}
\usepackage{xcolor}
\usepackage{listings}
\usepackage{adjustbox}

\definecolor{codegreen}{rgb}{0,0.6,0}
\definecolor{codegray}{rgb}{0.5,0.5,0.5}
\definecolor{codepurple}{rgb}{0.58,0,0.82}
\definecolor{backcolour}{rgb}{0.95,0.95,0.92}

\lstdefinestyle{mystyle}{
    backgroundcolor=\color{backcolour},   
    commentstyle=\color{codegreen},
    keywordstyle=\color{magenta},
    numberstyle=\tiny\color{codegray},
    stringstyle=\color{codepurple},
    basicstyle=\ttfamily\footnotesize,
    breakatwhitespace=true,         
    breaklines=true,                 
    captionpos=b,                    
    keepspaces=true,                 
    numbers=left,                    
    numbersep=1pt,                  
    showspaces=false,                
    showstringspaces=false,
    showtabs=false,                  
    tabsize=2
}
\usepackage{bm}
\usepackage{amsmath}
\usepackage{graphicx}
\usepackage{epstopdf, epsfig}
\usepackage{amsmath}
\usepackage{amssymb}
\usepackage{float}
\usepackage{subcaption}
\usepackage{hyperref}
\usepackage{booktabs}
\usepackage{todonotes}
\usepackage{physics}

\usepackage{soul}
\hypersetup{
  colorlinks   = true, 
  urlcolor     = blue, 
  linkcolor    = blue, 
  citecolor   = blue 
}
\usepackage{scalerel,tikz}
\newtheorem{theorem}{Theorem}

\newtheorem{lemma}[theorem]{Lemma}

\newtheorem{remark}{Remark}

\newtheorem{definition}{Definition}[section]

\newcommand{\revi}[1]{#1}
\newcommand{\revii}[1]{#1}

\journal{Journal of Computational Physics}

\begin{document}

\begin{frontmatter}



\title{Strongly stable dual-pairing summation by parts finite difference schemes for the vector invariant  nonlinear shallow water equations -- I: Numerical scheme and validation on the plane}


\author[inst1,inst2,inst3]{Justin Kin Jun Hew}

\affiliation[inst1]{organization={Mathematical Sciences Institute, Australian National University},
            city={ACT},
            postcode={2601},
            country={Australia}}
\affiliation[inst2]{organization={Research School of Physics, Australian National University},
            city={ACT},
            postcode={2601},
            country={Australia}}
\affiliation[inst3]{organization={Research School of Astronomy and Astrophysics, Australian National University},
            city={ACT},
            postcode={2611},
            country={Australia}}

\affiliation[inst4]{organization={Department of Mathematical Sciences , University of Texas at El Paso},
            city={El Paso, Texas},
            postcode={79968},
            country={USA}}

\author[inst1,inst4]{Kenneth Duru}
\author[inst1]{Stephen Roberts}
\author[inst1]{Christopher Zoppou}
\author[inst1]{Kieran Ricardo}

\begin{abstract}
We present an energy/entropy stable and high order accurate finite difference (FD) method for solving the nonlinear (rotating) shallow water equations (SWEs) in vector invariant form using the newly developed dual-pairing and dispersion-relation preserving summation by parts (SBP)  FD  operators. We derive new well-posed boundary conditions (BCs) for the SWE in one space dimension, formulated in terms of fluxes and applicable to linear and nonlinear SWEs. For the nonlinear \revii{vector invariant SWE in the subcritical regime}, \revii{where energy is an entropy functional, we find that} energy/entropy stability ensures the boundedness of numerical solution but does not guarantee convergence. Adequate amount of numerical dissipation is necessary to control high frequency errors which could negatively impact accuracy in the numerical simulations. Using the dual-pairing SBP framework, we derive high order accurate and nonlinear hyper-viscosity operator which dissipates entropy and enstrophy. The hyper-viscosity operator effectively minimises oscillations from shocks and discontinuities, and eliminates high frequency grid-scale errors. The numerical method is most suitable for the simulations of subcritical flows typically observed in atmospheric and geostrophic flow problems. We prove both nonlinear and local linear stability  results, as well as a priori error estimates for the semi-discrete approximations of both linear and nonlinear SWEs.  Convergence, accuracy, and well-balanced properties are verified via the method of manufactured solutions and canonical test problems such as the dam break and lake at rest. Numerical simulations in two-dimensions are presented which include the rotating and merging vortex problem and barotropic shear instability, with fully developed turbulence.
\end{abstract}




\end{frontmatter}

\section{Introduction}
The depth-averaged shallow water equations (SWE) are a set of hyperbolic partial differential equations (PDE) governing incompressible fluid flow subject to small wavelength amplitude disturbances. They were first derived by Saint-Venant in 1871 using the Navier-Stokes (momentum) equation \revii{with the assumption that the depth of the fluid layer is negligible compared to the horizontal scale of motion}. \revii{This core assumption makes the model well-suited to simulate various types of shallow flows, such as riverine flooding \cite{tan1992shallow}, dam breaks \cite{cozzolino2018solution}, coastal floods \cite{mignot2006modeling}, estuary and lake flows \cite{sadourny1975dynamics}, and even for simulating weak oblique shock-wave reflection \cite{defina2008numerical}}.  Moreover, they are also widely used to model atmospheric flows~\cite{williamson1992standard,behrens1998atmospheric}, particularly with respect to the motion constrained near the thin fluid layer on Earth, which is well-known to behave as a nearly incompressible fluid, with only minor changes in density. \revii{The nonlinear SWEs are often  derived from first principles in conservative form,  using conservation of mass and momentum. However, under the assumption of the smoothness of the solutions the equations can also be re-written, in the vector invariant form, to evolve the primitive variables, that is mass and velocity.} One unique aspect of these equations compared to the traditional quasi-geostrophic equations is that it allows for high-fidelity prediction of linear dispersive waves in high-altitude atmospheric flows, often characterised by large Rossby number and the presence of large-scale baroclinic vorticity. \revi{These effects, are well-known} to have a major influence on weather (see e.g., the role of atmospheric Rossby waves in stratospheric warming \cite{chandran2013secondary}), and must therefore be accounted for in weather forecast models. Moreover, the SWE are    also  important for the simulation of Rossby waves in astrophysical flows (see e.g., review by \cite{zaqarashvili2021rossby}).
Therefore, provably accurate and efficient numerical methods for reliable simulations of the nonlinear SWE is important in many applications.

The properties of fluids that are described by the SWE can be characterised by a dimensionless number called the \textit{Froude number}, defined as $\textrm{Fr} = |\mathbf{u}|/\sqrt{gh}$, where $\mathbf{u}$ is the depth-averaged fluid velocity field, $g$ is gravitational acceleration and $h$ is the water depth, \revii{and from physical considerations, $h$ must be positive, $h>0$}. The fluid flow is called subcritical when $\textrm{Fr} < 1$, critical when $\textrm{Fr} = 1$ and super-critical when $\textrm{Fr} >1$. The different regimes change the fundamental nature of flow behaviour. Therefore, numerical schemes solving the SWE between regimes require careful treatments in order to capture important features of a target application. However, in this study we are particularly interested in  subcritical flows, with $\textrm{Fr} < 1$, which are commonly observed in atmospheric and geostrophic flow problems.
 
Global-scale atmospheric flows modelled by the SWE are often posed on the surface of the sphere where periodic BCs (BCs) suffice for a well-posed model and stable numerical implementations using cubed sphere meshes \cite{thuburn2012framework,shipton2018higher,lee2018mixed}. For regional-scale atmospheric models \cite{williamson1992standard}, oceanic flow models \cite{zeitlin2018geophysical} and many other modelling scenarios \cite{o2011consistent,richards2011appropriate} non-periodic and well-posed   BCs with stable numerical implementations are crucial for accurate and reliable simulations. \revii{For instance, in ocean models, Kelvin waves are known to occur on lateral boundaries. As such, lateral BCs are required in order to simulate such effects, as well as to resolve the associated vortical motion  \cite{hsieh1983free,beletsky1997numerical}.  A typical requirement to resolve these eddies is a resolution of the order of the Rossby radius of deformation, which are $\sim 100$ km in the atmosphere and $\sim 10$ km in the ocean} \revii{\cite{wang2002kelvin}}. \revii{Thus,} the derivation of nonlinear and entropy stable BCs is a necessary step towards the development of provably stable numerical methods for the nonlinear SWE.   
 
 Many prior studies \cite{lundgren2020efficient,GHADER20141,shamsnia2022comparative,nordstrom2022linear} have performed linear analysis, where the assumption of smooth solutions is then used to simulate nonlinear problems. The analysis of nonlinear hyperbolic IBVPs, without linearisation, is technically daunting. \revii{There are recent efforts, however, (see, e.g., \cite{nordstrom2022linear,nordstrom2022nonlinear,nordstrom2024nonlinear}) where nonlinear analysis are considered for the SWE subject to nonlinear BCs. However, to the best of our knowledge there is no  numerical evidence yet verifying the theory.
 One objective of this study is the derivation of  well-posed nonlinear BCs and stable numerical implementations for the nonlinear SWE in vector invariant form \cite{gilman2000magnetohydrodynamic,bihlo2012invariant,korn2018conservative,ricardo2023conservation,shashkin2023summation}}. It is also significantly noteworthy that a nonlinear analysis must be consistent with the linear theory, since for smooth solutions the error equation required for  convergence proofs is a linear system of equations. \revi{We note however that this may contradict some of the results derived in \cite{nordstrom2022linear} where it was shown that the number of  BCs for the nonlinear SWEs  is independent of the magnitude of the flow, which disagrees with the linear theory \cite{GHADER20141} for critical and super-critical flows. Nevertheless, for sub-critical flows considered in this study the number of BCs is in agreement with the results of \cite{nordstrom2022linear}.}
 For subcritical flows in one space dimension (1D) we derive new well-posed BCs for the vector invariant form of the SWE formulated in terms of fluxes, and applicable to linear and nonlinear SWE.


Traditional approaches to solving the SWE rely on either Godunov-type approximate Riemann solvers \cite{lu2020well}, which are typically first or second order accurate in space, and low-order finite difference (FD)  methods \cite{ranvcic1996global,thuburn2009numerical}. These schemes are low-order accurate and cannot be expected to accurately capture highly dispersive wave-dominated phenomena, e.g., gravitational or Rossby wave propagation. Recent  trends for the development of robust and high order accurate numerical methods for the nonlinear SWE are variational based schemes such as discontinuous Galerkin (DG) and continuous Galerkin (CG) finite/spectral element methods \cite{gassner2016well, giraldo2020introduction, ricardo2023conservation}. While DG and CG finite/spectral element schemes are flexible for resolving complex geometries, however, they are often designed for the conservative form of the equations (with the unknown variables $h$ and $h\mathbf{u}$) and may not discretely preserve some important invariants such as potential/absolute  vorticity and enstrophy which are critical for accurate simulations of atmospheric/oceanic flows.  Meanwhile, recent developments show that DG and CG finite/spectral element schemes formulated for the vector invariant form of the SWE (with the unknown variables  $h$ and $\mathbf{u}$) can be designed to be mimetic, that is they preserve several important invariants imposed by the system \cite{lee2018discrete, ricardo2023conservation, lee2022comparison,giraldo2002nodal,taylor2010compatible}. However, in smooth and piecewise smooth domains, stable and high order accurate FD (FD) schemes are attractive because they are more computationally efficient, see e.g., \cite{duru2022dual}. Nevertheless, the development of provably stable, high order accurate  and mimetic FD schemes for the nonlinear SWE is a challenge. Similar to DG and CG schemes, several FD schemes for the SWE are derived for the conservative form of the equations and may not discretely preserve some important invariants, such as vorticity and linear geostrophic balance, in the system \cite{lundgren2020efficient}. The ultimate goal of this initial paper is to develop provably accurate and mimetic  FD methods for the vector invariant form of the SWE, in 1D and  2D  plane, with solid mathematical support. Part II of this paper will extend the method to 2D curvilinear meshes and manifolds $\Omega \subset \mathbb{R}^2$ embedded in $\mathbb{R}^3$.

A necessary ingredient for the development of robust numerical methods is the summation by parts (SBP) principle \cite{svard2014review,fernandez2014review}. An important property of SBP methods is its mimetic structure, which allows provably stable methods to be constructed at the semi-discrete level, provided stable and accurate treatment of the BCs can be derived. Traditional SBP FD operators are central FD stencils on co-located grids with one-sided boundary stencil closures designed such that the FD operator obey a discrete integration-by-parts principle.  Generally central FD stencils suffer from spurious unresolved wave-modes, which can destroy the accuracy of numerical solutions for a nonlinear first-order  hyperbolic PDE, such as the SWE.
The dual-pairing (DP) SBP framework \cite{mattsson2017diagonal,duru2022dual,williams2021provably} was recently introduced to improve the accuracy of FD methods for wave problems in complex geometries.  The DP-SBP operators are a pair of high order (backward and forward) difference stencils which together preserve the desired properties of SBP methods. The DP-SBP operators  can also include additional degrees of freedom that can be tuned to diminish numerical dispersion errors, yielding the so-called dispersion-relation-preserving (DRP) DP-SBP FD operators \cite{williams2021provably}.
However, the DP-SBP operators are asymmetric and dissipative, can potentially destroy symmetries that exist in the continuum problem.

In this paper we derive provably  stable, high order accurate and mimetic FD schemes for the nonlinear vector invariant SWE using traditional SBP, DP-SBP and DRP  DP-SBP FD operators \cite{mattsson2017diagonal,williams2021provably}.   
To begin, we carefully combine the DP SBP operator pairs  so that we preserve the discrete skew-symmetric property and several invariants of the SWE. \revii{As already mentioned earlier, the invariants of interest in the atmospheric flow context are local conservation of mass, absolute vorticity and energy \cite{thuburn2008some}, where the latter defines an entropy functional in the SWE.  We further implement BCs} weakly using penalties such that semi-discrete approximation is energy/entropy stable.  In periodic domains, the proof of energy/entropy conservation follows in a straightforward manner.
For nonlinear problems, although energy/entropy conservation may ensure the boundedness of numerical solutions, however, it does not guarantee convergence of numerical errors with mesh refinement \revi{\cite{gassner2022stability}}. Suitable amount of numerical dissipation is necessary to control high frequency oscillations in the numerical simulations. Using the DP SBP framework, we design high order accurate and nonlinear hyper-viscosity operator which dissipates entropy. The hyper-viscosity operator effectively controls oscillations from shocks and discontinuities and eliminates high frequency grid-scale errors. %
We prove both nonlinear (energy/entropy) and local linear stability  results, and derive a priori error estimates for smooth solutions. The stability analysis is verified by computing the numerical eigenvalues of the evolution operator, for both linear and nonlinear discrete problems. 
The scheme is suitable for the simulations of  subcritical flow regimes typical prevalent in  atmospheric and geostrophic flow problems.
The method has been validated with a number of canonical test cases as well as with the method of manufactured smooth solutions, including 2D merging vortex problems and barotropic shear instability \cite{galewsky2004initial}, with fully developed turbulence.

The rest of the paper is organised as follows. In Section~\ref{sec:contanalysis}, we perform continuous analysis of the rotational SWEs in vector invariant form, before introducing numerical preliminaries in Section~\ref{sec:sbp} and semi-discrete analysis in Section~\ref{sec:semiapprox}. Then, we derive error convergence rates in Section~\ref{sec:apriori} and energy-stable hyper-viscosity in Section ~\ref{sec:hyperviscosity}. Finally, we conclude with numerical experiments in one and two space dimensions (Sec.~\ref{sec:numexp}), verifying our theoretical analysis. Section~\ref{sec:conclusion} summarises the main results of the study.

\section{Continuous analysis}\label{sec:contanalysis}
Here we will analyse the SWEs in 1D, in vector invariant form. Both linear and nonlinear analysis will be performed. We then provide relevant definitions and derive well-posed linear and nonlinear BCs. Much of the arguments for the linear analysis can be deduced from standard texts \cite{gustafsson1995time}. The new results here are new well-posed BCs for the vector invariant form of the SWE formulated in terms of fluxes, and applicable to linear and nonlinear SWEs. We prove that the BCs are well-posed and energy/entropy stable.
%


%

\subsection{Rotating SWE in vector invariant form}
The vector invariant form of the nonlinear rotating SWE in 2D are:

\begin{equation}\label{eqn:nlSWE2D}
\begin{aligned}
\frac{\partial h}{\partial t } + \nabla \cdot (
    \mathbf{u}h) = G_h, \quad \frac{\partial \mathbf{u}}{\partial t} + \omega \mathbf{u}^\perp + \nabla \biggl(\frac{|\mathbf{u}|^2}{2} +gh\biggl) = \mathbf{G}_u\textcolor{blue}{,} \quad \omega = \frac{\partial v}{\partial x } - \frac{\partial u}{\partial y } + f_c, 
 \quad (x,y) \in \Omega \subset \mathbb{R}^2, \quad t\ge 0,
\end{aligned}
\end{equation}
where $(x,y)$ are the spatial variables, $t$ denotes time, $h>0$ and $\mathbf{u} = [u, v]^T$ are primitive variables denoting water height and the velocity vector, $\omega$ is the absolute vorticity, $f_c$ is the Coriolis frequency and $\mathbf{u}^\perp = [-v, u]^T$. The vector $\mathbf{G} =[{G}_h, \mathbf{G}_u]$ is a source term, and can be used to model non-flat bathymetry, with $\mathbf{G} =[0, -g\nabla b]$ where $b$ is the profile of the bathymetry.
In 1D, with $\mathbf{u} = [u, 0]^T$ and $\partial/\partial y =0$, \eqref{eqn:nlSWE2D} reduces to
\begin{equation}\label{eqn:nlSWE}
    \frac{\partial h}{\partial t } + \frac{\partial}{\partial x} (uh) = G_h, \quad \frac{\partial u}{\partial t} + \frac{\partial}{\partial x } \biggl(\frac{u^2}{2} +gh\biggl) = G_u\textcolor{blue}{,}
\quad x \in \Omega \subset \mathbb{R}, \quad t\ge 0,
\end{equation}
which can be written as
\begin{equation}\label{eqn:vecinvnl}
    \frac{\partial \mathbf{q}}{\partial t} + \frac{\partial \mathbf{F}}{\partial x} = \mathbf{G}\textcolor{blue}{,}
    \quad x \in \Omega \subset \mathbb{R}, \quad t\ge 0,
\end{equation}
where  $\mathbf{G} =[{G}_h, {G}_u]^T$, $\mathbf{q} = [h, u ]^T$, and $\mathbf{F} = [{F}_{1}, {F}_{2}]^T$  are the nonlinear fluxes with ${F}_1 = uh$ and ${F}_2 = u^2/2 + gh$.  We note that unlike the flux form \cite{lundgren2020efficient} where the integrals of the  conserved variables $[h, \mathbf{u}h]^T$ are invariants, when $\mathbf{G}=0$, here, in \eqref{eqn:vecinvnl}  the integrals of the vector of primitive variables $\mathbf{q} = [h, \mathbf{u}]^T$ are invariants, when $\mathbf{G}=0$.

We will also consider the linear version of the SWE \eqref{eqn:vecinvnl}. Introducing zero-mean quantities in the form of perturbed variables, i.e., $u = U + \widetilde{u}$ and $h = H +\widetilde{h}$ with the mean states $H > 0$, $U$, and discarding nonlinear terms of order $\mathcal{O}(\widetilde{u}^2,\widetilde{h}^2,\widetilde{u}\widetilde{h})$, we obtain the linear vector-invariant SWE \eqref{eqn:vecinvnl}
where {$\mathbf{F} = [{F}_{1}, {F}_{2}]^T$ are the linear fluxes,  with ${F}_{1} = Uh + Hu$ and ${F}_{2} = Uu+gh$, and we have dropped the tilde on fluctuating quantities for convenience.

\subsection{Linear analysis}\label{sec:linear}
As mentioned earlier, here, we consider first the  analysis of the linear SWE.
To begin we note that in the linear case, the flux $\mathbf{F}$ has a more simple expression, namely
\begin{equation}\label{eq:LinearFlux}
   \mathbf{F} = M \mathbf{q}, \quad  M = \begin{bmatrix}
U  & H \\
g  & U 
\end{bmatrix}. 
\end{equation}
To simplify the linear analysis, we will assume that the  mean states  are independent of the time variable $t$ but can depend on the spatial variable $x$, that is $U(x)$, $H(x)$. However, the analysis is valid when the mean states depend on both space and time, that is $U(x,t)$, $H(x,t)$, see Remark \ref{rem:W_t}.

\subsubsection{Well-posedness}\label{sec:linear-well-posedness}
We consider the linear  SWE \eqref{eqn:vecinvnl}--\eqref{eq:LinearFlux} with  a forcing function
\begin{equation}\label{eq:linear_SWE}
{\frac{\partial \mathbf{q}}{\partial t } = D\mathbf{q}  + \mathbf{G}, \quad \revi{D: = - \frac{\partial } {\partial x }(M \, \cdot)}},
\end{equation}
\revi{so that $D \mathbf{q} = - {\partial (M \, \mathbf{q})} /{\partial x }$ and $D \mathbf{q} = - M \,{\partial \mathbf{q}} /{\partial x } $ for a spatially constant coefficients matrix $M$}.
Further, we consider \eqref{eq:linear_SWE} in the bounded domain $\Omega = [0, L]$ with the boundary points $ \Gamma = \{0, L\}$,  and augment \eqref{eq:linear_SWE} with the  initial condition and homogeneous BCs
\begin{subequations}\label{eq:IC_BC}
\begin{equation}\label{eq:IC}
 \mathbf{q}(x,t =0) =  \mathbf{f}(x) \in L^2(\Omega), \quad x \in \Omega,
\end{equation}
\begin{equation}\label{eq:BC}
	\mathcal{B}\mathbf{q} = 0, \quad x \in \Gamma.
\end{equation}
\end{subequations}
Here, $\mathcal{B}$ is a linear boundary operator, that will be determined, with homogeneous boundary data. However, the analysis can be extended to non-homogeneous boundary data, but this would complicate the algebra.

In the following, we will introduce the relevant notation required to derive the boundary operator $\mathcal{B}$ and prove the well-posedness of the IBVP, \eqref{eq:linear_SWE}--\eqref{eq:IC_BC}. 
To begin, we define the weighted $L^2(\Omega)$ inner product and norm
\begin{equation}\label{eqn:L2norm}
 (\mathbf{p},\mathbf{q})_W = \int_{\Omega } \mathbf{p}^TW\mathbf{q} \; dx, \quad  \|\mathbf{q}\|^2_W = (\mathbf{q},\mathbf{q})_W,   
\end{equation}
where $W = W^T $ and $\mathbf{q}^TW\mathbf{q}>0 \quad \forall \mathbf{q}  \in \mathbb{R}^2 \backslash \{\mathbf{0}\}$ (i.e., symmetric and positive definite). 
For the linear SWE, we choose specifically the  symmetric weight matrix 
\begin{equation}\label{eqn:WeightMatrix}
W = \begin{bmatrix}
g & U \\
U & H
\end{bmatrix}. 
\end{equation}
The eigenvalues 
 of $W$ are $\lambda_W^{\pm} = \frac{1}{2}(g+H) \pm \sqrt{(g+H)^2 - 4gH(1-\textrm{Fr}^2})$, which are positive for subcritical flows, $\textrm{Fr}<1$, so that  $W = W^T $ and $\mathbf{q}^TW\mathbf{q}>0, \quad \forall \mathbf{q}  \in \mathbb{R}^2 \backslash \{\mathbf{0}\}$.
\begin{definition}\label{eqn:wellposed}
The IBVP \eqref{eq:linear_SWE}--\eqref{eq:IC_BC} is strongly well-posed if a unique solution, $\mathbf{q} = [h, u]^T$ satisfies $$\| \mathbf{q}\|_W \leq K e^{\mu t}\biggl(\| \mathbf{f}\|_W +  \max\limits_{\tau \in [0, t]}\|\mathbf{G(.~, \tau})\|_W\biggl),$$ for every forcing function $\mathbf{G} \in L^2(\Omega)$ and some constants $K>0$, $\mu$ $\in \mathbb{R}$ independent of $\mathbf{f}$ and  $\mathbf{G}$.
\end{definition}

The well-posedness of the IBVP \eqref{eq:linear_SWE}--\eqref{eq:IC_BC} can be related to the boundedness of the differential operator $D$. We introduce the function space 
\begin{equation}\label{eq:functionspace_V}
\mathbb{V}=\left\{\mathbf{p} \mid \, \mathbf{p}(x) \in \mathbb{R}^2, \quad\|\mathbf{p}\|_W<\infty, \quad x \in \Omega, \quad\{\mathcal{B} \mathbf{p}=0, x \in \Gamma\}\right\}.
\end{equation}
The following definition will be useful.
\begin{definition}\label{def:semibounded}
The differential operator $D$ is semi-bounded in the function space \revi{$\mathbb{V}$ if $\exists \, \mu \in \mathbb{R}$ independent of $\mathbf{q} \in \mathbb{V}$ such that
$$ 
(\mathbf{q}, D\mathbf{q})_W \leq \mu \|\mathbf{q}\|^2_W, \; \forall \mathbf{q} \in \mathbb{V}.
$$
}
\end{definition}
\begin{lemma}
Consider the linear differential operator \revi{ $D = -\frac{\partial}{\partial x}\left(M\,\cdot\right)$} given in \eqref{eq:linear_SWE} subject to the BCs, $\mathcal{B}\mathbf{q} = 0$, \eqref{eq:BC}, where $M$ is the constant coefficient matrix defined in \eqref{eq:LinearFlux}.
Let $W$ be the symmetric positive definite matrix given by  \eqref{eqn:WeightMatrix}, defining the weighted $L^2$-norm \eqref{eqn:L2norm}, then the  matrix product $\widetilde{M} = WM$ is symmetric. Further, if $\mathcal{B}\mathbf{q} = 0$ is such that     $\biggl(\mathbf{q}^T \widetilde{M} \mathbf{q}\biggl)\biggl\rvert^L_0 = \biggl(\mathbf{q}^T \widetilde{M} \mathbf{q}\biggl)\biggl\rvert_L- \biggl(\mathbf{q}^T \widetilde{M} \mathbf{q}\biggl)\biggl\rvert_0 \ge 0$, then $D$ is semi-bounded.
\end{lemma}
\begin{proof}
First, we consider the matrix product,
 $\widetilde{M} =WM$, 
 \begin{equation}\label{eqn:wmlin}
\widetilde{M} =~\begin{bmatrix}
2gU & U^2 + gH  \\
U^2 + gH & 2HU
\end{bmatrix}.
\end{equation} 
Obviously $\widetilde{M}$ is symmetric.
Now consider $(\mathbf{q}, D\mathbf{q})_{W  }$ and integrate by parts, we have 
$$(\mathbf{q}, D\mathbf{q})_{W  } =- \frac{1}{2}\int_\Omega \frac{\partial  }{\partial x}\left(\mathbf{q}^T \widetilde{M}\mathbf{q}\right) = -\frac{1}{2}\biggl(\mathbf{q}^T \widetilde{M} \mathbf{q}\biggl)\biggl\rvert^L_0 .$$
So for the boundary operator $\mathcal{B}\mathbf{q} = 0$, if $\biggl(\mathbf{q}^T \widetilde{M} \mathbf{q}\biggl)\biggl\rvert^L_0 \ge 0$ then 
$$(\mathbf{q}, D\mathbf{q})_{W  } = -\frac{1}{2}\biggl(\mathbf{q}^T \widetilde{M} \mathbf{q}\biggl)\biggl\rvert^L_0 \le 0.$$
 An upper bound is the case of $\mu=0$. Clearly for $(\mathbf{q}, D\mathbf{q})_{W  } =0$, and $D$ is semi-bounded by Definition \eqref{def:semibounded}.
\end{proof}
It is often possible to formulate BCs $\mathcal{B}\mathbf{q} =0$ such that $\biggl(\mathbf{q}^T \widetilde{M} \mathbf{q}\biggl)\biggl\rvert^L_0 \ge 0$. An immediate example is the case of periodic BCs, where $\biggl(\mathbf{q}^T \widetilde{M} \mathbf{q}\biggl)\biggl\rvert^L_0 = 0$. However, for non-periodic BCs, it is imperative that the boundary operator must not destroy the  existence of solutions. We will now introduce the definition of maximally semi-boundedness of $D$ which will ensure well-posedness of the IBVP \eqref{eq:linear_SWE}--\eqref{eq:IC_BC}.
\begin{definition}\label{eqn:maximal}
The differential operator $D$ defined in \eqref{eq:linear_SWE} is maximally semi-bounded if it is semi-bounded in the function space $\mathbb{V}$
but not semi-bounded in any function space with fewer BCs.
\end{definition}

The maximally semi-boundedness property  is intrinsically connected to well-posedness of the IBVP. We will formulate this result in the following theorem. The reader can consult \cite{gustafsson1995time} for more elaborate discussions.

\begin{theorem}\label{thm:swp}
    Consider the IBVP \eqref{eq:linear_SWE}--\eqref{eq:IC_BC} if the differential operator $D$ is maximally semi-bounded, $(\mathbf{q}, D\mathbf{q})_W \leq \mu \|\mathbf{q}\|^2_W$, then it is strongly well-posed. That is, there is a unique solution $\mathbf{q}$ satisfying the estimate
    $$\| \mathbf{q}\|_W \leq K e^{\mu t}\biggl(\| \mathbf{f}\|_W +  \max\limits_{\tau \in [0, t]}\|\mathbf{G(.~, \tau})\|_W\biggl), \quad {K = \max\left(1, {(1-e^{-\mu t})}/{\mu}\right)}.$$
\end{theorem}
\begin{proof}
    We consider 
\begin{align}
\begin{split}  
    \frac{1}{2} \frac{d}{dt} \| \mathbf{q}\|_W^2   = \left(\mathbf{q},\frac{\partial\mathbf{q}}{\partial t}\right)_W =\left(\mathbf{q},D\mathbf{q}\right)_W + {\left(\mathbf{q},\mathbf{G}\right)_W}.  
\end{split}
    \end{align}
    Semi-boundedness and Cauchy-Schwartz inequality yield
\begin{align}
\begin{split}  
    \frac{1}{2} \frac{d}{dt} \| \mathbf{q}\|_W^2   \leq \mu \| \mathbf{q} \|^2_W + \|\mathbf{q}\|_W \|\mathbf{G}\|_W \iff
    \frac{d}{dt} \| \mathbf{q} \|_W \leq  \mu \| \mathbf{q}\|_W  + \|\mathbf{G}\|_W. \\
\end{split}
    \end{align}
    
    Combining Gr\"onwall's Lemma and Duhammel's principle gives
    \begin{align}
    \|\mathbf{q}\|_W \leq e^{\mu t}\|\mathbf{f}\|_W + \int_0^t {e^{\mu \left(t-\tau\right)}} \|\mathbf{G(.~, \tau)}\|_W d\tau \leq  K e^{{\mu } t }\biggl(\|\mathbf{f}\|_W + \max_{\tau \in[0,t]}\|\mathbf{G}(.~,\tau)\|_W\biggl), 
\end{align}
where $ {K = \max\left(1, {(1-e^{-\mu t})}/{\mu}\right)}$. When $\mu = 0$, $K = \max\left(1, t\right)$ where $t>0$ is the final time, and we have the required result of well-posedness given by Definition \ref{eqn:wellposed}.
\end{proof}
\revi{
It is significantly important to note that while we have considered homogeneous boundary data here the analysis can be extended to non-homogeneous boundary data, in particular when $\mu < 0$. For more detailed discussions and  examples see  \cite{nordstrom2022linear,GHADER20141,NORDSTROM2020112857} and the references therein. Furthermore, numerical experiments performed later in this study confirms that our results extend to non-homogeneous boundary data. 
}
Theorem \ref{thm:swp} assumes that the weight matrix $W$ defined in \eqref{eqn:WeightMatrix} is time-independent. Nevertheless, the theorem also holds when $W$ is time-dependent. We formulate this as the remark.
\begin{remark}\label{rem:W_t}
   Note that if $W(x,t)$ then
\begin{align}
 \frac{1}{2} \frac{d}{dt} \| \mathbf{q}\|_W^2=\left(\mathbf{q},\frac{\partial\mathbf{q}}{\partial t}\right)_W + \frac{1}{2}\left(\mathbf{q},{W}_t\mathbf{q}\right) = \left(\mathbf{q},D\mathbf{q}\right)_W + {\left(\mathbf{q},\mathbf{G}\right)_W}+ \frac{1}{2}\left(\mathbf{q},{W}_t\mathbf{q}\right).
 \end{align} 
On the right hand side we use the fact that $D$ is maximally semi-bounded, and the last two terms can be bounded by Cauchy-Schwartz inequality giving
\begin{align}
\begin{split}  
    \frac{1}{2} \frac{d}{dt} \| \mathbf{q}\|_W^2   \leq (\mu+\alpha) \| \mathbf{q} \|^2_W + \|\mathbf{q}\|_W \|\mathbf{G}\|_W \iff
    \frac{d}{dt} \| \mathbf{q} \|_W \leq  (\mu+\alpha) \| \mathbf{q}\|_W  + \|\mathbf{G}\|_W,
\end{split}
    \end{align}
    where $0\le \alpha \le \max_{\tau \in [0,t]}(|W_t|/\min(2|\lambda_W^{-}|, 2|\lambda_W^{+}|))$ is a constant.
    Again, combining Gr\"onwall's Lemma and Duhammel's principle gives
\begin{align}
\|\mathbf{q}\|_W &\leq e^{\mu_0 t}\|\mathbf{f}\|_W + \int_0^t {e^{\mu_0 \left(t-\tau\right)}} \|\mathbf{G(.~, \tau)}\|_W d\tau \quad \leq  K e^{{\mu_0 } t }\biggl(\|\mathbf{f}\|_W + \max_{\tau \in[0,t]}\|\mathbf{G}(.~,\tau)\|_W\biggl),\notag\\  K &= \max\left(1, {(1-e^{-\mu_0 t})}/{\mu_0}\right),  \quad \mu_0 = \mu+\alpha. 
\end{align}

\end{remark}

    \subsubsection{Well-posed linear BCs}
    We will now formulate well-posed BCs for the IBVP \eqref{eq:linear_SWE}--\eqref{eq:IC_BC}.
Well-posed BCs require that the  differential operator $D$ to be  maximally semi-bounded in the function space $\mathbb{V}$. That is, we need a minimal number of BCs so that $D$ is semi-bounded in $\mathbb{V}$.  In general, the number of BCs must be equal to the number of nonzero eigenvalues of $\widetilde{M}$ and the location of the BCs would depend on the signs of the eigenvalues. In particular, the number of BCs at $x =0$ must be equal to the number of positive eigenvalues of $\widetilde{M}$ and  the number of BCs at $x =L$  equal  the number of negative eigenvalues of $\widetilde{M}$. The matrix $\widetilde{M}$ has two nonzero eigenvalues, namely 
$$
\begin{aligned}
\lambda^{+} = U(H + g) + \sqrt{U^2(H + g)^2 + (U^2 + gH)^2}, \\ \quad \lambda^{-} = U(H + g) - \sqrt{U^2(H + g)^2 + (U^2 + gH)^2}.
\end{aligned}
$$
The eigenvalues are real and have opposite signs, with  $\lambda^{+} >0$ and $\lambda^{-}<0$.
 This implies that, for subcritical flow, we  require one BC at $x= 0 $, and one BC at $x= L$.

To enable effective numerical treatments, for both the linear and nonlinear SWE, we will formulate the BC in terms of fluxes. First we note that
$
\frac{1}{2}\mathbf{q}^T \widetilde{M} \mathbf{q} = {F}_{1}{F}_{2},
$
where ${F}_{1} = Uh + Hu$ and ${F}_{2} = gh + Uu$ are linear mass flux and linear velocity flux, respectively.
The BC is given by
 \begin{align}\label{eq:BC-subcritical}
    \{\mathcal{B}\textbf{p}=0,  \ x\in \Gamma \} \equiv \left\{ { \alpha_1F_1 + \alpha_2F_2 =0, \ x =0};  \ {\beta_1F_1 - \beta_2F_2 =0, \ x = L}\right\},
    \end{align}
where $\alpha_j\ge 0$ and $\beta_j\ge 0$ (for $j = 1, 2$) are real constants that do not vanish, that is $|\alpha|>0$ and $|\beta|>0$. The BC \eqref{eq:BC-subcritical} can model different physical situations. For example 
\begin{enumerate}
    \item linear mass flux: $\alpha_1 > 0,\quad \alpha_2 =0$ \revi{;} \quad $\beta_1 > 0, \quad \beta_2 =0$,
    \item linear velocity flux (pressure BC): $\alpha_1 =0, \quad \alpha_2 >0$ \revi{;} \quad $\beta_1 =0, \beta_2 >0$,
    \item linear transmissive BC: $\alpha_1 =1, \quad \alpha_2 =\sqrt{H/g}$ \revi{;} \quad $\beta_1 = 1, \quad \beta_2 =\sqrt{H/g}$.
\end{enumerate}
The linear transmissive BC is equivalent to setting the incoming linear Riemann invariants to zero on the boundary, that is
\begin{align}\label{eq:BC-subcriticalTransmissive}
     \{\mathcal{B}\textbf{p}=0,  \ x\in \Gamma \} \equiv \left\{ { \sqrt{\frac{g}{H}}h + u =0, \ x =0};  \ {\sqrt{\frac{g}{H}}h - u = 0, \ x = L}\right\}.
     \end{align}
\begin{lemma}\label{lem:maximally_semi_bounded}
Consider the linear differential operator $D\mathbf{q} = -\frac{\partial }{\partial x}(M\mathbf{q})$ given in \eqref{eq:linear_SWE} subject to the BCs, $\mathcal{B}\mathbf{q} = 0$, \eqref{eq:BC-subcritical} with  $\alpha_j\ge 0$ and $\beta_j\ge 0$ (for $j = 1, 2$). For subcritical flows, with $\textrm{Fr} = |U|/\sqrt{gH} < 1$, the  differential operator 
 $D$ is maximally semi-bounded in $\mathbb{V}$.
\end{lemma}
\begin{proof}
As required, there is  one BC at $x= 0 $, and one BC at $x= L$.
It suffices to show that the differential operator $D$ is semi-bounded $(\mathbf{q}, D\mathbf{q})_{W  } \le 0$. Now consider $(\mathbf{q}, D\mathbf{q})_{W  } $ and integrate by parts, we have
$$
(\mathbf{q}, D\mathbf{q})_{W  }=-\frac{1}{2}\biggl(\mathbf{q}^T \widetilde{M} \mathbf{q}\biggl)\biggl\rvert^L_0= -{F}_{1}{F}_{2}\biggl\rvert^L_0 = {F}_{1}{F}_{2}\biggl\rvert_0 -{F}_{1}{F}_{2}\biggl\rvert_L.
$$
Note that if  $\alpha_1 =0, \quad \alpha_2 >0$ then $F_2 =0$ at $x =0$, and if $\beta_1 =0, \beta_2 >0$ then $F_2 =0$ at $x = L$, and we have $(\mathbf{q}, D\mathbf{q})_{W  } =0$. Furthermore, if  $\alpha_1 >0$ then $F_1 =-\alpha_2/\alpha_1 F_2$ at $x = 0$ and if  $\beta_1 >0$ then $F_1 =\beta_2/\beta_1 F_2$ at $x = L$, giving $(\mathbf{q}, D\mathbf{q})_{W  }= -\alpha_2/\alpha_1 F_2^2|_{x =0} -\beta_2/\beta_1 F_2^2|_{L} \leq0$. Thus for all $\alpha_j\ge 0$ and $\beta_j\ge 0$ we must have 
$$
(\mathbf{q}, D\mathbf{q})_{W  } = {F}_{1}{F}_{2}\biggl\rvert_0 -{F}_{1}{F}_{2}\biggl\rvert_L \le 0.
$$
As before, an upper bound is the case of $\mu=0$, which satisfies $(\mathbf{q}, D\mathbf{q})_{W  } =0$.
\end{proof}
 
We introduce the boundary term $\operatorname{BT}$ defined by 
\begin{align}\label{eq:BT_cont}
  \operatorname{BT} = {F}_{1}{F}_{2}\biggl\rvert_0 -{F}_{1}{F}_{2}\biggl\rvert_L \le 0.  
\end{align}
Note that by Lemma \ref{lem:maximally_semi_bounded} we must have $(\mathbf{q}, D\mathbf{q})_{W  } \le \mu \| \mathbf{q}\|_W^2\le 0$, where $\mu = \max_{\tau \in[0,t]}{\frac{\operatorname{BT}}{\| \mathbf{q}\|_{W}^2}} \le 0$.
\begin{theorem}\label{thm:energycontlin}
Consider the vector invariant form of the linear SWE \eqref{eq:linear_SWE}, at subcritical flows, with $\textrm{Fr} = |U|/\sqrt{gH} < 1$, subject to the initial condition \eqref{eq:IC} and the BC \eqref{eq:BC}, $\mathcal{B}\mathbf{q}=0$, where $\mathcal{B}\mathbf{q} $ is given by \eqref{eq:BC-subcritical} with  $\alpha_j\ge 0$ and $\beta_j\ge 0$ (for $j = 1, 2$) and $\mathrm{BT}\leq0$ given by \eqref{eq:BT_cont}. The corresponding IBVP \eqref{eq:linear_SWE}-\eqref{eq:IC_BC} is strongly well-posed. 
That is, there is a unique $\mathbf{q}$ satisfying the estimate,
$$
\begin{aligned}
\| \mathbf{q}\|_W \leq K e^{{\mu } t }\biggl(\|\mathbf{f}\|_W + \max_{\tau \in[0,t]}\|\mathbf{G}(.~,\tau)\|_W\biggl),
\quad
\mu = \max_{\tau \in[0,t]}{\frac{\operatorname{BT}}{\| \mathbf{q}\|_{W}^2}} \leq 0, \\  \quad {K = \max\left(1, {(1-e^{-\mu t})}/{\mu}\right)}.
\end{aligned}
$$

\end{theorem}
\begin{proof}
Invoking Lemma \ref{lem:maximally_semi_bounded} and  Theorem~\ref{thm:swp} yields the required result. 
\end{proof}
\subsection{Nonlinear analysis}\label{sec:nlanalysis}
We will now extend the linear analysis performed in section \ref{sec:linear} to the nonlinear vector invariant SWE. It is necessary that the nonlinear analysis must not contradict the conclusions drawn from the linear analysis. Otherwise, any valid linearisation will violate the linear analysis and would render the nonlinear theory ineffective. In particular, the number and location of BCs must be consistent with the linear theory.  The consistency between the linear and nonlinear theory is also necessary for proving the convergence of the numerical method for nonlinear problems, see Section \ref{sec:apriori} for details. 

\revi{
\begin{remark}
As noted in the introduction, remark that the consistency of linear theory and nonlinear theory advocated for here may contradict some of the results derived in \cite{nordstrom2022linear} where it was shown that the number of  BCs for the nonlinear SWEs  is independent of the magnitude of the flow, which disagrees with the linear theory \cite{GHADER20141} for critical and super-critical flows. Nevertheless, for sub-critical flows considered in this study the number of BCs is in agreement with the results of the nonlinear theory presented in \cite{nordstrom2022linear}. However, critical and super-critical flows are not considered in the current work.
\end{remark}
}
A main contribution of this study is the proof of strong well-posedness and strict stability of numerical approximations for the IBVP for nonlinear SWE in vector invariant form for subcritical flows, $\textrm{Fr}  < 1$. Many prior studies \cite{lundgren2020efficient} have proven linear stability, where the assumption of smooth solutions is then used to simulate nonlinear problems. There are a few exceptions, however, (see, e.g., \cite{nordstrom2024nonlinear, nordstrom2022linear}) where nonlinear analysis are considered for the IBVP. However, to the best of our knowledge there is no  numerical evidence yet verifying the theory.

The proof we present here is very similar to that given in Section~\ref{sec:linear}, so we shall keep it brief. To begin, we extend the linear BC  \eqref{eq:BC-subcritical} to the nonlinear regime.
The nonlinear BC is given by
 \begin{align}\label{eq:BC-subcritical_nonlinear}
    \{\mathcal{B}\textbf{p}=0,  \ x\in \Gamma \} \equiv \left\{ { \alpha_1F_1 + \alpha_2F_2 =0, \ x =0};  \ {\beta_1F_1 - \beta_2F_2 =0, \ x = L}\right\},
    \end{align}
where $F_1 = uh$ and $F_2 = u^2/2 + gh$ are nonlinear mass flux and velocity flux for \eqref{eqn:vecinvnl}, and $\alpha_j\ge 0$ and $\beta_j\ge 0$ (for $j = 1, 2$) are nonlinear coefficients that do not vanish.  As above the nonlinear BC \eqref{eq:BC-subcritical_nonlinear} can model different physical situations. We also give the nonlinear examples 
\begin{enumerate}
    \item Mass flux: $\alpha_1 > 0,\quad \alpha_2 =0$ \revi{;} $\beta_1 > 0, \quad \beta_2 =0$.
    \item Velocity flux (pressure BC): $\alpha_1 =0, \quad \alpha_2 >0$ \revi{;} $\beta_1 =0, \beta_2 >0$.
    \item Transmissive BC: $\alpha_1 =1, \quad \alpha_2 =\sqrt{h/g}(\sqrt{gh}-0.5u)/(\sqrt{gh}-u)>0$ \revi{;} \\  $\beta_1 = 1, \quad \beta_2 =\sqrt{h/g}(\sqrt{gh}+0.5u)/(\sqrt{gh}+u)>0$.
\end{enumerate}
%
%
With some tedious but straightforward algebraic manipulations, it can be shown that the nonlinear transmissive BC is equivalent to setting the incoming nonlinear Riemann invariants to zero on the boundary,
\begin{align}\label{eq:BC-subcriticalTransmissive_nonlinear}
     \{\mathcal{B}\textbf{p}=0,  \ x\in \Gamma \} \equiv \left\{ {  u + 2\sqrt{gh} =0, \ x =0};  \ { u - 2\sqrt{gh} = 0, \ x = L}\right\}.
     \end{align}
For subcritical flows, with  $\textrm{Fr}=|u|/\sqrt{gh}< 1$,  the coefficients $\alpha_2 =\sqrt{h/g}(\sqrt{gh}-0.5u)/(\sqrt{gh}-u)>0$ and $\beta_2 =\sqrt{h/g}(\sqrt{gh}+0.5u)/(\sqrt{gh}+u) >0$ for the nonlinear transmissive BCs are positive and depend nonlinearly on $u$ and $h>0$. This is  opposed to the linear case where $\alpha_2 =\sqrt{H/g}$, $\beta_2 =\sqrt{H/g}$, and $H>0$ are known constants.

We introduce the nonlinear differential operator
\begin{equation}\label{eqn:nlSWE_DiffOp}
\mathcal{D}\mathbf{q}: = 
-\begin{bmatrix}
    \frac{\partial}{\partial x} (uh)
    \\
    \frac{\partial}{\partial x } \biggl(\frac{u^2}{2} +gh\biggl)
\end{bmatrix}
=  -\frac{\partial \mathbf{F}}{\partial x }.
\end{equation}
We also define
\begin{equation}\label{eqn:L2norm_nln}
 (\mathbf{p},\mathbf{q})_W = \int_{\Omega } \mathbf{p}^TW\mathbf{q} \; dx, \quad  \|\mathbf{q}\|^2_W = (\mathbf{q},\mathbf{q})_W,   \quad W = \begin{bmatrix}
g & \frac{u}{2} \\
\frac{u}{2} & \frac{h}{2}
\end{bmatrix}
, \quad h>0.
\end{equation}
The eigenvalues 
 of $W$ are $\lambda_W^{\pm} = \frac{1}{2}(g+\frac{1}{2}h) \pm \frac{1}{2}\sqrt{(g+\frac{1}{2}h)^2 - {gh} (2-\textrm{Fr}^2) }$, which are positive for subcritical flows, $\textrm{Fr}<1<2$, so that  $W = W^T $ and $\mathbf{q}^TW\mathbf{q}>0 \quad \forall \mathbf{q}  \in \mathbb{R}^2 \backslash \{\mathbf{0}\}$.
 We will prove the nonlinear equivalence of Lemma \ref{lem:maximally_semi_bounded}.
\begin{lemma}\label{lem:maximally_semi_bounded_nln}
Consider the nonlinear differential operator $\mathcal{D}$ given in \eqref{eqn:nlSWE_DiffOp} subject to the nonlinear BC  \eqref{eq:BC-subcritical_nonlinear} with  $\alpha_j\ge 0$ and $\beta_j\ge 0$ (for $j = 1, 2$). For subcritical flows, with $\textrm{Fr} = |u|/\sqrt{gh} < 1$, the  differential operator 
 $\mathcal{D}$ is maximally semi-bounded in $\mathbb{V}$.
\end{lemma}
\begin{proof}
There is  one BC at $x= 0 $, and one BC at $x= L$ which agrees with the linear theory. So
it suffices to show that the differential operator $\mathcal{D}$ is semi-bounded $(\mathbf{q}, \mathcal{D}\mathbf{q})_{W  } \le 0$.
We consider $(\mathbf{q}, \mathcal{D}\mathbf{q})_{W  }$ and integrate by parts, we have
$$
(\mathbf{q}, \mathcal{D}\mathbf{q})_{W  }= -{F}_{1}{F}_{2}\biggl\rvert^L_0 = {F}_{1}{F}_{2}\biggl\rvert_0 -{F}_{1}{F}_{2}\biggl\rvert_L.
$$
We have to show that ${F}_{1}{F}_{2}\biggl\rvert_0 
\le 0$ and ${F}_{1}{F}_{2}\biggl\rvert_L \ge 0$,
which follows from the proof of Lemma \ref{lem:maximally_semi_bounded}. 
Therefore, for all $\alpha_j\ge 0$ and $\beta_j\ge 0$ we must have 
$
 (\mathbf{q}, \mathcal{D}\mathbf{q})_{W  }  \le 0.
$
\end{proof}

To prove the nonlinear equivalence of Theorem \ref{thm:energycontlin}, we introduce the nonlinear weighted $L^2$ norm
\begin{equation}\label{eqn:L2norm_nonlinear}
 \|\mathbf{q}\|^2_{W^{\prime}}  = \int_{\Omega } \mathbf{q}^T{W^{\prime}}\mathbf{q} dx,    \quad W^{\prime} = \begin{bmatrix}
g &0 \\
0 & {h}
\end{bmatrix}, \quad h>0.
\end{equation}
For subcritical flows, with $\textrm{Fr} = |u|/\sqrt{gh} < 1$, the two nonlinear  norms, $\|\cdot\|_{W}$ and $\|\cdot\|_{W^{\prime}}$,  are equivalent, that is  for positive real numbers $0<\gamma_1 \le 1$,  $\gamma_2 \ge 2$ we have
\begin{align}\label{eq:equivalence_l2_norm}
\gamma_1\|\cdot\|_{W^{\prime}}\le \|\cdot\|_{W} \le \gamma_2\|\cdot\|_{W^{\prime}}.
\end{align}
Note that for $\mathbf{q} = [h, u]^T$, the quantity 
\begin{equation}\label{eqn:entropy_nonlinear}
e=\frac{1}{2}\mathbf{q}^T W^{\prime}\mathbf{q}= \frac{1}{2} ( gh^2 +h u^2), 
\end{equation} 
is the elemental energy. For subcritical flows, with $\textrm{Fr} = |u|/\sqrt{gh} < 1$, the elemental energy $e$ is a convex function (in terms of the prognostic variables $u$ and $h$) and defines an entropy, see \cite{ricardo2023conservation}. When $\mathbf{G} =0$, the entropy/energy is conserved for smooth solutions, and dissipated across shocks and discontinuous solutions \cite{gustafsson1995time}.

Let $\mathrm{BT}$ be given by \eqref{eq:BT_cont}. Note that by Lemma \ref{lem:maximally_semi_bounded_nln} we must have $(\mathbf{q}, \mathcal{D}\mathbf{q})_{W  } \le \mu \| \mathbf{q}\|_{W^{\prime}}^2\le 0$, where $\mu = \max_{\tau \in[0,t]}{\frac{BT}{\| \mathbf{q}\|_{W^{\prime}}^2}} \le 0$.
We are now ready to state the theorem which proves the nonlinear equivalence of Theorem \ref{thm:energycontlin}.
\begin{theorem}\label{thm:nlstable}
Consider the vector invariant form of the nonlinear SWE \eqref{eqn:vecinvnl}, subcritical flow regime with $\textrm{Fr} = |u|/\sqrt{gh} < 1$, subject to the initial condition \eqref{eq:IC} and the nonlinear BC  \eqref{eq:BC-subcritical_nonlinear} with  $\alpha_j\ge 0$ and $\beta_j\ge 0$ (for $j = 1, 2$) and $\mathrm{BT}\le0$ be given by \eqref{eq:BT_cont}.  The solutions $\mathbf{q}$ of the corresponding IBVP \eqref{eqn:vecinvnl}, \eqref{eq:IC} and \eqref{eq:BC-subcritical_nonlinear}   satisfy the nonlinear energy estimate
    $$
\begin{aligned}
    \| \mathbf{q}\|_{W^{\prime}}\leq Ke^{\mu t}\biggl(\| \mathbf{f}\|_{W^{\prime}} +  \max\limits_{\tau \in [0, t]}\|\mathbf{G(.~, \tau})\|_{W^{\prime}}\biggl), \quad \mu = \max_{\tau \in[0,t]}{\frac{BT}{\| \mathbf{q}\|_{W^{\prime}}^2}} \le 0, \quad {K = \max\left(1, {(1-e^{-\mu t})}/{\mu}\right)}.
\end{aligned}
    $$
\end{theorem}
\begin{proof}
As before, from the left, we  multiply the vector invariant nonlinear SWE  \eqref{eqn:vecinvnl} with $ [\mathbf{F}_2, \ \mathbf{F}_1] = \mathbf{q}^TW$ and integrate over the domain $\Omega$, where here $W$ is defined by \eqref{eqn:L2norm_nln}.
This yields, 
\begin{align}
\begin{split}  
    \left(\mathbf{q},\frac{\partial\mathbf{q}}{\partial t}\right)_W =\left(\mathbf{q},\mathcal{D}\mathbf{q}\right)_W + {\left(\mathbf{q},\mathbf{G}\right)_W}.  
\end{split}
    \end{align}
We note that different from the linear case, we have, 
$$
 \mathbf{q}^T W \frac{\partial \mathbf{q}}{\partial t} = \frac{1}{2}\frac{\partial }{\partial t} \left( \mathbf{q}^T W^{\prime}\mathbf{q}\right) ,    \quad W^{\prime} = \begin{bmatrix}
g &0 \\
0 & {h}
\end{bmatrix}, \quad h>0.
$$
The semi-boundedness of $\mathcal{D}$ yields, 
\begin{equation}\label{eqn:BTnl}
    \frac{1}{2}\frac{d}{dt}\|\mathbf{q}\|_{W^{\prime}}^2  \le \mu \| \mathbf{q}\|_{W^{\prime}}^2 + (\mathbf{q}, \widetilde{\mathbf{G}})_{W^{\prime}}, \quad \widetilde{\mathbf{G}} = {(W^{\prime})}^{-1}WG,
\end{equation}
 where $\mu = \max_{\tau \in[0,t]}{\frac{BT}{\| \mathbf{q}\|_{W^{\prime}}^2}} \le 0$.
Note that at subcritical flow regime we have $|\widetilde{\mathbf{G}}| \le  |\mathbf{G}|$.
 We apply Cauchy-Schwarz inequality to the right hand side of \eqref{eqn:BTnl} giving
\begin{equation}\label{eqn:BTnl_0}
    \frac{1}{2}\frac{d}{dt}\|\mathbf{q}\|_{W^{\prime}}^2  \leq \mu \| \mathbf{q}\|_{W^{\prime}}^2 + \|\mathbf{q}\|_{W^{\prime}}\|\widetilde{\mathbf{G}}\|_{W^{\prime}}
    \leq \mu \| \mathbf{q}\|_{W^{\prime}}^2+ \|\mathbf{q}\|_{W^{\prime}}\|{\mathbf{G}}\|_{W^{\prime}}.
\end{equation}
Making use of Gr\"onwall's Lemma and Duhammel's principle gives
    \begin{align}
    \| \mathbf{q}\|_{W^{\prime}}\leq Ke^{\mu t}\biggl(\| \mathbf{f}\|_{W^{\prime}} +  \max\limits_{\tau \in [0, t]}\|\mathbf{G(.~, \tau})\|_{W^{\prime}}\biggl). 
\end{align}

\end{proof}
In the coming sections we will introduce numerical approximations for the IBVPs using DP-SBP operators. We will prove both linear and nonlinear numerical stability by deriving discrete energy estimates analogous to Theorems \ref{thm:energycontlin}--\ref{thm:nlstable}.
\section{Dual-pairing SBP  FD operators}\label{sec:sbp}
In order to introduce the relevant notation and keep our study self-contained, we will firstly give a quick introduction to the standard DP-SBP FD operators derived by \cite{mattsson2017diagonal}, as an extension to the traditional (centred difference) SBP (see e.g., reviews by \cite{svard2014review,fernandez2014review}). Then we discuss the recently derived dispersion-relation preserving DP operators of \cite{williams2021provably}.

 Consider a finite 1D interval, $\Omega = [0, L]$, and discretise it into $N+1$ grid points $x_j$ with contant spatial step $\Delta{x} >0$, we have
\begin{equation}
     x_j = j\Delta x, \quad j \in \{ 0, 1, \hdots ,N \}, \quad \Delta x = \frac{L}{N}, \quad \mathbf{x}:=(x_0, x_1,\ldots, x_{N})^T \in \mathbb{R}^{N+1}.
\end{equation}
We will use DP upwind, $\alpha$-DRP-DP upwind and traditional SBP operators \cite{mattsson2017diagonal,duru2022dual,williams2021provably,svard2014review} to approximate the spatial derivatives, $\partial/{\partial x}$.
Let 
\begin{equation}\label{eq:SBP-norm}
    P:=\operatorname{diag}\left(p_0, p_1, \ldots, p_N\right), \quad  p_j>0, \quad \forall j
\end{equation}
be the weights of a composite quadrature rule such that $\sum_{j = 0}^N p_j \zeta(x_j)   \approx \int_0^L \zeta(x) dx$ for an integrable function $\zeta(x)$. Then, we introduce the DP first derivative operators \cite{williams2021provably,mattsson2017diagonal} $D_{\pm}: \mathbb{R}^{n+1} \mapsto \mathbb{R}^{n+1}$ with $(D_{\pm}\boldsymbol{f})_j \approx \partial f/\partial x|_{x=x_j}$ so that the SBP property holds:
\begin{equation}\label{sbp_property}
\boldsymbol{g}^T P\left(D_{+} \boldsymbol{f}\right) +\boldsymbol{f}^T P \left(D_{-} \boldsymbol{g}\right)=f\left(x_{N}\right) g\left(x_{N}\right)-f\left(x_0\right) g\left(x_0\right),
\end{equation}
where $\boldsymbol{f}=\left(f\left(x_0\right), \ldots, f\left(x_{N}\right)\right)^T$, $\boldsymbol{g}=\left(g\left(x_0\right), \ldots, g\left(x_N\right)\right)^T$ are vectors sampled from differentiable functions, $f, g \in C^1([0,L])$. 
Now we introduce the matrix operators,
\begin{equation}\label{eq:linearoperators}
    \begin{gathered}
Q_{\pm}=P D_{\pm}, \quad A_{\pm}=Q_\pm - \frac{1}{2} B, \quad B=e_N e_N^T-e_\revi{0}e_\revi{0}^T, \\ \quad \revi{e_0}:=(1,0, \ldots, 0)^T, \quad e_N:=(0, \ldots, 0,1)^T \in \mathbb{R}^{N+1},
\end{gathered}
\end{equation}
so that from \eqref{sbp_property} we have $ Q_+ + Q_{-}^T =B$ and $ A_+ + A_{-}^T =0.$
\begin{definition}
Let $D_{\pm}, Q_{\pm}, A_{\pm}: \mathbb{R}^{N+1} \mapsto \mathbb{R}^{N+1}$ be linear operators that solve \eqref{sbp_property} and \eqref{eq:linearoperators} and for the diagonal norm $P \in \mathbb{R}^{n \times n}$. If the matrix $S_{+}:=A_{+}+A_{+}^T$ is negative semi-definite and $S_{-}:=A_{-}+A_{-}^T$ is positive semi-definite, then the 3-tuple $\left(P, D_{-}, D_{+}\right)$ is called an upwind diagonal-norm DP-SBP operator.
We call $\left(P, D_{-}, D_{+}\right)$ an upwind diagonal-norm DP-SBP operator of order $q$ if the accuracy conditions
$$
\left(D_{\pm}\mathbf{x}^i\right)_j=i {x}_j^{i-1} , i \in \{0,1 \hdots q\},
$$
is satisfied within the interior points $x_j$, $n_s < j < N-n_s+1$ for some $n_s \ge 1$. The FD stencils for boundary points $x_j$, for $0\le j\le n_s$ and $N-n_s+1 \le j \le N$, satisfy the same property but have either $q/2$ order accurate stencils if $q$ is even, or $(q-1)/2$ if $q$ is odd.
\end{definition}
\section{Semi-discrete approximation}\label{sec:semiapprox}
Now we consider the semi-discrete vector-invariant SWE, and  discuss numerical  stability and convergence  properties of the semi-discrete approximation.
To begin, we consider the semi-discrete approximation of the SWE, using the DP-SBP framework \cite{williams2021provably}
\begin{equation}\label{eq:semi-discrete-SWE}
   \frac{d\mathbf{q}}{dt} + \mathbf{D}_x \mathbf{F} = \mathbf{G} + \mathbf{SAT}, \quad \mathbf{q}(0) =  \mathbf{f}, \quad \mathbf{D}_x = \begin{bmatrix}
{D}_{+}  & \mathbf{0} \\
\mathbf{0}  & {D}_{-} 
\end{bmatrix}, \quad \mathbf{SAT} = \begin{bmatrix}
\textrm{SAT}_1  \\
\textrm{SAT}_2 
\end{bmatrix},
\end{equation}
 where $\mathbf{q} = [\mathbf{h} \quad \mathbf{u} ]^T$ are grid functions, and $\mathbf{F} = [\mathbf{F}_{1} \quad \mathbf{F}_{2}]^T$  are  fluxes, with $\mathbf{F}_1 = \mathbf{u}\mathbf{h}$ and $\mathbf{F}_2 = \mathbf{u}^2/2 + g\mathbf{h}$ for the nonlinear SWE and $\mathbf{F}_{1} = U\mathbf{h} + H\mathbf{u}$ and $\mathbf{F}_{2} = U\mathbf{u}+g\mathbf{h}$ for the linear SWE. Here $\mathbf{SAT}$ are penalty terms defined by 
 \revi{
\begin{equation}\label{eq:SAT-term}
\begin{aligned}
    \textrm{SAT}_1 = -\tau_{11}P^{-1}\mathbf{e}_1\left(\alpha_1F_{10} + \alpha_2F_{20}-g_0(t)\right) - \tau_{12}P^{-1}\mathbf{e}_N\left(\beta_1F_{1N} - \beta_2F_{2N}-g_L(t)\right), \\
    \textrm{SAT}_2 = -\tau_{21}P^{-1}\mathbf{e}_1\left(\alpha_1F_{10} + \alpha_2F_{20}-g_0(t)\right) - \tau_{22}P^{-1}\mathbf{e}_N\left(\beta_1F_{1N} - \beta_2F_{2N}-g_L(t)\right),
\end{aligned}
\end{equation}
weakly implementing the  linear BC \eqref{eq:BC-subcritical} or  nonlinear BC \eqref{eq:BC-subcritical_nonlinear} modified by non-homogeneous boundary. Here $g_0(t)$ and $g_L(t)$ are the non-homogeneous boundary data at $x =0$ and $x=L$. When $g_0(t)=0$ and $g_L(t)=0$ we recover the homogeneous BCs BC \eqref{eq:BC-subcritical} or  nonlinear BC \eqref{eq:BC-subcritical_nonlinear}.
}
The real coefficients $\tau_{ij}$ are penalty parameters to be determined by requiring stability. Note that the spatial derivative in the continuity equation is approximated with $D_+$ and the spatial derivative in the momentum equation is approximated with the dual-pair $D_-$. \revi{Note that using the forward and backward difference operators in this order is a particular choice. There is much more flexibility with the design of numerical approximations using DP SBP operators. For example the theory and results also hold if, when we reverse the order, the continuity equation is approximated with $D_-$ and the spatial derivative in the momentum equation is approximated with the dual-pair $D_+$.  This has been separately verified numerically in different tests}.

We will now analyse the semi-discrete approximation \eqref{eq:semi-discrete-SWE}. \revi{As in the continuous setting, we will consider homogeneous boundary data $g_0(t)=0$ and $g_L(t)=0$,  begin with the linear analysis and then proceed to the nonlinear analysis.}

\subsection{Linear semi-discrete analysis}\label{sec:linearsat}

We will prove the stability of the semi-discrete numerical approximation \eqref{eq:semi-discrete-SWE} for linear fluxes $\mathbf{F}_{1} = U\mathbf{h} + H\mathbf{u}$ and $\mathbf{F}_{2} = U\mathbf{u}+g\mathbf{h}$.
Now, we introduce the discrete weighted $L_2$-norm,
\begin{equation}\label{def:wpnorm}
\| {\textbf{q}} \|_{WP}^2:= {\textbf{q}}^T\left(W\otimes P\right){\textbf{q}}  > 0, \quad \forall  \textbf{q}\ne \textbf{0}.
\end{equation}
where the linear weight matrix $W$ is given by \eqref{eqn:WeightMatrix} and $P$ is the diagonal SBP norm defined in \eqref{eq:SBP-norm}. The main idea is to choose penalty parameters $\tau_{ij}$ such that we can prove a discrete analogue of Theorem \ref{thm:energycontlin}. To be precise, we introduce the following definition.
\begin{definition}\label{eqn:strong_stabilty}
The semi-discrete approximation  \eqref{eq:semi-discrete-SWE} is strongly stable  if the solution $\mathbf{q}$ satisfies $$\| \mathbf{q}\|_{WP} \leq K e^{\mu t}\biggl(\| \mathbf{f}\|_{WP} +  \max\limits_{t \in [0, \tau]}\|\mathbf{G(\tau})\|_{WP}\biggl),$$ for some constants  $K>0$, $\mu$ $\in \mathbb{R}$.
\end{definition}

We introduce the boundary term $\textrm{BT}$ given by
\begin{equation}\label{eq:BT_linear}
\begin{aligned}
   \textrm{BT} & = 
    {F}_{10}{F}_{20} - {F}_{1N}{F}_{2N} -\tau_{11}F_{20}\left(\alpha_1F_{10} + \alpha_2F_{20}\right) - \tau_{12}F_{2N}\left(\beta_1F_{1N} - \beta_2F_{2N}\right)\\
   &-\tau_{21}F_{10}\left(\alpha_1F_{10} + \alpha_2F_{20}\right) - \tau_{22}F_{1N}\left(\beta_1F_{1N} - \beta_2F_{2N}\right),
\end{aligned}
\end{equation}
\revi{with homogeneous boundary data $g_0(t)=0$ and $g_L(t)=0$.}
As we will see below,  the numerical approximation  \eqref{eq:semi-discrete-SWE} can be shown to be stable if $\textrm{BT} \le 0$.  
We can prove the following Lemma
\begin{lemma} \label{lem:stable_penalty}
    Consider the boundary term $\textrm{BT}$ given by
\eqref{eq:BT_linear}. If the penalty parameters $\tau_{ij}$ are chosen such that

$
\tau_{11} \ge 0, \quad \tau_{21} = {1}/{\alpha_2}, 
$
$
\tau_{12} \le 0, \quad \tau_{22} = {1}/{\beta_2},  
$
for $\alpha_1 = 0$, $\alpha_2 >0$, $ \beta_1 = 0$, $\beta_2 >0$, and 

$
\tau_{11} = {1}/{\alpha_1}, \quad \tau_{21} = 0, 
$
$
\tau_{12} = - {1}/{\beta_1}, \quad \tau_{22} = 0,  
$
for $\alpha_1 > 0$, $\alpha_2 \ge 0$, $ \beta_1 > 0$, $\beta_2 \ge 0$,

then $\textrm{BT} \le 0$.
\end{lemma}
\begin{proof}
   If 
$
\tau_{11} \ge 0, \quad \tau_{21} = {1}/{\alpha_2}, 
$
$
\tau_{12} \le 0, \quad \tau_{22} = {1}/{\beta_2},  
$
for $\alpha_1 = 0$, $\alpha_2 >0$, $ \beta_1 = 0$, $\beta_2 >0$, then we have
$$
\textrm{BT}  = 
    -\tau_{11} \alpha_2F_{20}^2 + \tau_{12}\beta_2F_{2N}^2\le 0.
$$
If $
\tau_{11} = {1}/{\alpha_1}, \quad \tau_{21} = 0, 
$
$
\tau_{12} =  -{1}/{\beta_1}, \quad \tau_{22} = 0,  
$
for $\alpha_1 > 0$, $\alpha_2 \ge 0$, $ \beta_1 > 0$, $\beta_2 \ge 0$, then we also have
$$
\textrm{BT}  = 
     -\frac{\alpha_2}{\alpha_1}F_{20}^2 -  \frac{\beta_2}{\beta_1}F_{2N}^2 \le 0.
$$
\end{proof}
The penalty parameters $\tau_{ij}$ ($i= 1, 2$, $j= 1, 2$) given by Lemma \ref{lem:stable_penalty} cover all well-posed BC parameters $\alpha_j \ge 0$, $\beta_j \ge 0$. However, the penalty parameters $\tau_{ij}$  ($i= 1, 2$, $j= 1, 2$) given above are not exhaustive.
We will now state the theorem which proves the stability of the semi-discrete approximation \eqref{eq:semi-discrete-SWE} for linear fluxes. 
\begin{theorem}\label{thm:st_well_posed_semi}
Consider the semi-discrete approximation \eqref{eq:semi-discrete-SWE} with $\mathbf{F}_{1} = U\mathbf{h} + H\mathbf{u}$ and $\mathbf{F}_{2} = U\mathbf{u}+g\mathbf{h}$, at subcritical flows, with $\textrm{Fr} = |U|/\sqrt{gH} < 1$ \revi{and  homogeneous boundary data $g_0(t)=0$, $g_L(t)=0$}. If the penalty parameters  $\tau_{ij}$ are chosen as in Lemma \ref{lem:stable_penalty} such that $\textrm{BT} \leq 0$, where $\textrm{BT}$ is given by \eqref{eq:BT_linear}, then  the semi-discrete approximation \eqref{eq:semi-discrete-SWE}  is strongly stable. 
That is, the numerical solution $\mathbf{q}$ satisfies the estimate
$$
\begin{aligned}
\| \mathbf{q}\|_{WP} \leq K e^{\mu t}\biggl(\| \mathbf{f}\|_{WP} +  \max\limits_{\tau \in [0, t]}\|\mathbf{G(\tau})\|_{WP}\biggl), \quad \mu = \max_{\tau \in[0,t]}{\frac{BT}{\| \mathbf{q}\|_{WP}^2}} \leq 0,  \quad {K = \max\left(1, {(1-e^{-\mu t})}/{\mu}\right)}.
\end{aligned}
$$
\end{theorem}
\begin{proof}
 Similar to the continuous setting, we consider 
\begin{align}
\begin{split}  
    \frac{1}{2} \frac{d}{dt} \| \mathbf{q}\|_{WP}^2   &= \left(\mathbf{q},\frac{ d \mathbf{q}}{d t}\right)_{WP} =-\left(\mathbf{q},\mathbf{D}_x\mathbf{F}\right)_{WP} + {\left(\mathbf{q},\mathbf{G} + \mathbf{SAT}\right)_{WP}}\\
    &= \underbrace{-\left(\mathbf{F}_{2}^T P \left(D_{+} \mathbf{F}_{1}\right) + \mathbf{F}_{1}^T P \left(D_{-} \mathbf{F}_{2}\right)\right)}_{\text{by DP-SBP \eqref{sbp_property}: } {F}_{10}{F}_{20} - {F}_{1N}{F}_{2N}} + {\left(\mathbf{q},\mathbf{G} + \mathbf{SAT}\right)_{WP}}.
\end{split}
    \end{align}
    Using the DP-SBP property \eqref{sbp_property} and the definition \eqref{eq:SAT-term} for $\mathbf{SAT}$ gives
    \begin{align}
\begin{split}  
    \frac{1}{2} \frac{d}{dt} \| \mathbf{q}\|_{WP}^2   
    = \textrm{BT} + {\left(\mathbf{q},\mathbf{G} \right)_{WP}} ,
\end{split}
    \end{align}
    where $\mathrm{BT}$ is given by \eqref{eq:BT_linear}.
     Cauchy-Schwartz inequality yields
\begin{align}
\begin{split}  
    \frac{1}{2} \frac{d}{dt} \| \mathbf{q}\|_{WP}^2   \leq \mu  \| \mathbf{q}\|_{WP}^2+ \|\mathbf{q}\|_W \|\mathbf{G}\|_{WP} \iff
    \frac{d}{dt} \| \mathbf{q} \|_{WP} \leq  \mu  \| \mathbf{q}\|_{WP}+ \|\mathbf{G}\|_{WP},
\end{split}
    \end{align}
    with
    $
    \mu = \max_{\tau \in[0,t]}{\frac{BT}{\| \mathbf{q}\|_{WP}^2}} \leq 0.
    $
   Again, Gr\"onwall's Lemma and Duhammel's principle give
$$
\| \mathbf{q}\|_{WP} \leq K e^{\mu t}\biggl(\| \mathbf{f}\|_{WP} +  \max\limits_{\tau \in [0, t]}\|\mathbf{G(\tau})\|_{WP}\biggl).
$$
\end{proof}
\subsection{Nonlinear semi-discrete analysis}
We will now consider the nonlinear fluxes, with $\mathbf{F}_1 = \mathbf{u}\mathbf{h}$ and $\mathbf{F}_2 = \mathbf{u}^2/2 + g\mathbf{h}$, and prove nonlinear stability for the semi-discrete numerical approximation \eqref{eq:semi-discrete-SWE}. We will follow closely  the steps taken in the last section to prove linear stability. To begin, we note that  the penalty parameters $\tau_{ij}$ given by Lemma \ref{lem:stable_penalty} are applicable to the nonlinear BCs \eqref{eq:BC-subcritical_nonlinear}. However, in this case the well-posed boundary parameters $\alpha_1 \ge 0$, $\beta_1 \ge 0$, $\alpha_2 \ge 0$, $\beta_2 \ge 0$ as well as  the penalty parameters $\tau_{ij}$  depend nonlinearly on the solution.

As in the continuous setting we will consider the nonlinear weight matrices $W$ and $W^{\prime}$, given by \eqref{eqn:L2norm_nln}--\eqref{eqn:L2norm_nonlinear},   projected on the grid yielding
\begin{equation}\label{def:wpnorm_weight}
\quad \mathbf{W}^{\prime} = \begin{bmatrix}
g\mathbf{I} & \textbf{0} \\
\textbf{0} & \mathrm{diag}(\mathbf{h})
\end{bmatrix},
\quad 
\mathbf{W} = \begin{bmatrix}
g\mathbf{I} & \frac{1}{2}\mathrm{diag}(\mathbf{u}) \\
\frac{1}{2}\mathrm{diag}(\mathbf{u})  & \frac{1}{2}\mathrm{diag}(\mathbf{h}) 
\end{bmatrix}
, \quad \mathbf{h}>0.
\end{equation}
Now, we introduce the discrete nonlinearly-weighted $L_2$-norm,
\begin{equation}\label{def:wpnorm}
\| {\textbf{q}} \|_{W^{\prime}P}^2:= {\textbf{q}}^T\mathbf{W}^{\prime}\left(I_2\otimes P\right){\textbf{q}}  > 0, \quad \forall  \textbf{q}\ne \textbf{0},
\end{equation}
 as well as the definition for nonlinear stability.
\begin{definition}\label{eqn:nl_strong_stabilty}
The semi-discrete approximation \eqref{eq:semi-discrete-SWE}, with nonlinear fluxes $\mathbf{F}_1 = \mathbf{u}\mathbf{h}$ and $\mathbf{F}_2 = \mathbf{u}^2/2 + g\mathbf{h}$, is strongly stable  if the numerical solution $\mathbf{q}$ satisfies 
$$
\| \mathbf{q}\|_{W^{\prime}P} \leq K e^{\mu t}\biggl(\| \mathbf{f}\|_{W^{\prime}P} +  \max\limits_{t \in [0, \tau]}\|\mathbf{G(\tau})\|_{W^{\prime}P}\biggl), 
$$ 
for some constants  $K>0$, $\mu$ $\in \mathbb{R}$.
\end{definition}
The following theorem  proves the nonlinear equivalence of Theorem \ref{thm:st_well_posed_semi}.
\begin{theorem}\label{thm:nl_semi_stable}
Consider the semi-discrete approximation \eqref{eq:semi-discrete-SWE} with $\mathbf{F}_1 = \mathbf{u}\mathbf{h}$ and $\mathbf{F}_2 = \mathbf{u}^2/2 + g\mathbf{h}$, at subcritical flows with $\textrm{Fr} = |u|/\sqrt{gh} < 1$. If the penalty parameters  $\tau_{ij}$ are chosen as in Lemma \ref{lem:stable_penalty} such that $\textrm{BT} \leq 0$, where $\textrm{BT}$ is given by \eqref{eq:BT_linear} \revi{with  homogeneous boundary data $g_0(t)=0$, $g_L(t)=0$}, then  the semi-discrete approximation \eqref{eq:semi-discrete-SWE}  is strongly stable. 
That is, the numerical solution $\mathbf{q}$ satisfies the estimate
    $$
\begin{aligned}
    \| \mathbf{q}\|_{W^{\prime}P}\leq Ke^{\mu t}\biggl(\| \mathbf{f}\|_{W^{\prime}P} +  \max\limits_{\tau \in [0, t]}\|\mathbf{G(\tau})\|_{W^{\prime}P}\biggl), \quad  \mu = \max_{\tau \in[0,t]}{\frac{BT}{\| \mathbf{q}\|_{WP}^2}} \leq 0, \quad {K = \max\left(1, {(1-e^{-\mu t})}/{\mu}\right)}.
\end{aligned}
    $$
\end{theorem}
\begin{proof}
This proof follows analogously from the continuous setting (Theorem~\ref{thm:nlstable}) once again, so we shall keep it brief. Multiplying \eqref{eq:semi-discrete-SWE} with $\mathbf{q}^T\mathbf{W}\left(I_2\otimes P\right)= [\mathbf{F}_2, \mathbf{F}_1]^T \left(I_2\otimes P\right)$ from the left yields, 
\begin{align}
\begin{split}  
    \left(\mathbf{q},\frac{d \mathbf{q}}{d t}\right)_{WP} &=-\left(\mathbf{q},\mathbf{D}_x\mathbf{F}\right)_{WP} + {\left(\mathbf{q},\mathbf{G} + \mathbf{SAT}\right)_{WP}}\\
    &= \underbrace{-\left(\mathbf{F}_{2}^T P \left(D_{+} \mathbf{F}_{1}\right) + \mathbf{F}_{1}^T P \left(D_{-} \mathbf{F}_{2}\right)\right)}_{\text{by DP-SBP \eqref{sbp_property}: } {F}_{10}{F}_{20} - {F}_{1N}{F}_{2N}} + {\left(\mathbf{q},\mathbf{G} + \mathbf{SAT}\right)_{WP}}.
\end{split}
    \end{align}
   Similar to the continuous analysis,  we also note that as opposed to the linear case, we have, 
$$
 \mathbf{q}^T \mathbf{W} \frac{d \mathbf{q}}{d t} = \frac{1}{2}\frac{d }{d t} \left( \mathbf{q}^T \mathbf{W}^{\prime}\mathbf{q}\right) ,    \quad \mathbf{W}^{\prime} = \begin{bmatrix}
g\mathbf{I} &0 \\
0 & \mathrm{diag}(\mathbf{h})
\end{bmatrix}, \quad \mathbf{h}>0.
$$
Using the DP-SBP property \eqref{sbp_property} and the definition \eqref{eq:SAT-term} for $\mathbf{SAT}$ gives
    \begin{align}
\begin{split}  
      \frac{1}{2}\frac{d }{dt} \|\mathbf{q}\|_{W^{\prime}P}^2  
    = \textrm{BT} + {\left(\mathbf{q},\mathbf{G} \right)_{WP}} \le \mu \|\mathbf{q}\|_{W^{\prime}P}^2 + {\left(\mathbf{q},\mathbf{G} \right)_{WP}} = \mu \|\mathbf{q}\|_{W^{\prime}P}^2 + (\mathbf{q},\widetilde{\mathbf{G}})_{W^{\prime}P}
\end{split}
    \end{align}
    where 
    $
    \mu = \max_{\tau \in[0,t]}{\frac{BT}{\| \mathbf{q}\|_{W^{\prime}P}^2}} \leq 0,
    $ $\widetilde{\mathbf{G}} = {(\mathbf{W}^{\prime})}^{-1}\mathbf{W}\mathbf{G}$ and at subcritical flow regime we have $|\widetilde{\mathbf{G}}| \le  |\mathbf{G}|$.

    Cauchy-Schwartz inequality yields
\begin{align}
\begin{split}  
    \frac{1}{2} \frac{d}{dt} \| \mathbf{q}\|_{W^{\prime}P}^2   \leq \mu \|\mathbf{q}\|_{W^{\prime}P}^2 + \|\mathbf{q}\|_{W^{\prime}} \|\mathbf{G}\|_{W^{\prime}P} \iff
    \frac{d}{dt} \| \mathbf{q} \|_{W^{\prime}P} \leq \mu \|\mathbf{q}\|_{W^{\prime}P} +  \|\mathbf{G}\|_{W^{\prime}P}. \\
\end{split}
    \end{align}
    Combining Gr\"onwall's Lemma and Duhamel's principle gives,
$$
\| \mathbf{q}\|_{W^{\prime}P} \leq Ke^{\mu t}\biggl(\| \mathbf{f}\|_{W^{\prime}P} +  \max\limits_{\tau \in [0, t]}\|\mathbf{G(\tau})\|_{W^{\prime}P}\biggl).
$$
\end{proof}
\section{A priori numerical error analysis}\label{sec:apriori}
We will now derive a priori error estimate for the semi-discrete approximation \eqref{eq:semi-discrete-SWE} with linear fluxes and nonlinear solution assuming smooth solutions.
Consider $\mathbf{q}_e(x_j,t)$ as the exact solution to the IBVP at $x_j = j \Delta x$, and $\mathbf{q}_j(t)$ a numerical solution obtained on the grid. Let $\mathbf{\mathcal{E}}_j(t) = [\mathcal{E}_{jh}(t), \mathcal{E}_{ju}(t)]^T:= \mathbf{q}_j(t)-\mathbf{q}_e(x_j,t)$ denote the error vector. Then, the error of  semi-discrete approximation \eqref{eq:semi-discrete-SWE} for the linear SWE~\eqref{eq:linear_SWE} satisfies the evolution equation:
\begin{equation}\label{eqn:error_eqn}
   \frac{d\mathbf{\mathcal{E}}}{dt} + \mathbf{D}_x 
    \mathcal{\mathbf{\mathbf{F}(\mathbf{\mathcal{E}})}} =  \textbf{\textrm{SAT}}(\mathbf{\mathcal{E}}) + \mathbb{T}, \quad \mathcal{E}(0) = \mathbf{0}, \quad \mathbf{D}_x = \begin{bmatrix}
\mathbf{D}_{+}  & \mathbf{0} \\
\mathbf{0}  & \mathbf{D}_{-} 
\end{bmatrix},
\end{equation}
where $\mathbf{F}_{1} = U\mathcal{E}_h + H\mathcal{E}_u$ and $\mathbf{F}_{2} = U\mathcal{E}_u+g\mathcal{E}_h$ are linear fluxes with corresponding \textbf{SAT} vectors. Here, $\mathbb{T}$ denotes the truncation error of the FD operators, which is: 
\begin{equation}\label{eq:truncation_error}
    \mathbb{T}_{j}=\left\{\begin{array}{ll}
\left.\Delta{x}^\gamma \chi_j \frac{\partial^{\gamma+1} \mathbf{q}_e}{\partial x^{\gamma+1}}\right|_{x_j}, & \text { if boundary, } \\
\left.\Delta{x}^\zeta \chi_j \frac{\partial^{\zeta+1} \mathbf{q}_e}{\partial x^{\zeta+1}}\right|_{x_j}, & \text { if interior, }
\end{array}\right.
\end{equation}
where $\chi_j$ are grid independent constants, $\Delta{x}$ is the grid spacing, $\gamma$ denotes the order of accuracy of the boundary stencils, and $\zeta$ is the order of accuracy of the interior stencils. 
For the traditional SBP schemes based on central difference stencils for the interior, the order $(\gamma, \zeta) = (p, 2p)$, $p \in \{1, 2, \cdots \}$ \cite{kreiss1974finite,gustafsson1995time}. The order of accuracy of the interior stencil of traditional SBP operators  is always even. For DP and DRP SBP operators satisfying  diagonal norm property but utilising upwind stencils can support even and odd order interior stencils. DP-SBP operators with even  order accuracy in the interior satisfy the same boundary-interior order of accuracy as the traditional SBP, $(\gamma, \zeta) = (p, 2p)$, while odd order interior DP-SBP operators have $(\gamma, \zeta) = (p, 2p + 1  )$. 

We will assume that the exact solution $\mathbf{q}_e$ is sufficiently smooth such that the truncation error $\mathbb{T}_{j}$ given in \eqref{eq:truncation_error} is defined for all grid points $j\in \{0, 1, \cdots, N\}$. The following theorem proves the convergence of the semi-discrete approximation \eqref{eq:semi-discrete-SWE}.
\begin{theorem}\label{thm:error_well_posed}
Consider the error equation \eqref{eqn:error_eqn}. The error satisfies the estimate,
    $$
    \| \mathcal{E}\|_{WP} \leq  K\biggl(\max\limits_{\tau \in [0, t]}\|\mathbb{T}(\tau)\|_{WP}\biggl), \quad  \mu = \max_{\tau \in[0,t]}{\frac{BT(\mathcal{E})}{\| \mathcal{E}\|_{WP}^2}} \leq 0,  \quad K = \frac{e^{\mu t}-1}{\mu}.
    $$
\end{theorem}
\begin{proof}
The proof is analogous to the proof of Theorem~\ref{thm:st_well_posed_semi}.
\end{proof}
Theorem~\ref{thm:error_well_posed} proves the convergence of the semi-discrete approximation \eqref{eq:semi-discrete-SWE}. In particular the theorem shows that the weighted $l_2$ error is bounded above by the truncation error $\mathbb{T}$ of the operators and converges to zero at the rate $O(\Delta{x}^p)$. However, this error estimate is not sharp, and it is less optimal. Nearly-optimal $O(\Delta{x}^{p+1/2})$ and optimal $O(\Delta{x}^{p+1})$ convergence rates can be proven using the Laplace transform technique, (see e.g., \cite{hicken2013summation,gustafsson1995time,svard2019convergence, gustafsson1975convergence,gustafsson1981convergence,svard2019convergence}).

\begin{remark}
    It suffices to only consider a  linear flux, $\mathbf{F} = [\mathbf{F}_{1},\mathbf{F}_{2}]^T$ in the error equation \eqref{eqn:error_eqn}.
    Note that for the nonlinear flux, by the mean-value theorem, we have
    $$
\mathbf{F}(\mathbf{q}) - \mathbf{F}(\mathbf{q}_e) = M \mathcal{E}, \quad
M = \begin{bmatrix}
U  & H \\
g  & U 
\end{bmatrix},
    $$
    where $M$ is the Jacobian of $\mathbf{F}$ evaluated at $\mathbf{q}_0 = [H, U]^T$, where $\mathbf{q}_0(x,t)$ generally depends on space and time, that is $H(x,t)$, $U(x, t)$ and $W(x,t)=\begin{bmatrix}
g & U \\
U & H
\end{bmatrix}$.
    This extends the convergence analysis to the nonlinear case, under the assumption of smoothness of the exact solution $\mathbf{q}_e$  and the flux $\mathbf{F}$ in $[0,L]$, with the error estimate,
    \begin{align*}
    \|\mathcal{E}\|_W \leq 
    Ke^{\mu_0 t}\biggl(\max\limits_{\tau \in [0, t]}\|\mathbb{T}(\tau)\|_{WP}\biggl), \quad  \mu = \max_{\tau \in[0,t]}{\frac{BT(\mathcal{E})}{\| \mathcal{E}\|_{WP}^2}} \leq 0, \\  \quad {K = \max\left(1, {(1-e^{-\mu_0 t})}/{\mu_0}\right)},
    \quad \mu_0 = \mu+\alpha,
\end{align*}
where $0\le \alpha \le \max_{\tau \in [0,t]}(|W_t|/\min(2|\lambda_W^{-}|, 2|\lambda_W^{+}|))$.
\end{remark}
\section{Energy/entropy stable hyper-viscosity}\label{sec:hyperviscosity}
The continuous and numerical analyses for nonlinear energy stability in the previous sections show that, away from boundaries, energy/entropy is conserved when $\mathbf{G}=0$. For physically consistent solutions, this is obviously true for smooth solutions. For non-smooth solutions energy/entropy must be dissipated across shocks and discontinuities. To minimise unwanted oscillations across shocks, here, we introduce nonlinear energy/entropy stable hyper-viscosity.

\subsection{Semi-discrete approximation with hyper-viscosity}
We denote the viscosity operator $\mathcal{K}$ and make the ansatz
\begin{equation}\label{eq:SWE_Diss}
\mathcal{K}=\mathbf{W}^{-1}\left(I_2\otimes P^{-1}\right)\left(I_2\otimes \mathcal{A}\right), \quad \mathcal{A} = \mathcal{A}^T, \quad \mathbf{u}^T\mathcal{A}\mathbf{u} \le 0, \quad \forall \mathbf{u}\in \mathbb{R}^{N+1}.
\end{equation}
Here $P$ is the diagonal  SBP norm \eqref{eq:SBP-norm}. For the linear problem  $\mathbf{W} = W\otimes\mathbf{I}$ where $W$ is the weight matrix given by  \eqref{eqn:WeightMatrix}, and for the nonlinear problem  $\mathbf{W}$ is given by \eqref{def:wpnorm_weight}. Note that the matrix products above commute, that is $\mathbf{W}^{-1}\left(I_2\otimes P^{-1}\right) = \left(I_2\otimes P^{-1}\right)\mathbf{W}^{-1}$ and $\mathbf{W}\left(I_2\otimes P\right) = \left(I_2\otimes P\right)\mathbf{W}$.
The semi-discrete approximation with hyper-viscosity is obtained by appending the (nonlinear) hyper-viscosity operator $\mathcal{K}$ to the right hand side of \eqref{eq:semi-discrete-SWE} giving
\begin{equation}\label{eq:semi-discrete-SWE_Diss}
   \frac{d\mathbf{q}}{dt} + \mathbf{D}_x \mathbf{F} = \mathcal{K} \mathbf{q} + \mathbf{G} + \mathbf{SAT}, \quad \mathbf{q}(0) =  \mathbf{f}.
\end{equation}
We can prove the following linear and nonlinear stability results for the semi-discrete approximation \eqref{eq:semi-discrete-SWE_Diss}.
\begin{theorem}\label{thm:swp_Diss_linear}
    Consider the semi-discrete approximation \eqref{eq:semi-discrete-SWE_Diss} where the difference operator $\mathbf{D}_x$ is given by \eqref{eq:semi-discrete-SWE}, the  SAT vector is given by \eqref{eq:SAT-term} and the hyper-viscosity operator $\mathcal{K}$ defined by \eqref{eq:SWE_Diss}.
    For linear fluxes  $\mathbf{F}_{1} = U\mathbf{h} + H\mathbf{u}$ and $\mathbf{F}_{2} = U\mathbf{u}+g\mathbf{h}$ and subcritical flows, with $\textrm{Fr} = |U|/\sqrt{gH} < 1$, if the penalty parameters  $\tau_{ij}$ are chosen as in Lemma \ref{lem:stable_penalty} such that $\textrm{BT} \leq 0$, where $\textrm{BT}$ is given by \eqref{eq:BT_linear}, then  the semi-discrete approximation \eqref{eq:semi-discrete-SWE_Diss}  is strongly stable. 
That is, the numerical solution $\mathbf{q}$ satisfies the estimate,
    $$
\begin{aligned}
    \| \mathbf{q}\|_{WP} \leq K e^{\mu t}\biggl(\| \mathbf{f}\|_{WP} +  \max\limits_{\tau \in [0, t]}\|\mathbf{G(\tau})\|_{WP}\biggl), \quad \mu = \max_{\tau \in[0,t]}{\frac{BT + \mathbf{q}^T\left(I_2\otimes \mathcal{A}\right)\mathbf{q}}{\| \mathbf{q}\|_{WP}^2}} \leq 0,  \quad {K = \max\left(1, {(1-e^{-\mu t})}/{\mu}\right)}.
\end{aligned}
$$
\end{theorem}

\begin{theorem}\label{thm:swp_Diss_nonlinear}
    Consider the semi-discrete approximation \eqref{eq:semi-discrete-SWE_Diss} where the difference operator $\mathbf{D}_x$ is given by \eqref{eq:semi-discrete-SWE}, the  SAT vector is given by \eqref{eq:SAT-term} and the nonlinear hyper-viscosity  operator $\mathcal{K}$ defined by \eqref{eq:SWE_Diss}.
    For nonlinear fluxes  $\mathbf{F}_1 = \mathbf{u}\mathbf{h}$ and $\mathbf{F}_2 = \mathbf{u}^2/2 + g\mathbf{h}$ and subcritical flows with $\textrm{Fr} = |u|/\sqrt{gh} < 1$, if the penalty parameters  $\tau_{ij}$ are chosen as in Lemma \ref{lem:stable_penalty} such that $\textrm{BT} \leq 0$, where $\textrm{BT}$ is given by \eqref{eq:BT_linear}, then  the semi-discrete approximation \eqref{eq:semi-discrete-SWE_Diss}  is strongly stable. 
That is, the numerical solution $\mathbf{q}$ satisfies the estimate
    $$
\begin{aligned}
    \| \mathbf{q}\|_{W^{\prime}P} \leq K e^{\mu t}\biggl(\| \mathbf{f}\|_{W^{\prime}P} +  \max\limits_{\tau \in [0, t]}\|\mathbf{G(\tau})\|_{W^{\prime}P}\biggl), \quad \mu = \max_{\tau \in[0,t]}{\frac{BT + \mathbf{q}^T\left(I_2\otimes \mathcal{A}\right)\mathbf{q}}{\| \mathbf{q}\|_{W^{\prime}P}^2}} \leq 0,  \quad {K = \max\left(1, {(1-e^{-\mu t})}/{\mu}\right)}.
\end{aligned}
$$
\end{theorem}

The proofs of Theorems \ref{thm:swp_Diss_linear} and \ref{thm:swp_Diss_nonlinear} follow similarly from the proofs Theorems \ref{thm:st_well_posed_semi} and \ref{thm:nl_semi_stable}.
\subsection{Hyper-viscosity operator}
In order to construct the hyper-viscosity operator $\mathcal{A}$, we consider the even order ($2p$, $p = 1, 2, \ldots$) derivative operator 
$$
 (-1)^{p-1}\alpha\frac{\partial^p }{\partial x^p} c(x)\frac{\partial^p }{\partial x^p}, \quad \alpha \ge 0, \quad c(x) \ge 0,
$$
and approximate it with the DP-SBP operators on the grid. Here $\alpha$ is a real constant and $c(x)$ is a positive smooth function  that vanishes at the boundaries, $x =0, L$ and whose derivatives (up to $(p-1)$th derivatives) also vanish at the boundaries, $x =0, L$. Figure~\ref{fig:boxcar} shows an example of the smooth function $c(x)$ which vanishes on the boundaries and its first and second derivatives also vanish on the boundaries. 
Let $\mathbf{c} = \operatorname{diag}([c(x_0), c_(x_1), \cdots, c(x_N)])$,
the $4$th order hyper-viscosity operator, is given by
$$
\begin{aligned}
 P^{-1}\mathcal{A}= - \alpha D_+ D_- \mathbf{c} D_+D_- =   -\alpha P^{-1}\left(D_-^T P D_- \left(\mathbf{c} P^{-1}\right)D_-^TPD_-\right), \\ \quad \alpha = \delta \Delta x^3,  \quad \delta \ge 0,
\end{aligned}
$$
and the $6$th order hyper-viscosity operator given by
$$
\begin{aligned}
 P^{-1}\mathcal{A}=  \alpha D_-D_+ D_- \mathbf{c} D_+D_-D_+ =   -\alpha P^{-1}\left(D_+^T P D_+ P^{-1}D_+^T \left(P\mathbf{c} \right)D_+P^{-1}D_+^TPD_+\right), \\ \quad \alpha = \delta \Delta x^5, \quad \delta \ge 0.
\end{aligned}
$$
On the right hand sides of $P^{-1}\mathcal{A}$, we have eliminated boundary contributions using the fact that the smooth function $c(x)$ and its first and second derivatives vanish on the boundaries, that is  $c(x)=c^{\prime}(x)=c^{\prime\prime}(x)=0$ at $x =0, L$.
Note that since $P$ and $\mathbf{c}$ are  diagonal matrices then the products $\mathbf{c} P^{-1}$ and $\mathbf{c} P$ are also diagonal matrices. Therefore, for $\alpha \ge0$ the dissipation operators 
$\mathcal{A} = -\alpha \left(D_-^T P D_- \left(\mathbf{c} P^{-1}\right)D_-^TPD_-\right)$  and $\mathcal{A} =-\alpha \left(D_+^T P D_+ P^{-1}D_+^T \left(P\mathbf{c} \right)D_+P^{-1}D_+^TPD_+\right)$ are  symmetric and negative semi-definite. Similar as above, higher order hyper-viscosity operators can also be derived using the DP-SBP operators. However, we will require that the corresponding higher derivatives of $c(x)$ should vanish at the boundaries.

\revi{
    The hyper-viscosity operators $P^{-1}\mathcal{A}$ are norm-compatible with the DP SBP operators and therefore enable the  prove both linear and nonlinear stability for discrete approximations of PDEs that employ DP SBP operators. However, when compared with the the artificial viscosity operators developed for traditional SBP methods by \cite{Mattsson_etal2004}, the hyper-viscosity operators presented in this paper has 2 additional stencil points for even order accurate operators. So, they are slightly more expensive than the artificial viscosity operators developed for traditional SBP methods \cite{Mattsson_etal2004}, but the computational cost does not increase by a large amount. 
    }
\begin{remark}
    Theorems \ref{thm:st_well_posed_semi} \ref{thm:nl_semi_stable}, \ref{thm:swp_Diss_linear}, \ref{thm:swp_Diss_nonlinear} for linear and nonlinear strong stability of the semi-discrete approximation \eqref{eq:semi-discrete-SWE} and the convergence result Theorem \ref{thm:error_well_posed} hold for standard DP upwind and $\alpha$-DRP operators closed with periodic FD  stencils where $\mathbf{SAT} \equiv 0$.
\end{remark}
\subsection{Discrete eigen-spectrum}
\revi{Further, to add credence to the linear and nonlinear stability analysis performed in the last sections we will compute the numerical eigen-spectrum of the evolution operator for the semi-discrete approximation \eqref{eq:semi-discrete-SWE} given by
\begin{equation}\label{eq:evolution_operator}
   \vb{RHS}(\mathbf{q}) = - \mathbf{D}_x \mathbf{F}  + \mathbf{SAT},  \quad \mathbf{D}_x = \begin{bmatrix}
{D}_{+}  & \mathbf{0} \\
\mathbf{0}  & {D}_{-} 
\end{bmatrix}, \quad \mathbf{SAT} = \begin{bmatrix}
\textrm{SAT}_1  \\
\textrm{SAT}_2 
\end{bmatrix}, \quad 
\mathbf{q} = \begin{bmatrix}
\vb{h}  \\
\vb{u} 
\end{bmatrix} \in \mathbb{R}^{2(N+1)}.
\end{equation}
For the linear problem the evolution operator can be written as a matrix vector product, that is $\vb{RHS}(\mathbf{q}) = \vb{A}\mathbf{q}$, where $\vb{A} \in \mathbb{R}^{m \times m}$ with $ m =2(N+1)$ is the evolution matrix. As in \cite{gassner2022stability}, the corresponding evolution matrix $\vb{A}$ for the linearised nonlinear discrete operator $\vb{RHS}(\mathbf{q})$  will be approximated by standard FD method.  
}
\revi{\subsubsection{Stability analysis for the  linear operator}}
To verify the linear stability  analysis we compute the eigenvalues of the semi-discrete spatial evolution operator including the SAT penalty terms which enforce the BCs. We consider specifically mass flux, velocity flux and transmissive BCs, with hyper-viscosity $\alpha >0$ and without hyper-viscosity $\alpha =0$ respectively. The numerical eigenvalues are displayed in Figure \ref{fig:eigenspectrum}. Note that there are no eigenvalues with positive real parts, which verifies the stability proofs of Theorems \ref{thm:st_well_posed_semi} and \ref{thm:swp_Diss_linear}. For the mass flux BC and the velocity flux BC,  without the hyper-viscosity the eigenvalues of the spatial evolution operators are purely imaginary, that is with zero real parts. The addition of hyper-viscosity moves the eigenvalues to the negative half-plane of the complex plane.  Note however, with hyper-viscosity, the magnitude of the real parts of the eigenvalues are several order of magnitude ($\sim 10^{-6}$) smaller than the imaginary parts of the eigenvalues. This also implies that the addition of hyper-viscosity has a negligible impact on the stable time-step for explicit numerical time-integration of the semi-discrete approximation \eqref{eq:semi-discrete-SWE}.

\begin{figure}
    \centering
    \includegraphics[width = 0.5\linewidth]{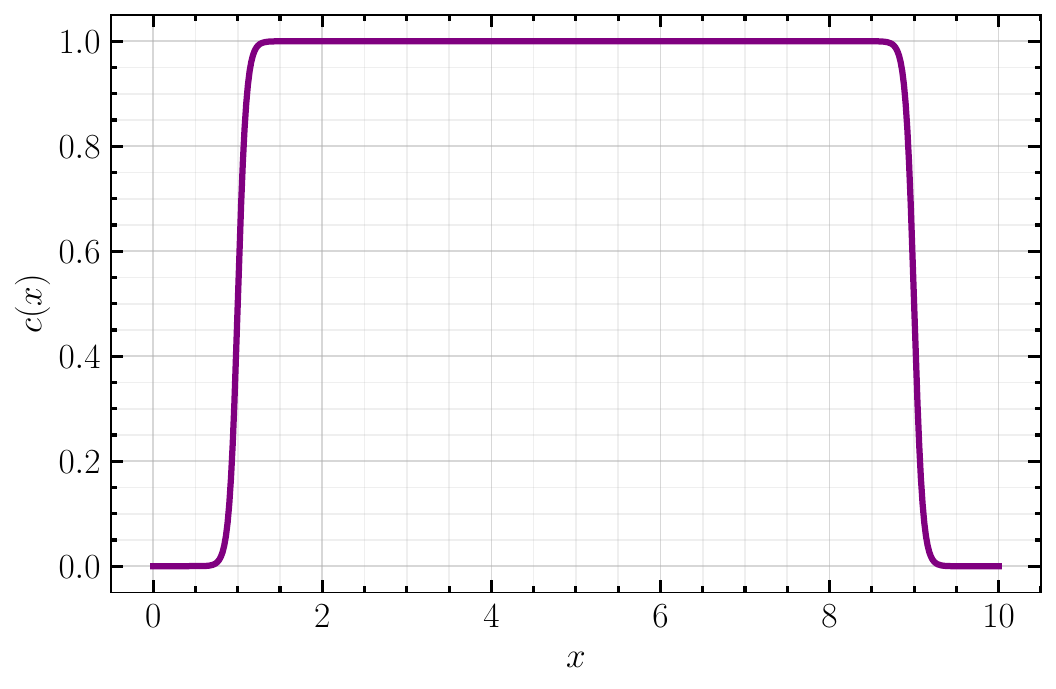}
    \caption{Smooth boxcar function used for the hyper-viscosity operator.}
    \label{fig:boxcar}
\end{figure}
\begin{figure}
    \begin{subfigure}{0.5\textwidth}
    \centering
    \includegraphics[ width = 1\linewidth]{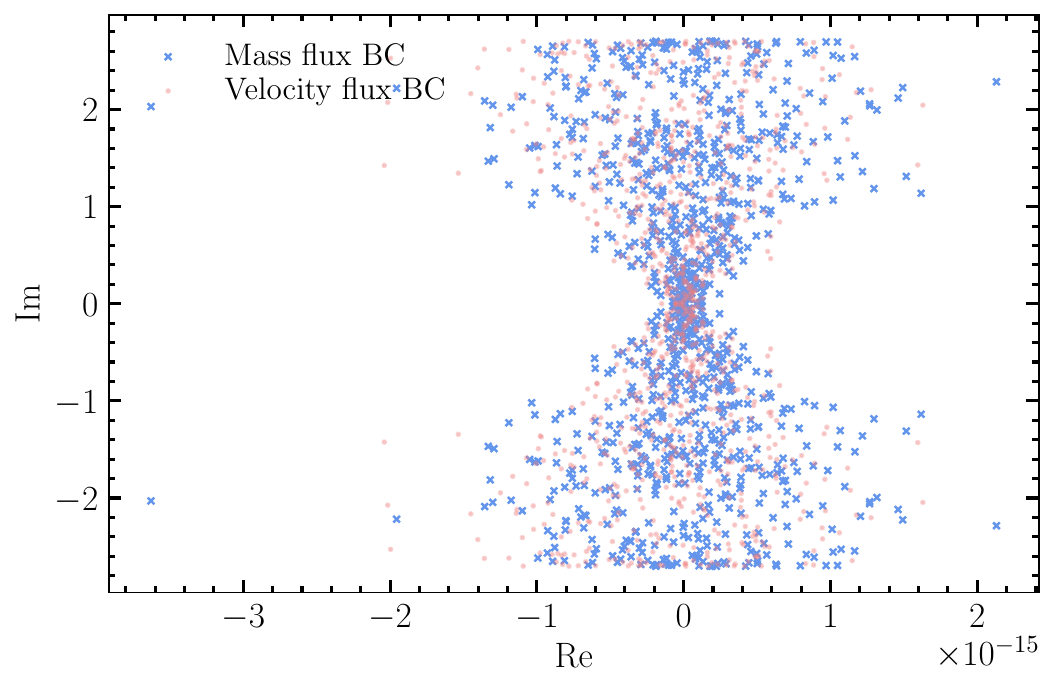}
    \caption{\label{fig:lake_small}$\alpha = 0$}
    \end{subfigure}
    \begin{subfigure}{0.5
\textwidth}
    \centering
    \includegraphics[width = 1\linewidth]{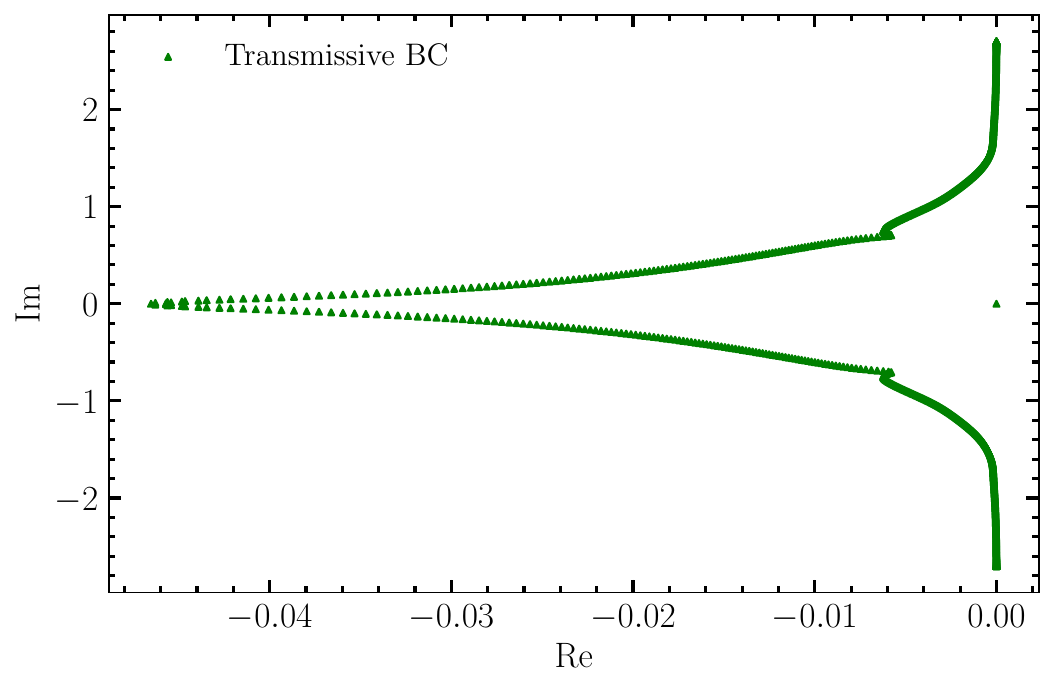}
    \caption{\label{fig:no_hyper_u}$\alpha = 0$}
    \end{subfigure}
        \begin{subfigure}{0.5\textwidth}
    \centering
    \includegraphics[ width = 1\linewidth]{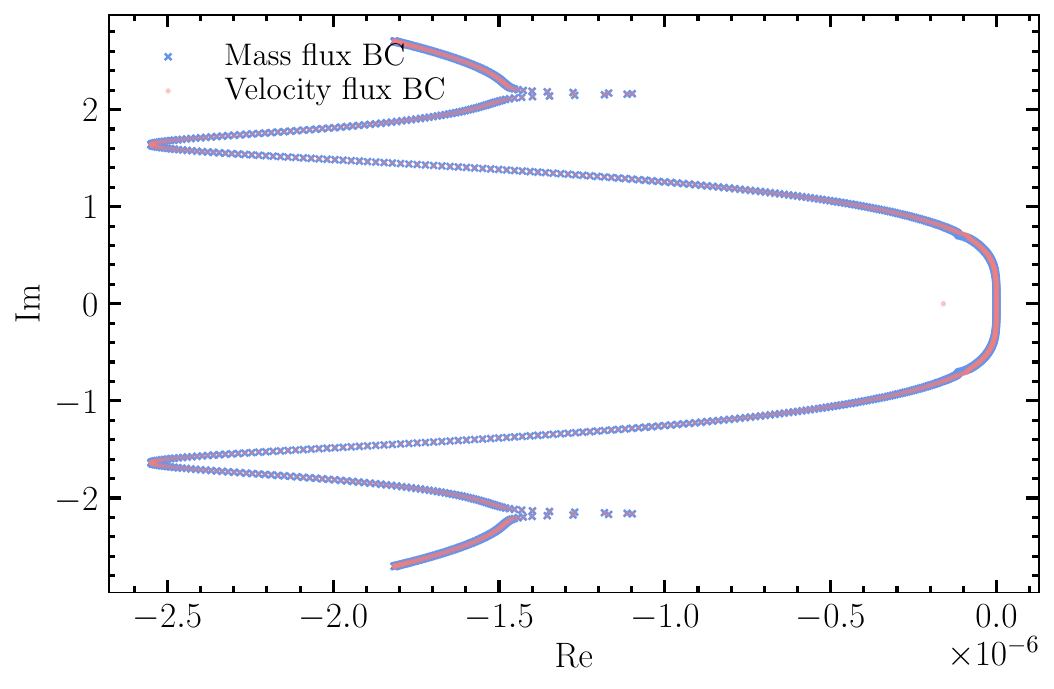}
    \caption{\label{fig:with_hyper_h}$\alpha = 0.1\Delta{x}^5$}
    \end{subfigure}
    \begin{subfigure}{0.5
\textwidth}
    \centering
    \includegraphics[width = 1\linewidth]{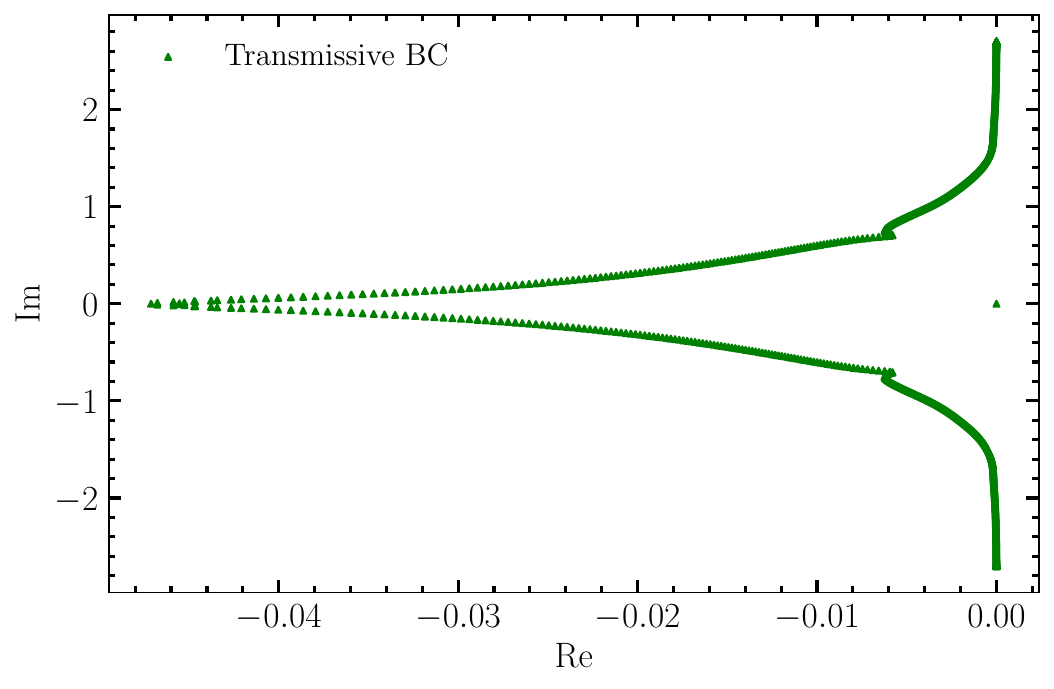}
    \caption{\label{fig:with_hyper_u}$\alpha = 0.1\Delta{x}^5$}
    \end{subfigure}
\caption{\label{fig:eigenspectrum} Numerical eigenspectra of the DP operators of order 6 with $N = 501$, $g=H = 1$, for the three main BCs considered using SAT terms: mass flux, velocity flux and transmissive BCs. The first two are energy conserving without hyper-viscosity and maintains zero real parts up to machine error, while adding hyper-viscosity and utilising transmissive BCs should dissipate energy, which yields  non-positive real parts. All spectra clearly indicate strict stability of the eigensystem.}
\end{figure}
\revi{\subsubsection{Stability analysis for the linearised nonlinear operator with nontrivial background}}
\revi{Here we also perform a local stability analysis of the linearised nonlinear operator under a non-trivial background state, based on the method of Gassner et al. \cite{gassner2022stability} for the Burgers' equation. 
We consider the smooth base flow $H(x)$ and $U(x)$ given by
\begin{equation}
H(x) = \theta_h \sin\bigl(2\pi(x + 0.7)\bigl) + 2, \quad 
U(x) = \theta_u \cos\bigl(2\pi(x - 0.7)\bigl), \quad x \in [0,1], \quad \theta_h = \theta_u=0.1
\end{equation}
Then, the Jacobian matrix $\mathbf{A}$ for the linearised operator can be approximated via a centred difference approach as,
\begin{equation}
\mathbf{A}_{:, j} \approx \frac{\mathbf{RHS}(\mathbf{q} + \varepsilon \mathbf{e}_j) - \mathbf{RHS}(\mathbf{q} - \varepsilon \mathbf{e}_j)}{2 \epsilon}, \quad 
\mathbf{q} = \begin{bmatrix}
\vb{H}  \\
\vb{U} 
\end{bmatrix} \in \mathbb{R}^{2(N+1)}.
\end{equation}
where $\mathbf{A}_{:,j}$ is the $j$-th column of $\mathbf{A}$, and $\mathbf{e}_j$ is the $j$-th unit vector, with $1$ at the \(j\)-th position and $0$ elsewhere, \(\varepsilon\) is a small perturbation value which we set to $1\times 10^{-6}$. The computed eigenvalues of the approximative Jacobi matrix give an approximate spectrum of linearised evolution operator for a smooth and nontrivial background flow.

The numerical eigenvalues are displayed in Figure \ref{fig:eigenspectrum_nonlinear}. Similar to the linear case, note that there are no eigenvalues with significant positive real parts, which verifies the local stability result and  convergence proof of Theorem \ref{thm:error_well_posed}. Again for the mass flux BC and the velocity flux BC with $\alpha =0$,  without the hyper-viscosity,  the eigenvalues of the spatial evolution operators are purely imaginary, that is with real parts $\sim 10^{-11}$, which is the order of rounding errors. For the transmissive BC with $\alpha =0$,  without the hyper-viscosity, there are a few eigenvalues with negative real parts but the  majority of the eigenvalues lie on the imaginary line in the complex plane. As expected, the addition of hyper-viscosity moves the eigenvalues to the negative half-plane of the complex plane.  Note in particular, with hyper-viscosity, the eigenvalues corresponding to unresolved high frequency modes are moved further into the negative complex plane.
}
\begin{figure}
    \begin{subfigure}{0.5\textwidth}
    \centering
    \includegraphics[ width = 1\linewidth]{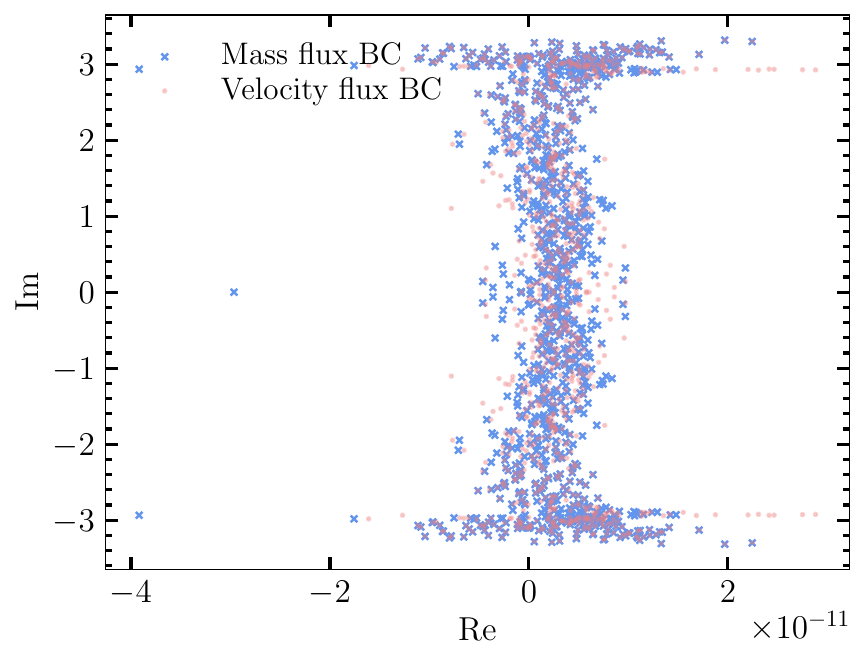}
    \caption{\label{fig:mass_flux_nonlinear}$\alpha = 0$}
    \end{subfigure}
    \begin{subfigure}{0.5
\textwidth}
    \centering
    \includegraphics[width = 1\linewidth]{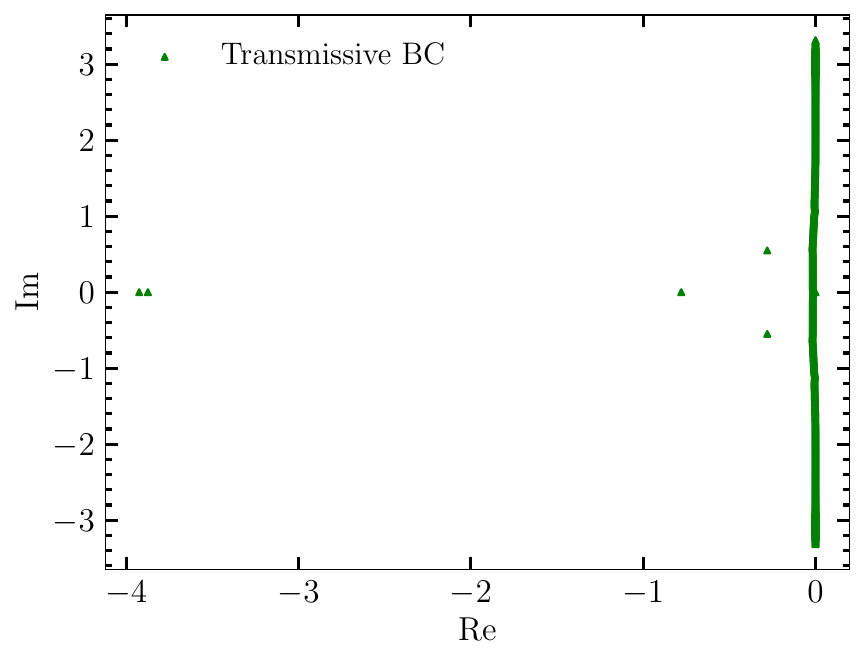}
    \caption{\label{fig:no_hyper_u_nl}$\alpha = 0$}
    \end{subfigure}
        \begin{subfigure}{0.5\textwidth}
    \centering
    \includegraphics[ width = 1\linewidth]{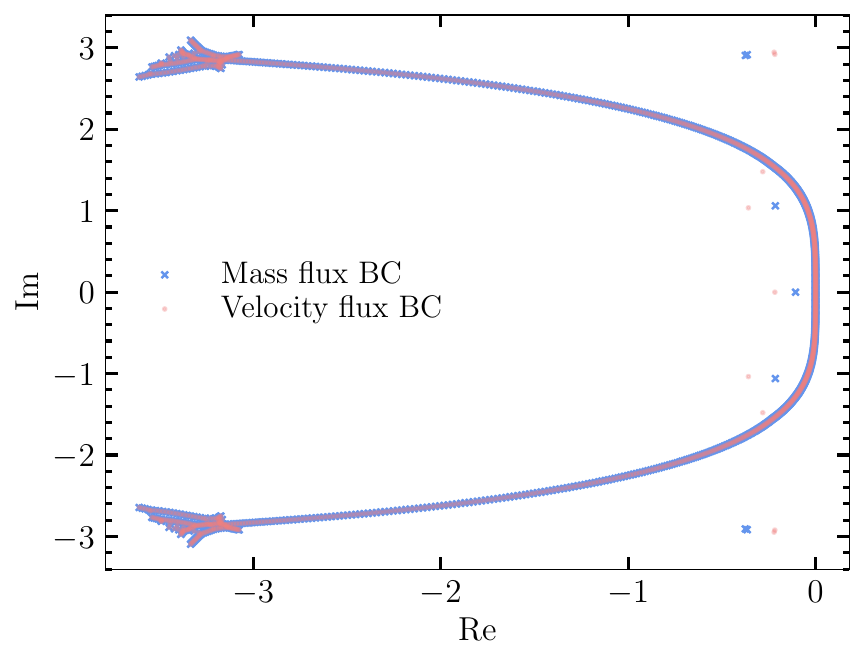}
    \caption{\label{fig:with_hyper_h_nl}$\alpha = 0.1\Delta{x}^5$}
    \end{subfigure}
    \begin{subfigure}{0.5
\textwidth}
    \centering
    \includegraphics[width = 1\linewidth]{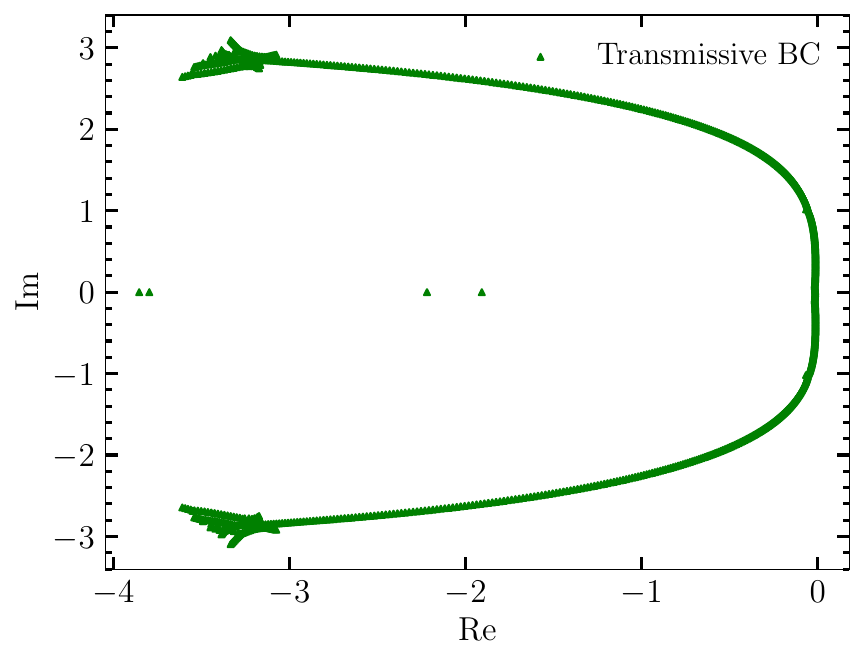}
    \caption{\label{fig:with_hyper_u_nl}$\alpha = 0.1\Delta{x}^5$}
    \end{subfigure}
\caption{\label{fig:eigenspectrum_nonlinear} \revi{Numerical eigenspectra of the linearised non-linear operator for the DP operators of order 6 with $N = 501$, $g=1$, with the three main BCs considered using SAT terms: mass flux, velocity flux and transmissive BCs. The first two are energy conserving without hyper-viscosity and maintains zero real parts up to machine error, while adding hyper-viscosity and utilising transmissive BCs should dissipate energy, which yields non-positive real parts. All spectra clearly indicate strict stability of the eigensystem.}}
\end{figure}

In the next section  we present numerical experiments in 1D verifying the theory derived in this paper. We will also present some 2D numerical simulations showing extensions of the method to the 2D nonlinear rotating SWEs.

\section{Numerical experiments}\label{sec:numexp}
In this section we present detailed numerical experiments verifying the theoretical analysis performed in the previous sections. In particular the experiments  are designed to verify stability, accuracy and convergence properties of the method, including verification of the nonlinear transmissive BCs. We also present canonical test problems such as the dam break with wet domain, and the well-balanced test called the lake at rest with non-smooth and smooth immersed bump. Numerical examples in 2D within a doubly periodic domain are presented showing extension of the method to 2D and effectiveness of the hyper-viscosity operator for simulating merging vortexes and fully developed turbulent flows.

For time discretisation, we use the  classical fourth order accurate explicit Runge-Kutta method.
Throughout the 1D numerical experiments, we set the global time step 
\begin{equation}
    \Delta t = \frac{\textrm{CFL} \times \Delta x}{\max_x{(|U| +\sqrt{gH})}},
\end{equation}
where $\Delta x = L/N$,  $N+1$ is the number of grid points and $\textrm{CFL} = 0.3$.  We use natural units, so that $g = 9.81$, $ H = 1 $, $U= -0.3\sqrt{gH}$, $|U| +\sqrt{gH}$ is the characteristic wave speed of the linear problems, and for the nonlinear problem we use $U=u(t=0,x)$ and $H=h(t=0,x)$. 


\subsection{Numerical experiments in 1D}
We first study the accuracy of the numerical method via the method of manufactured solutions (MMS), before moving to canonical test problems such as the dam break problem with wet domain, and the well-balanced test called the lake at rest with non-smooth and smooth immersed bump, as  proposed in \cite{delestre2013swashes}.

\subsubsection{Method of manufactured solutions}

\begin{table}

\centering

	\resizebox{\textwidth}{!}	{\begin{tabular}{c c c c c c c c c}
 \hline  \hline
	  \multicolumn{1}{c}{\texttt{DP, order 4, MMS} } \\	
& \multicolumn{4}{c}{Linear} 
&
\multicolumn{4}{c}{Nonlinear}  \\ 
		\multicolumn{1}{c}{$N$} & $\log_2 \| \textrm{error}\|_{2,u}$ & $\log_2\|\textrm{error}\|_{2,h}$ & $q_u$    & $q_h$     & $\log_2 \| \textrm{error}\|_{2,u}$ & $\log_2\|\textrm{error}\|_{2,h}$ & $q_u$    & $q_h$     \\
    \multicolumn{1}{c}{(1)} & (2) & (3) & (4) & (5) & (6) & (7) & (8) & (9)\\ \hline
		41  & -6.5665                              &  -6.5454                           & -        & -  &  -5.5362 & -5.7156 & - & -     \\ 
		81  & -10.8428                            & -10.8198                         & 4.2763 & 4.2743 & -8.9666 & -9.0302 & 3.4305 & 3.3146\\ 
		161  & -14.8953                            & -14.8755                         & 4.0525 & 4.0557 & -12.8307 & -13.0011 & 3.9284 & 4.0598  \\
		321 & -18.9189                            & -18.9011                          & 4.0236 & 4.0256 & -16.7592 & -17.0609 & 3.9284 & 4.0598\\ 
		641 & -22.9305                           & -22.9139                        & 4.0116 & 4.0127 & -20.7279 & -21.0780 & 3.9687 & 4.0171 \\ \hline \hline
	\end{tabular}}
 \resizebox{\textwidth}{!}{\begin{tabular}{c c c c c c c c c}

  \multicolumn{1}{c}{\texttt{DRP, order 4, MMS} } \\
	& \multicolumn{4}{c}{Linear} 
&
\multicolumn{4}{c}{Nonlinear}  \\ 
		\multicolumn{1}{c}{$N$} & $\log_2 \| \textrm{error}\|_{2,u}$ & $\log_2\|\textrm{error}\|_{2,h}$ & $q_u$    & $q_h$     & $\log_2 \| \textrm{error}\|_{2,u}$ & $\log_2\|\textrm{error}\|_{2,h}$ & $q_u$    & $q_h$     \\
    \multicolumn{1}{c}{(1)} & (2) & (3) & (4) & (5) & (6) & (7) & (8) & (9)\\ \hline
		41  & -4.2724                              &  -4.7107                           & -        & - & -4.2123 & -4.3524 & - & -          \\ 
		81  & -10.7031                            & -10.7447                         & 6.4307 & 6.0341 & -9.3048 & -8.5032 &  5.0925 & 4.1508\\ 
		161  & -10.7031                             & -10.7447         & 4.4595 & 4.3913 & -12.9814 & -12.8435 & 3.6765 & 4.3403 \\
		321 & -15.1626             & -15.1360                          & 4.4595 & 4.3913 & -16.8969 & -17.1903 & 3.9155 & 4.3468 \\ 
		641 & -23.1426                           & -23.1243       & 3.9983 & 4.0011  & -20.8564 & -21.2958 & 3.9596 & 4.1054 \\ \hline \hline
	\end{tabular}}

 	\resizebox{\textwidth}{!}{\begin{tabular}{c c c c c c c c c}
	
  \multicolumn{1}{c}{\texttt{DP, order 6, MMS} } \\
	&	\multicolumn{4}{c}{Linear} 
&
\multicolumn{4}{c}{Nonlinear}  \\ 
		\multicolumn{1}{c}{$N$} & $\log_2 \| \textrm{error}\|_{2,u}$ & $\log_2\|\textrm{error}\|_{2,h}$ & $q_u$    & $q_h$     & $\log_2 \| \textrm{error}\|_{2,u}$ & $\log_2\|\textrm{error}\|_{2,h}$ & $q_u$    & $q_h$     \\
    \multicolumn{1}{c}{(1)} & (2) & (3) & (4) & (5) & (6) & (7) & (8) & (9)\\ \hline
		41  & -5.9761                              &  -5.0590                          & -        & - &    -5.6406   & -4.7045 & - & - \\ 
		81  & -10.5185                            & -10.4705                         & 4.5424 & 5.4115 & -10.2011 & -10.5153 & 4.5605 & 5.8108 \\ 
		161  & -15.5532                            & -15.2167                        & 5.0348 & 4.7461 & -15.1654 & -14.7981 &  4.9643 & 4.2828\\ 
		321 & -20.1258                            & -19.7233                          & 4.5726 & 4.5067 & -19.7564 & -19.2424 & 4.5910  & 4.4443 \\ 
		641 & -24.7807                           & -24.2329                        & 4.6549 & 4.5095  & -24.4788 & -23.7511 & 4.7224 & 4.5087\\ \hline \hline
	\end{tabular}}
  \resizebox{\textwidth}{!}{\begin{tabular}{c c c c c c c c c}
  \multicolumn{1}{c}{\texttt{DRP, order 6, MMS} } \\
	&	\multicolumn{4}{c}{Linear} 
&
\multicolumn{4}{c}{Nonlinear}  \\ 
		\multicolumn{1}{c}{$N$} & $\log_2 \| \textrm{error}\|_{2,u}$ & $\log_2\|\textrm{error}\|_{2,h}$ & $q_u$    & $q_h$     & $\log_2 \| \textrm{error}\|_{2,u}$ & $\log_2\|\textrm{error}\|_{2,h}$ & $q_u$    & $q_h$     \\
    \multicolumn{1}{c}{(1)} & (2) & (3) & (4) & (5) & (6) & (7) & (8) & (9)\\ \hline
		41  & -3.9230      &  -5.2152                         & -        & -   &  -3.1407 & -4.1012 & - & -  \\ 
		81  & -8.3058                            & -9.0906                         & 4.3827 & 3.8754 & -8.0456 & -8.6623 & 4.9049 & 4.5611 \\ 
		161  & -14.3327                            & -14.1198                       & 6.0270 & 5.0292 & -13.7421 & -13.7835 &  5.6965 &  5.1212 \\ 
		321 & -18.9399                           & -18.3042                         & 4.6072  & 4.1844 & -18.3152 & -17.6874 & 4.5731 & 3.9039 \\
  		641 & -23.4452                           & -22.7525                         & 4.5053  & 4.4483  & -22.8357 & -22.1028 & 4.5205 & 4.4154\\ 
	 \hline \hline
	\end{tabular}}
 \caption{\label{table:1}Log$_2 L^2$ norm errors and convergence rate \revi{$q_u$ and $q_h$ of velocity and height} of various FD operators using MMS, for the primitive variables $u$, $h$.  }

\end{table}
 In this section, we use MMS to verify the convergence of our numerical scheme.
We consider both linear and nonlinear problems. The computational domain is $x \in [0, L]$ with  length $L = 10$. We force the system with the exact solution
\begin{equation}
\begin{aligned}
  u_e(x,t) = \exp(-(x-x_0-c_s t)^2) , \quad h_e(x,t) = u_e(x,t) + 10, \quad x_0 = 0.5L, \\ \quad c_s = \sqrt{gH},
  \end{aligned}
\end{equation}
while enforcing the linear mass flux BC $F_1(x,t):=Hu(x,t) + Uh(x,t)=Hu_e(x,t) + Uh_e(x,t)$ at $x = 0, L$ for the linear problem and the nonlinear mass flux BC $F_1(x,t):=u(x,t)h(x,t)=u_e(x,t)h_e(x,t)$ at $x = 0, L$ for the nonlinear problem. \revi{Note that $Hu_e(x,t) + Uh_e(x,t)$ and $u_e(x,t)h_e(x,t)$ yield nonzero boundary data, $g_0(t)=Hu_e(0,t) + Uh_e(0,t)$, $g_L(t)=Hu_e(L,t) + Uh_e(L,t)$ for the linear problem and $g_0(t)=u_e(0,t)h_e(0,t)$, $g_L(t)=u_e(L,t)h_e(L,t)$ for the nonlinear problem.}
%
We run the simulations until the final time $t = 0.5$. 

\paragraph{MMS without hyper-viscosity.}
We consider first grid convergence tests without  hyper-viscosity and proceed later to convergence tests  with hyper-viscosity. The numerical errors and convergence rates for different  operators, without  hyper-viscosity, for the linear and nonlinear cases  are  reported in Table~\ref{table:1}. As has already been noted in prior works \cite{mattsson2017diagonal,mattsson2018compatible,lundgren2020efficient}, the diagonal norm DP upwind operators have ``higher than expected" convergence rates. Our experiments also show that the DRP operators exhibit the same behaviour, higher than expected convergence rates, as  it also satisfy the upwind property. We do not understand the precise reason for this yet, further analysis is needed to unravel the super-convergence properties exhibited by the schemes. For additional convergence studies using traditional SBP operators, and odd-order upwind DP and DRP operators, see~\ref{appen:addconv}.  

\paragraph{MMS with hyper-viscosity.}
We now assess the accuracy of the method with  hyper-viscosity for linear and nonlinear problems. To do this  we use MMS and perform grid convergence studies for smooth solutions. We set $\alpha = \delta {\Delta{x}}^{p-1}$ and $\delta = 0.1$, where $p= 4,6$ for $4$th and $6$th order derivative hyper-viscosity operators. For smooth solutions the truncation error for the hyper-viscosity operators is $\mathbb{T} \sim {\Delta{x}}^{p-1} + \mathcal{O}({\Delta x}^p)$. We expect the convergence rate of the errors to be at most $\mathcal{O}({\Delta x}^{p-1})$.  Table~\ref{table:hyper} shows the error and convergence rates of the error for DP upwind and DRP operators of $4$ and $6$ interior accuracy, for linear and nonlinear problems. 
Hence, our DP and DRP schemes remain high order accurate with  hyper-viscosity for smooth solutions, for both linear and nonlinear problems.
\begin{table*}
\centering
 \resizebox{\textwidth}{!}	{\begin{tabular}{c c c c c c c c c}
		\hline  \hline
  \multicolumn{1}{c}{\texttt{DP order 4, hyper-viscosity, MMS} } \\
& \multicolumn{4}{c}{ Linear} 
&
\multicolumn{4}{c}{Nonlinear}  \\ 
		\multicolumn{1}{c}{$N$} & $\log_2 \| \textrm{error}\|_{2,u}$ & $\log_2\|\textrm{error}\|_{2,h}$ & $q_u$    & $q_h$     & $\log_2 \| \textrm{error}\|_{2,u}$ & $\log_2\|\textrm{error}\|_{2,h}$ & $q_u$    & $q_h$    \\
    \multicolumn{1}{c}{(1)} & (2) & (3) & (4) & (5) & (6) & (7) & (8) & (9)\\ \hline
		41  & -7.1027    &  -6.2709                          & 0        & 0 & -6.7094  &  -5.9420 & 0 & 0        \\ 
		81  & -11.3750                             & -11.0446                        &  4.2723 & 4.7738 & -11.1816 & -10.8197 & 4.4722 & 4.8777    \\ 
		161  & -15.8828                           & -15.6486      & 4.5078 &  4.6040 & -15.4499  & -15.0417 &  4.2683 & 4.2220 \\ 
		321 & -19.4923                         & -19.4737                         &  3.6095 &  3.8251 &  -18.8463 & -18.3501 & 3.3964 &  3.3084\\ 
   \hline \hline 
	\end{tabular}}
  \resizebox{\textwidth}{!}	{\begin{tabular}{c c c c c c c c c}
		
  \multicolumn{1}{c}{\texttt{DRP order 4, hyper-viscosity, MMS} } \\
& \multicolumn{4}{c}{ Linear} 
&
\multicolumn{4}{c}{Nonlinear}  \\ 
		\multicolumn{1}{c}{$N$} & $\log_2 \| \textrm{error}\|_{2,u}$ & $\log_2\|\textrm{error}\|_{2,h}$ & $q_u$    & $q_h$     & $\log_2 \| \textrm{error}\|_{2,u}$ & $\log_2\|\textrm{error}\|_{2,h}$ & $q_u$    & $q_h$    \\
    \multicolumn{1}{c}{(1)} & (2) & (3) & (4) & (5) & (6) & (7) & (8) & (9)\\ \hline
		41  & -7.1027    &  -6.2709                          & 0        & 0   & -3.8329 & -4.01740 & 0 & 0      \\ 
		81  & -11.3750                             & -11.0446                        &  4.2723 & 4.7738 & -8.7963 & -10.9147 & 4.9633 & 6.8973  \\ 
		161  & -15.8828                           & -15.6486      & 4.5078 &  4.6040 & -14.2213 & -14.0322 & 5.4250 & 3.1175 \\ 
		321 & -19.4923                         & -19.4737                         &  3.6095 &  3.8251 & -18.5059 & -17.9750 & 4.2846 & 3.9429\\ 
   \hline \hline 
	\end{tabular}}
  \resizebox{\textwidth}{!}	{\begin{tabular}{c c c c c c c c c}
		
  \multicolumn{1}{c}{\texttt{DP order 6, hyper-viscosity, MMS} } \\
& \multicolumn{4}{c}{ Linear} 
&
\multicolumn{4}{c}{Nonlinear}  \\ 
		\multicolumn{1}{c}{$N$} & $\log_2 \| \textrm{error}\|_{2,u}$ & $\log_2\|\textrm{error}\|_{2,h}$ & $q_u$    & $q_h$     & $\log_2 \| \textrm{error}\|_{2,u}$ & $\log_2\|\textrm{error}\|_{2,h}$ & $q_u$    & $q_h$    \\
    \multicolumn{1}{c}{(1)} & (2) & (3) & (4) & (5) & (6) & (7) & (8) & (9)\\ \hline
		41  & -7.0942    &  -6.2723                         & 0        & 0 & -6.6884 & -5.9477 & 0 & 0         \\ 
		81  &  -11.4008                             & -11.0658 &  4.2723  & 4.7738 & -11.2279 & -10.8404 & 4.5396 & 4.8927 \\ 
		161  & -16.2131                          & -15.8455    & 4.5078 &  4.6040 & -15.9767 & -15.4594 & 4.7487 & 4.6190\\ 
		321 & -20.9717                         & -20.5254                        & 4.7586  &  4.6798 & -20.6898 & -20.1138 & 4.7131 & 4.6544\\ 
   \hline \hline 
	\end{tabular}}
   \resizebox{\textwidth}{!}	{\begin{tabular}{c c c c c c c c c}
		
  \multicolumn{1}{c}{\texttt{DRP order 6, hyper-viscosity, MMS}} \\
& \multicolumn{4}{c}{ Linear} 
&
\multicolumn{4}{c}{Nonlinear}  \\ 
		\multicolumn{1}{c}{$N$} & $\log_2 \| \textrm{error}\|_{2,u}$ & $\log_2\|\textrm{error}\|_{2,h}$ & $q_u$    & $q_h$     & $\log_2 \| \textrm{error}\|_{2,u}$ & $\log_2\|\textrm{error}\|_{2,h}$ & $q_u$    & $q_h$    \\
    \multicolumn{1}{c}{(1)} & (2) & (3) & (4) & (5) & (6) & (7) & (8) & (9)\\ \hline
		41  & -4.2029   &  -4.7415                          & 0        & 0  & -3.8306 & -4.0142 & 0 & 0        \\ 
		81  & -9.2636                            & -10.1468            &  5.0607 & 5.4053  & -8.7986 & -11.1187 & 4.4.968 & 7.1045\\ 
		161  & -14.9135                          & -14.6628                         & 3.8338 & 3.9573 & -14.3236 & -14.0398 & 5.5251 & 2.9211 \\ 
		321 & -19.7461                            & -19.0091                          &  4.8326 & 4.3463 & -19.1884 & -18.4021 & 4.8648 & 4.3623\\ 
   \hline \hline 
	\end{tabular}}
\caption{\label{table:hyper}Log$_2$ $L^2$ norm errors and convergence rate \revi{$q_u$ and $q_h$ of velocity and height} of the DP and DRP operators of order $p =4, 6$ interior operators, with hyper-viscosity pre-factor $\delta = 0.1$.}
\end{table*}
\subsubsection{Lake at rest with immersed bump}
 We consider the canonical 1D lake at rest problem, modeled by the nonlinear SWE \eqref{eqn:nlSWE} with bathymetry (immersed bump), as described in~\cite{delestre2013swashes}. The length of the domain is $L=25$ and  the bathymetry defined as:
\begin{equation}\label{eqn:b(x)}
    b(x)= \begin{cases}0.2-0.05(x-10)^2 & \text { if } 8 <x<12,  \\ 0 & \text { $x \in \Omega \backslash \{8,12\} $},\end{cases}
\end{equation}
The bathymetry $b(x)$ enters the nonlinear SWE \eqref{eqn:nlSWE} through its derivative
\begin{equation}\label{eqn:b_prime}
    {b^{\prime}(x)}= \begin{cases} -0.1(x-10) & \text { if } 8 <x<12,  \\ 0 & \text { $x \in \Omega \backslash \{8,12\} $},\end{cases}
\end{equation}
which appears as a source term in the momentum equation. Note that $b^{\prime}(x)$ is discontinuous at $x = \{8, 12\}$.

\paragraph{No perturbation.}
We will now verify the well-balanced property of the numerical method.
We prescribe the initial condition
$$
h(x,0)+b = 0.5, \quad u(x,0) = 0,
$$
 throughout the domain. As the initial conditions solve a steady state problem, a good numerical scheme should  maintain the prescribed initial condition for all time.
Following \cite{lundgren2020efficient}, we prescribe periodic BCs at $x = 0, L$, with periodic FD stencils. 
No hyper-viscosity is used for this test case, despite the fact that $b(x)$ is non-smooth. 
 As shown in Fig.~\ref{fig:lake_at_rest} and Table~\ref{fig:lake_at_rest}, our scheme maintains both stage, ($h+b$)  and velocity at the analytical value up to machine precision \revi{Note that we have omitted the stage error since it is  identically zero.} These are independent of the grid resolution, and confirms the well-balanced property of the numerical scheme.
\begin{figure}
    \centering
    \includegraphics[width = 0.6\linewidth]{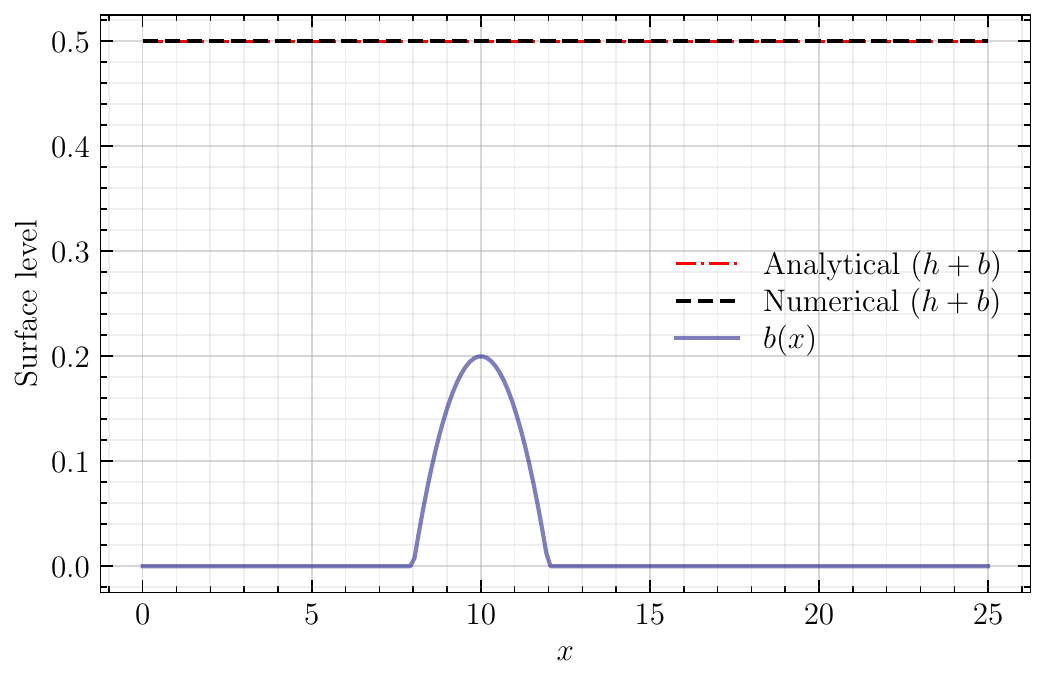}

    \caption{Plot of the height for the canonical lake at rest problem at end time $t = 5$, where the bathymetry $b(x)$ is defined via \eqref{eqn:b(x)}, which is non-smooth about $x \in (8,12)$. The simulation is initialised with $h + b = 0.5$, $N = 201$, $L = 25$. Our numerical scheme maintains the initial condition up to machine precision.   }
    \label{fig:lake_at_rest}
\end{figure}
\begin{table}
\centering 
\begin{tabular}{c c c c  }
		\hline \hline
  \multicolumn{1}{c}{\texttt{FD, order 6, Lake at Rest} } \\
		\multicolumn{1}{c}{} & $\log_{10} \| \textrm{error}\|_{2,u}$ & $\log_{10}\|\textrm{error}\|_{2,u}$ & $\log_{10}\|\textrm{error}\|_{2,u}$        \\
    \multicolumn{1}{c}{$N$} & (SBP) & (DP) & (DRP)\\ \hline
		51  & -14.6270   &  -14.0611 & -15.4622                                  \\ 
		101  & -14.3264 & -13.7623                              & -15.1600 
                     \\ 
		151  & -14.1503 & -13.5869                            & -14.9836                          \\ 
		201 & -14.0253 & -12.7646                            &
 -14.8581                           \\ 
   \hline \hline 
	\end{tabular}

\caption{\label{table:lake}Log$_{10}$ error of the lake at rest's velocity profile, it can be seen that all operators achieve values near machine precision, regardless of grid resolution.}
\end{table}
\paragraph{Small perturbations.}
The aim of this experiment is to verify the implementation and accuracy  of the nonlinear transmissive BCs \eqref{eq:BC-subcriticalTransmissive_nonlinear}.
As above we consider the 1D lake at rest problem, with zero initial condition for velocity, $u(x,0) = 0$, and we add small perturbations to the initial condition for height,
\begin{equation*}\label{eqn:gaussianform}
 h(x,0) = 0.5 - b(x) + \delta{b}(x), \quad   \delta{b}(x) = 0.1 \max_{x}{|b(x)|} \times e^{-\frac{\left(x - 10\right)^2}{0.3}}.
\end{equation*}
The perturbation will generate variations in both height and velocity which will propagate through the domain. However,  the variations introduced by the perturbation are expected to leave the domain through the transmissive boundaries without reflections.
%
\begin{figure}
    \begin{subfigure}{0.5\textwidth}
    \centering
    \includegraphics[ width = 1\linewidth]{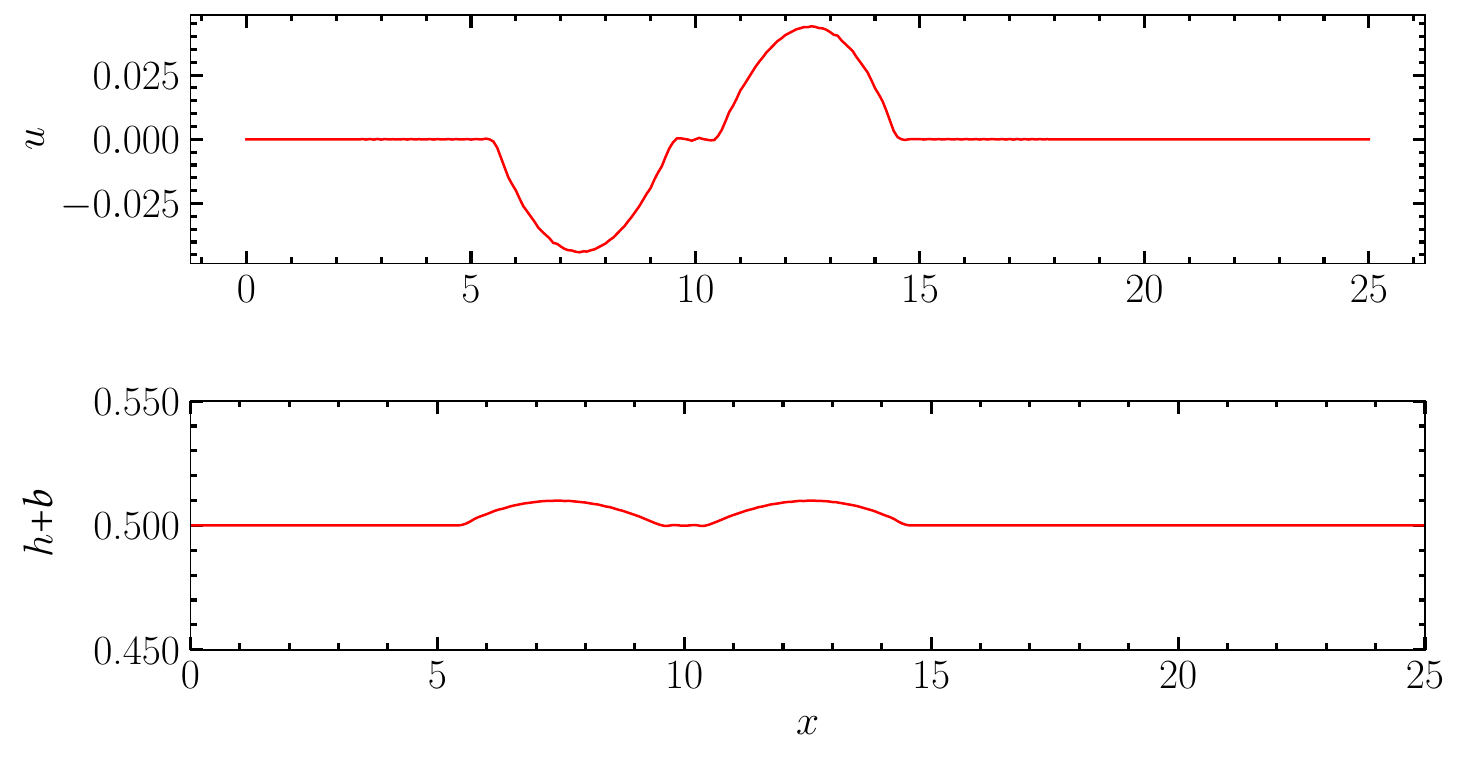}
    \caption{\label{fig:lake_small}$t = 100\Delta t $}
    \end{subfigure}
    \begin{subfigure}{0.5
\textwidth}
    \centering
    \includegraphics[width = 1\linewidth]{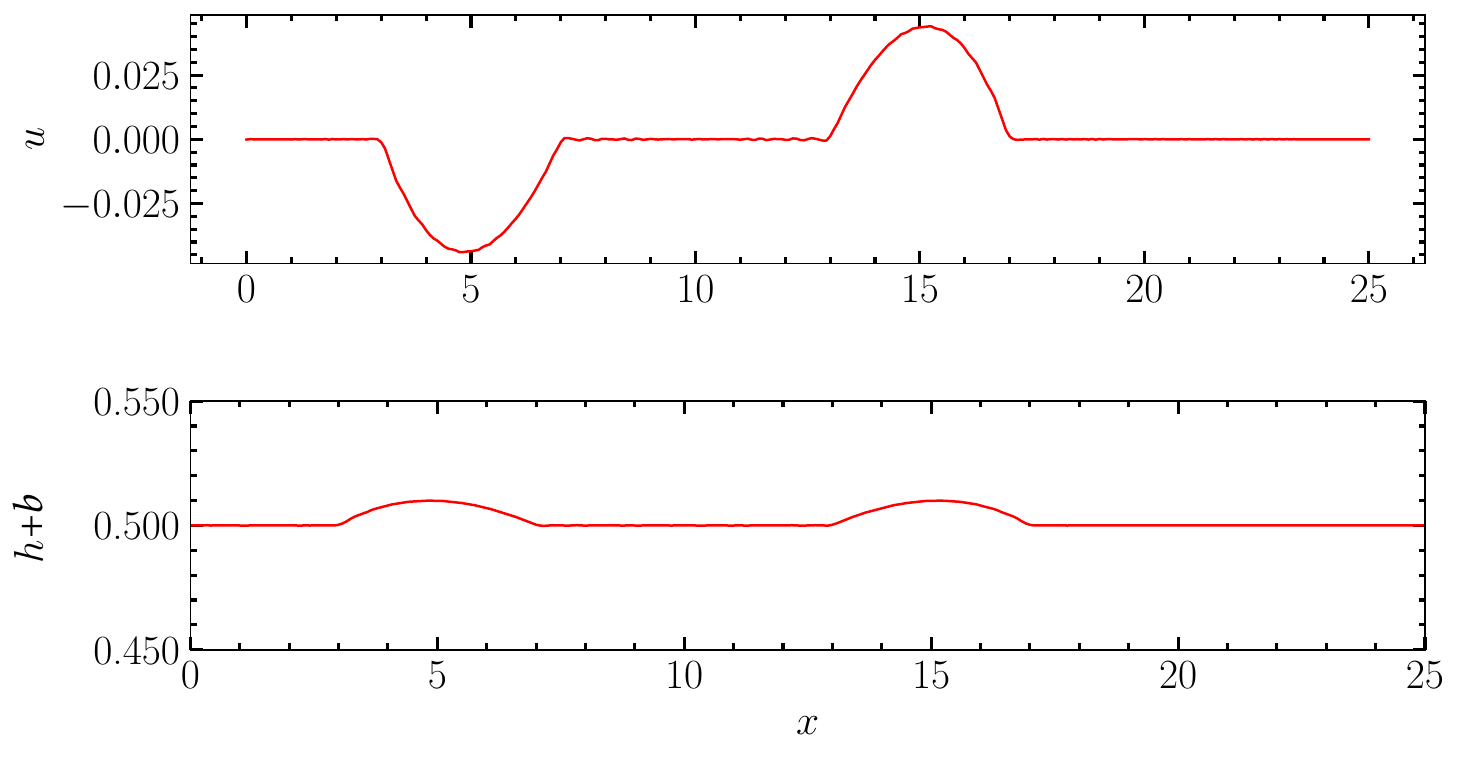}
    \caption{\label{fig:no_hyper_u}$t = 200\Delta t $}
    \end{subfigure}
        \begin{subfigure}{0.5\textwidth}
    \centering
    \includegraphics[ width = 1\linewidth]{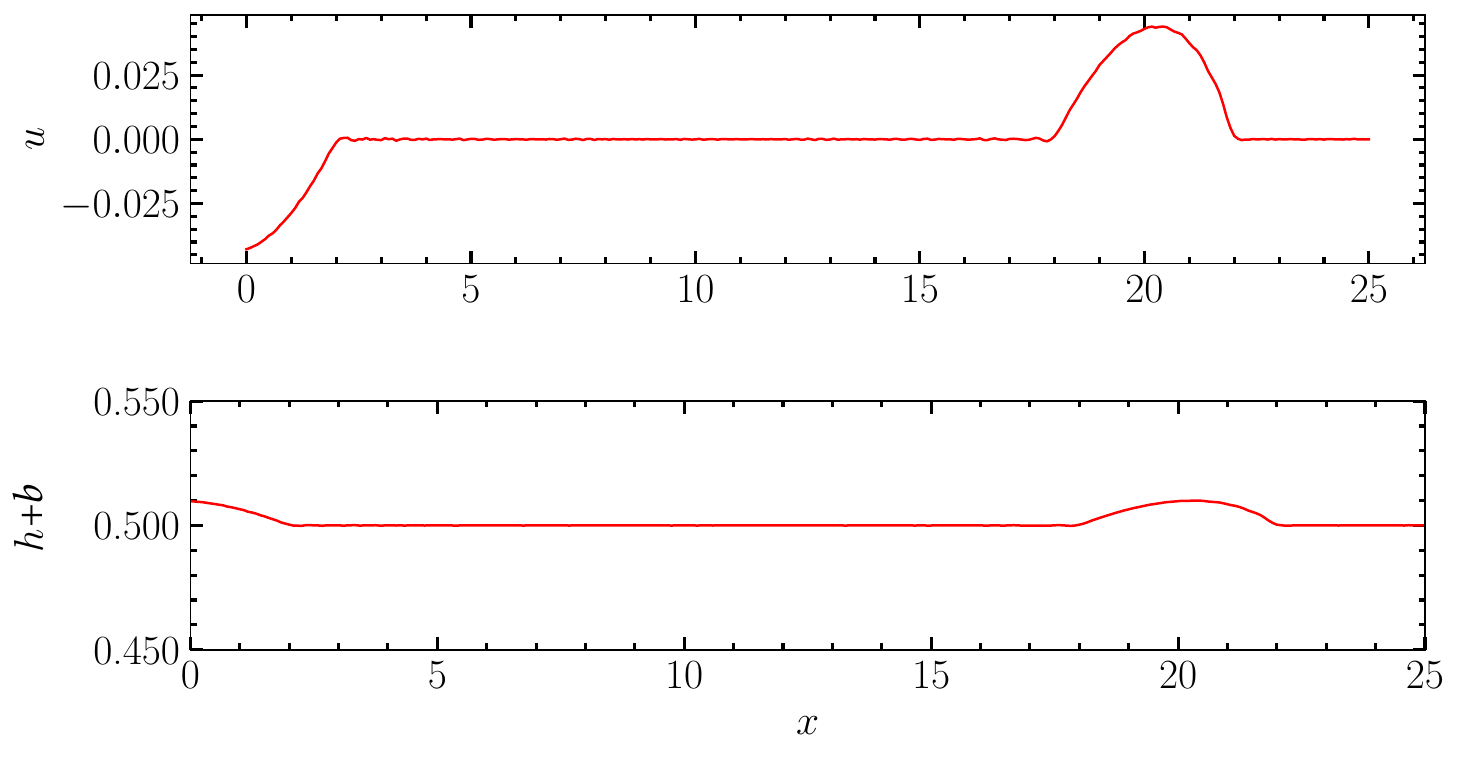}
    \caption{\label{fig:with_hyper_h}$t = 400\Delta t $}
    \end{subfigure}
    \begin{subfigure}{0.5
\textwidth}
    \centering
    \includegraphics[width = 1\linewidth]{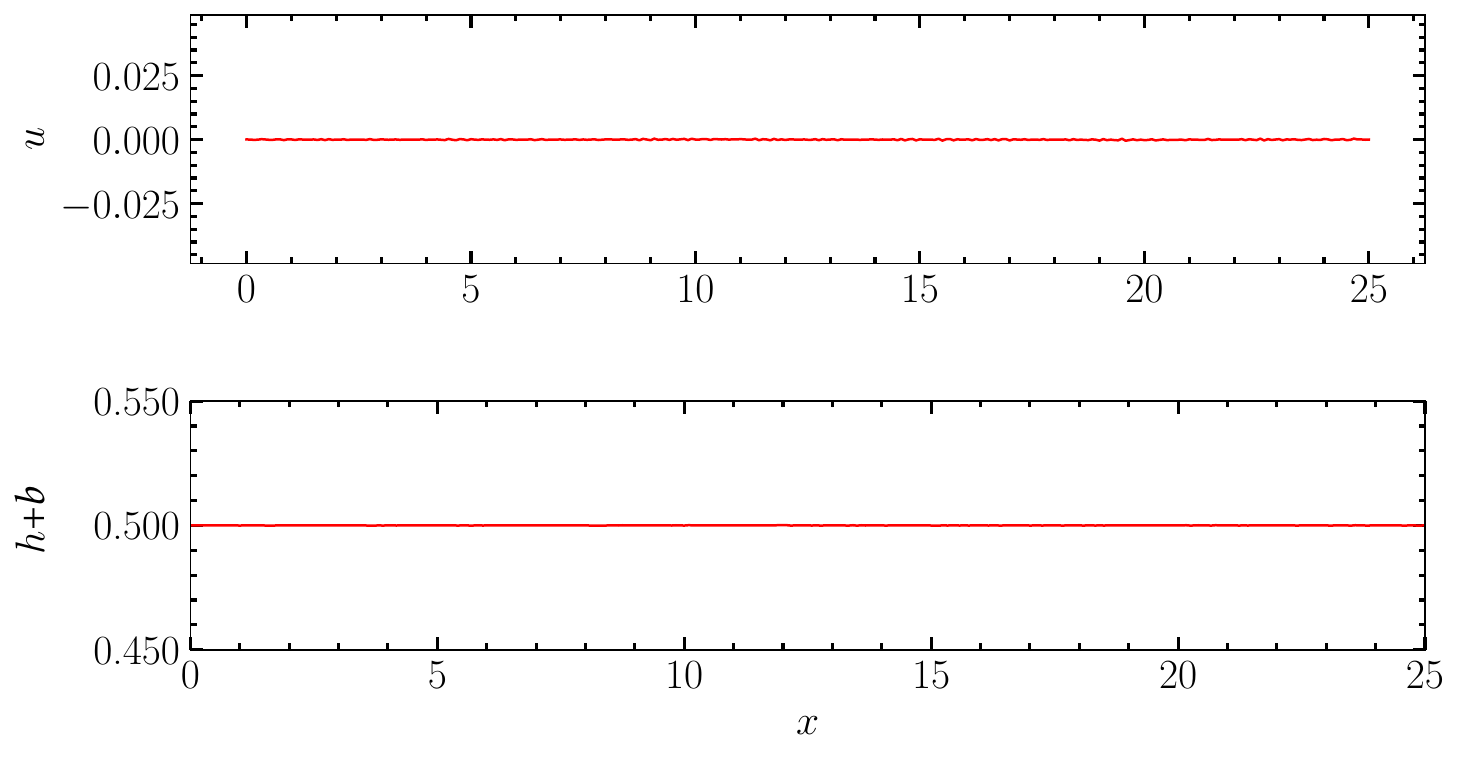}
    \caption{\label{fig:with_hyper_u}$t = 800\Delta t $}
    \end{subfigure}
\caption{\label{fig:lake_at_rest_small} Time evolution of the lake at rest profile with 10\% perturbations on $h$ with transmissive BC. The waves transmit out the right and left boundaries to maintain an overall steady state flowfield similar to the canonical lake at rest problem.}
\end{figure}
\begin{table}
\centering 
\begin{tabular}{c c c c c c c c c  }
		\hline \hline
  \multicolumn{1}{c}{\texttt{FD, order 6, Lake at Rest, smooth perturbations} } \\
    \multicolumn{1}{c}{$N$} &  \multicolumn{2}{c}{(SBP)}  &  \multicolumn{2}{c}{(DP)} & \multicolumn{2}{c}{(DRP)} \\
 & $q_u$       & $q_h$       & $q_u$       & $q_h $     & $q_u  $     & $q_h  $     &  \\ \hline 

		101  & -  &  - & -                             & - & -     \\ 
		201  & 2.0126 & 1.8438 &   2.6958  & 2.8914                           & 2.2873 & 1.6657
                     \\ 
		401  & 2.9720 & 2.9958                            & 3.6705 & 3.5118   & 2.7088 & 2.8301                      \\ 
		801 & 3.9608 & 3.8156 & 5.5908                      & 5.5869 & 3.9935 & 4.1396
                             \\ 
   \hline \hline 
	\end{tabular}
%
\caption{\label{table:lake_smooth}Convergence rate \revi{$q_u$ and $q_h$ of velocity and height} for the lake at rest problem with small perturbations on $h$ ($10\%$), with smooth bathymetry of the Gaussian form, see \eqref{eqn:gaussianform}. \revi{This has been computed at the end time of the simulation, $t = 800\Delta t$ using the analytical solution of the traditional lake at rest problem}, satisfactory convergence rates are achieved with smooth $b(x)$.}
\end{table}
%
%
%
In Fig.~\ref{fig:lake_at_rest_small} we show the time evolution of the perturbed velocity and height profiles. Clearly the perturbations leave the domain without reflections, and we recover the steady state solutions.
We have also performed grid convergence studies, see Table~\ref{table:lake_smooth} for the convergence rates of the  numerical error. The errors converge  to zero at an optimal rate for order 6 FD operators.

\subsubsection{Dam break with wet domain}
Next, we investigate the efficacy of the method in the presence of nonlinear shocks.
We consider the canonical dam break problem \cite{delestre2013swashes} in a wet domain with the following initial conditions:
\begin{equation}
    h(x,t=0)= \begin{cases}h_l & \text { if } 0  \leq x \leq x_0, \\ h_r & \text { if } x_0<x \leq L, \end{cases}
    \qquad u(x,t=0) = 0,
\end{equation}
where we use $x_0 = 5$, $h_l =1$, $h_r = 0.5$, $L =10$ and $Fr \approx 0.3458$, so the flow is subcritical. Note that the initial condition for the water height is discontinuous at $x=x_0$. The exact solutions consists of a right-going  shock front and a left-going rarefaction fan. Using the method of characteristics the  exact solution has the closed form expression   \cite{delestre2013swashes}
\begin{equation*}
\begin{aligned}
    h(t, x)=\left\{\begin{array}{l}h_l \\ \frac{4}{9 g}\left(\sqrt{g h_l}-\frac{x-x_0}{2 t}\right)^2 \\ \frac{c_m^2}{g} \\ h_r\end{array}\right., \; u(t, x)= \begin{cases}0  & \text { if } x \leq x_A(t), \\ \frac{2}{3}\left(\frac{x-x_0}{t}+\sqrt{g h_l}\right) & \text { if } x_A(t) \leq x \leq x_B(t), \\ 2\left(\sqrt{g h_l}-c_m\right) & \text { if } x_B(t) \leq x \leq x_C(t), \\ 0  & \text { if } x_C(t) \leq x,\end{cases}
\end{aligned}
\end{equation*}
where subscripts $l$ and $r$ denote left and right, and $x_A$, $x_B$, $x_C$ satisfy the following evolution equations, with $c_m$ defined as a solution to the algebraic equation,
$$
-8 g h_r c_m^2\left(\sqrt{g h_l}-c_m\right)^2+\left(c_m{ }^2-g h_r\right)^2\left(c_m{ }^2+g h_r\right)=0 
$$
and the coordinates, $x_A, x_B, x_C$ are defined by
\begin{equation}
\begin{aligned}
 x_A(t)=x_0-t \sqrt{g h_l}, \quad x_B(t)=x_0+t\left(2 \sqrt{g h_l}-3 c_m\right) \text { and } \, x_C(t)=x_0+t \frac{2 c_m^2\left(\sqrt{g h_l}-c_m\right)}{c_m^2-g h_r}.
\end{aligned}
\end{equation}

\begin{figure}
    \begin{subfigure}{0.5\textwidth}
    \centering
    \includegraphics[ width = 1\linewidth]{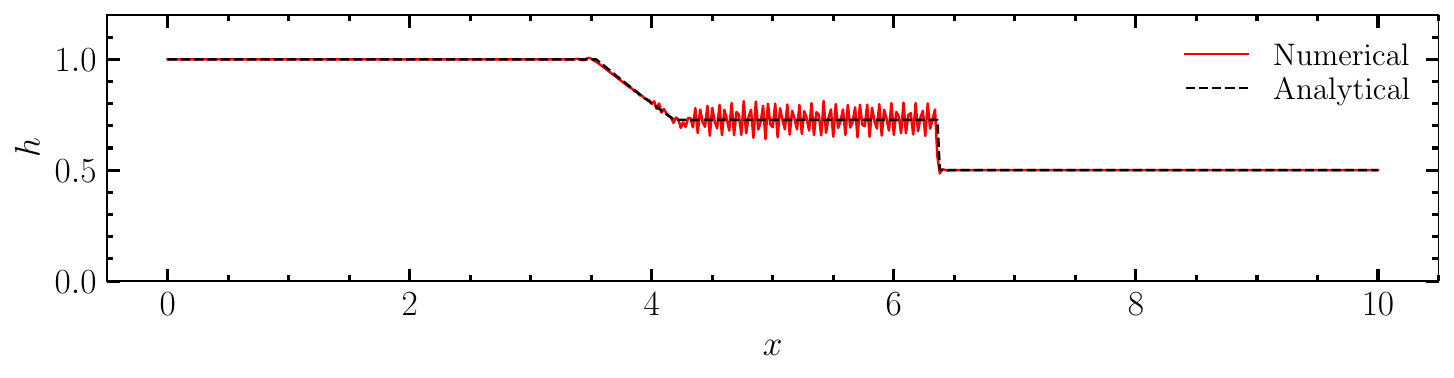}
    \caption{\label{fig:order41}Without hyper-viscosity}
    \end{subfigure}
    \begin{subfigure}{0.5
\textwidth}
    \centering
    \includegraphics[width = 1\linewidth]{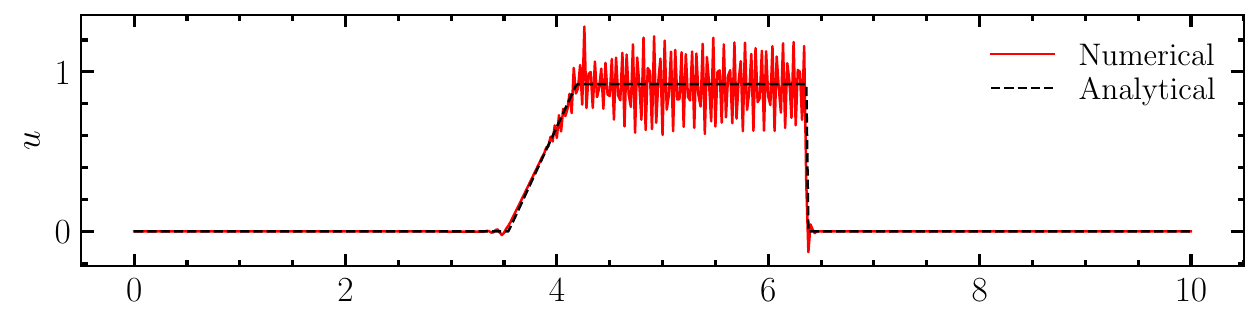}
    \caption{\label{fig:sym2}Without hyper-viscosity}
    \end{subfigure}
        \begin{subfigure}{0.5\textwidth}
    \centering
    \includegraphics[ width = 1\linewidth]{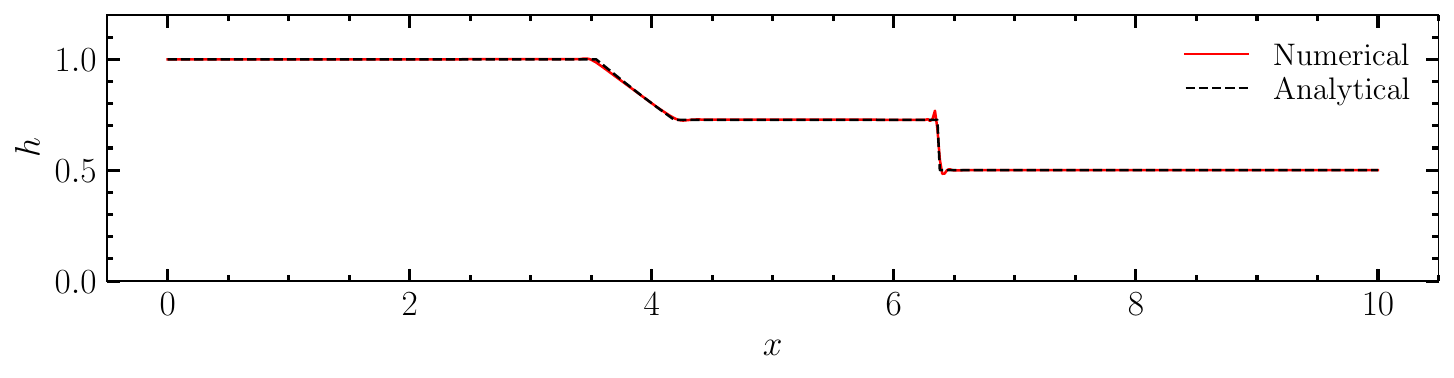}
    \caption{\label{fig:order41}With hyper-viscosity, $\delta = 0.1$}
    \end{subfigure}
    \begin{subfigure}{0.5
\textwidth}
    \centering
    \includegraphics[width = 1\linewidth]{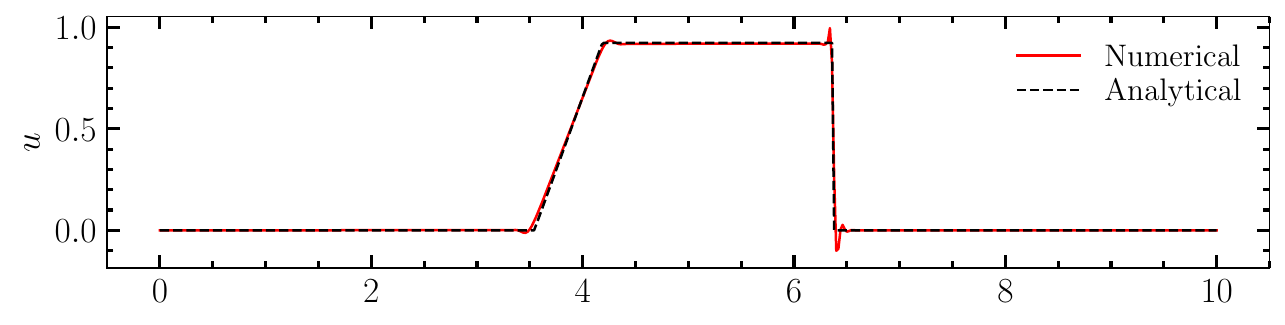}
    \caption{\label{fig:sym2}With hyper-viscosity, $\delta = 0.1$}
    \end{subfigure}
\caption{\label{fig:height_dam_break_hyper}Time evolution of the height and velocity profile of the dam break with wet domain problem, with the overlaid analytical solution using the DP order 6 operators, $N = 501$, showing time snapshot at $t = 1000\Delta t$, comparing the solutions with and without hyper-viscosity. Clearly the oscillations are heavily amplified in those without.  }
\end{figure}

To perform numerical simulations we close the boundaries with transmissive BCs and discretise the computational domain with $N=1001$ grid points. We run the simulation until the final time $t = 2$ using the $6$th order accurate DP-SBP operators, with hyper-viscosity $\delta = 0.1$ and without hyper-viscosity $\delta = 0$. In Figure \ref{fig:height_dam_break_hyper} we compare numerical solutions with the exact solutions at $t = 1000\Delta{t}$. Note that the numerical solution is stable and the shock speed and the rarefaction fan are well resolved by our high order numerical method. However, without hyper-viscosity there are spurious oscillations from the shock front which pollute the numerical solutions. The addition of the hyper-viscosity, with $\delta = 0.1$, eliminates the oscillations without destroying the high order accuracy of the solution in smooth regions. Furthermore, in 
Figs.~\ref{fig:height_dam_break} and \ref{fig:velocity_dam_break} 
we show the time evolution of the height and velocity, with hyper-viscosity. It is clearly demonstrated that the numerical and analytical solutions agree well, despite the presence of shocks and the highly nonlinear nature of the flow problem.
%
\begin{figure}[H]
    \begin{subfigure}{0.5\textwidth}
    \centering
    \includegraphics[ width = 1\linewidth]{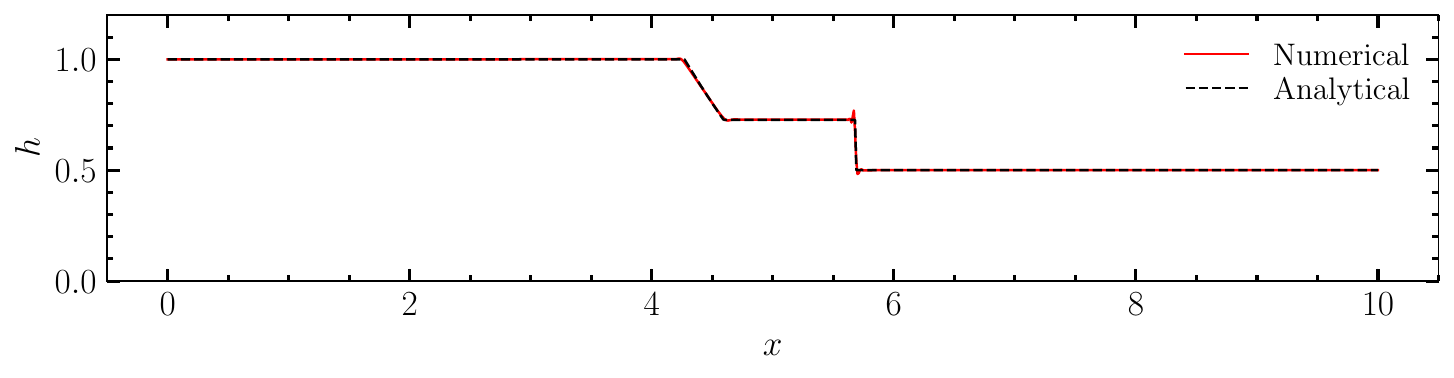}
    \caption{\label{fig:order41}$t = 1000\Delta t $}
    \end{subfigure}
    \begin{subfigure}{0.5
\textwidth}
    \centering
    \includegraphics[width = 1\linewidth]{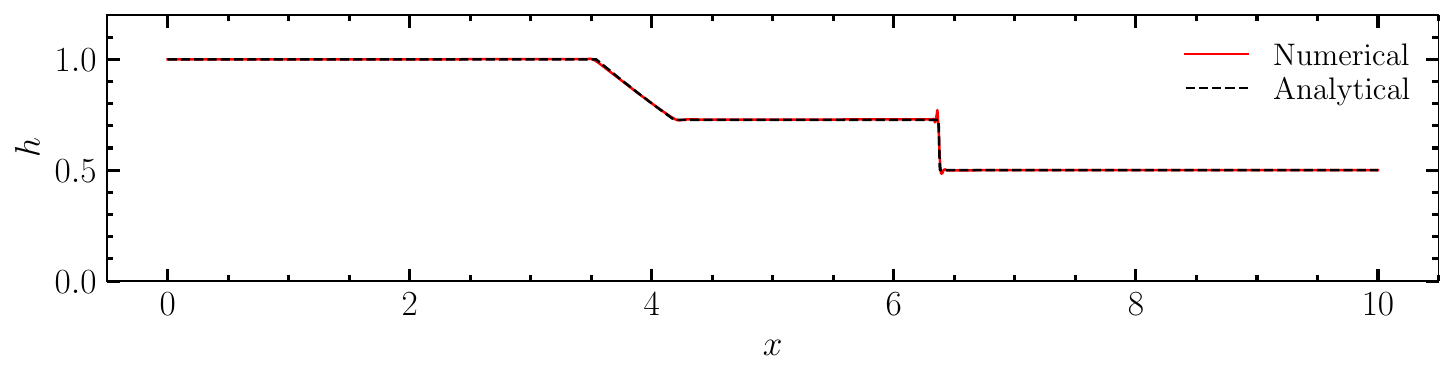}
    \caption{\label{fig:sym2}$t = 2000\Delta t $}
    \end{subfigure}
        \begin{subfigure}{0.5\textwidth}
    \centering
    \includegraphics[ width = 1\linewidth]{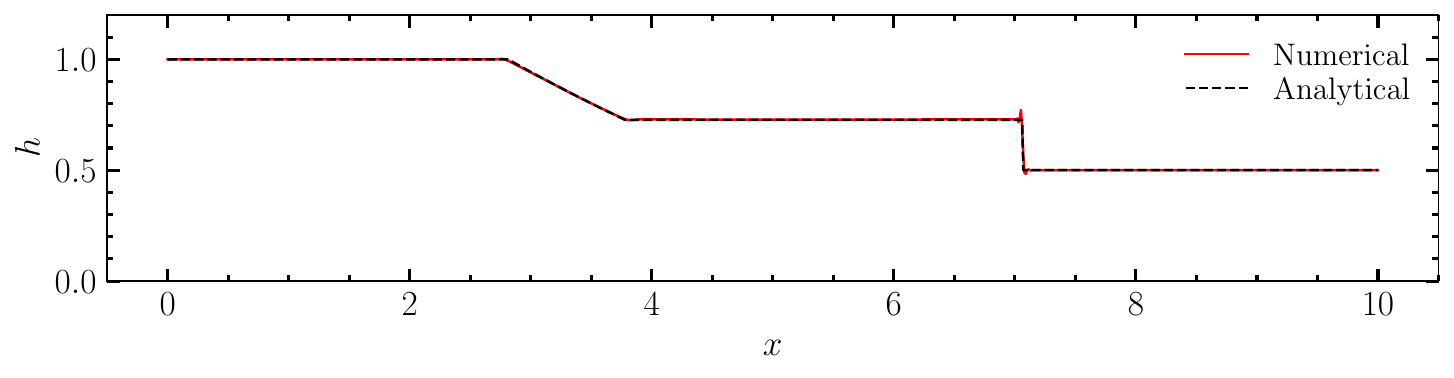}
    \caption{\label{fig:order41}$t = 3000\Delta t $}
    \end{subfigure}
    \begin{subfigure}{0.5
\textwidth}
    \centering
    \includegraphics[width = 1\linewidth]{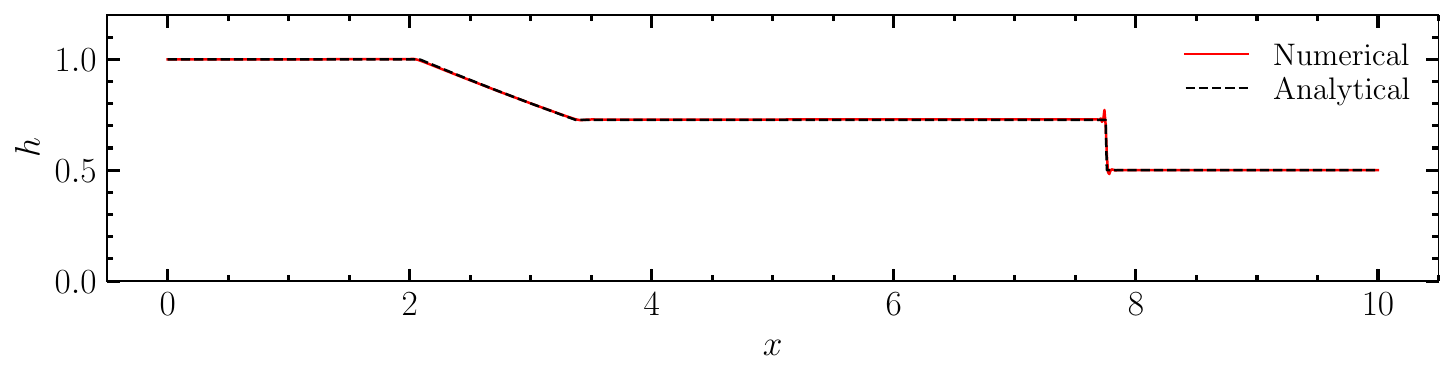}
    \caption{\label{fig:sym2}$t = 4000\Delta t $}
    \end{subfigure}
\caption{\label{fig:height_dam_break}Time evolution of the height profile of the dam break with wet domain problem, with the overlaid analytical solution using the DRP order 6 operators, $N = 1001$, with final time $t = 1$. The hyper-viscosity parameter is set to $\delta = 0.1$ }
\end{figure}
\begin{figure}[H]
    \begin{subfigure}{0.5\textwidth}
    \centering
    \includegraphics[ width = 1\linewidth]{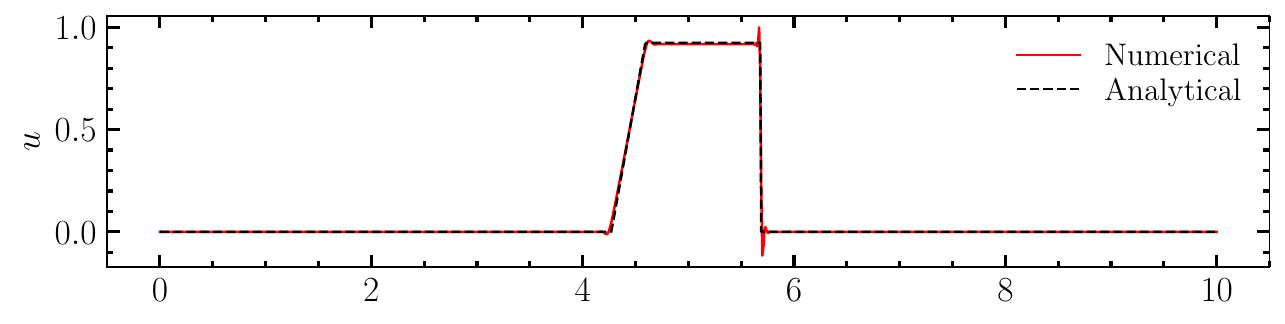}
    \caption{\label{fig:no_hyper_h}$t = 1000\Delta t $}
    \end{subfigure}
    \begin{subfigure}{0.5
\textwidth}
    \centering
    \includegraphics[width = 1\linewidth]{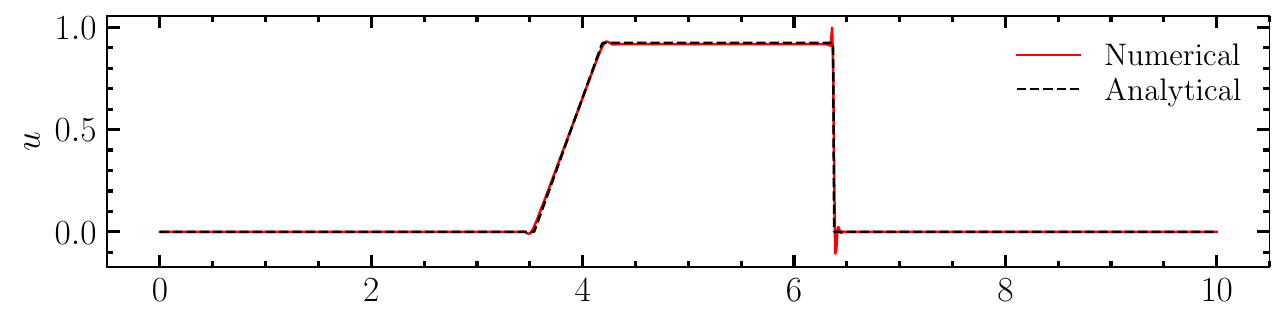}
    \caption{\label{fig:no_hyper_u}$t = 2000\Delta t $}
    \end{subfigure}
        \begin{subfigure}{0.5\textwidth}
    \centering
    \includegraphics[ width = 1\linewidth]{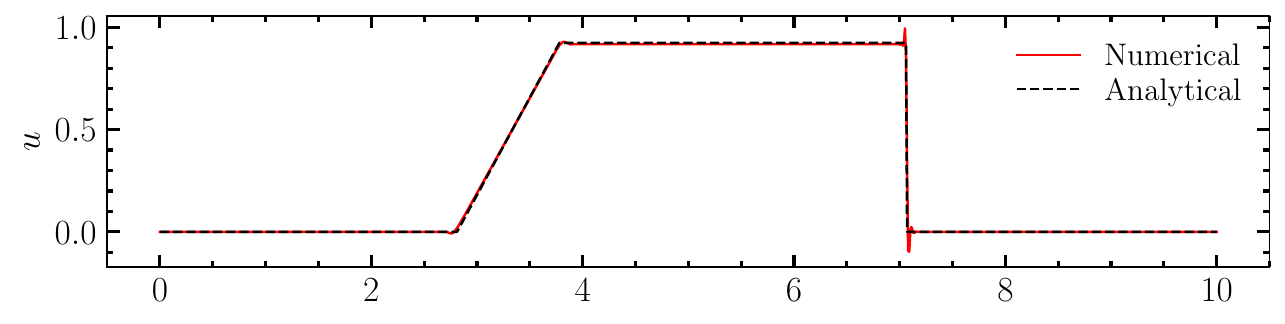}
    \caption{\label{fig:with_hyper_h}$t = 3000\Delta t $}
    \end{subfigure}
    \begin{subfigure}{0.5
\textwidth}
    \centering
    \includegraphics[width = 1\linewidth]{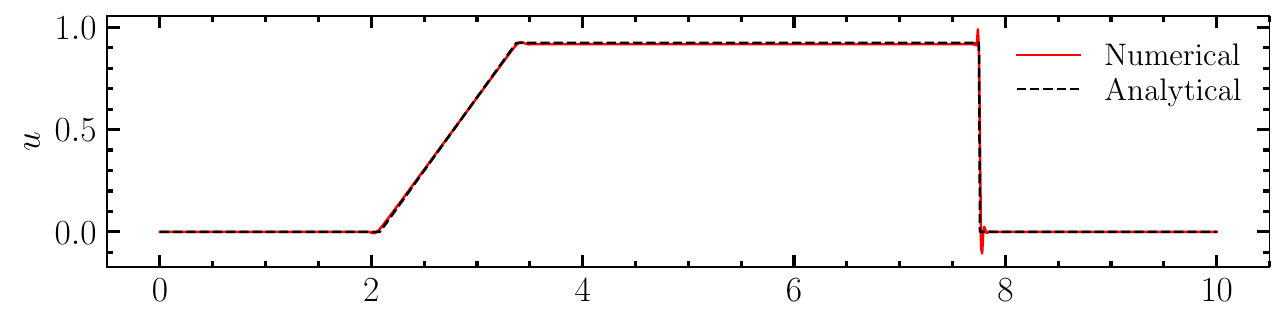}
    \caption{\label{fig:with_hyper_u}$t = 4000\Delta t $}
    \end{subfigure}
\caption{\label{fig:velocity_dam_break} Same as Fig.~\ref{fig:height_dam_break} but for the time evolution of the velocity profile.}
\end{figure}

\subsection{Numerical experiments in 2D}
Here, we now consider numerical experiments in 2D. First we extend the numerical method to be applicable to 2D (see details in ~\ref{sec:hyper2d}) with periodic BCs. We again verify our 2D analysis through a series of numerical experiments such as the method of manufactured solutions, a merging vortex test case which maintains geostrophic balance, and the classical barotropic instability problem. The numerical experiments show accuracy, stability and convergence  properties of our robust numerical scheme.
\subsubsection{Method of manufactured solutions}

Here we verify accuracy for 2D simulations using forced solutions in a doubly periodic domain. For the 2D rotating SWE \eqref{eqn:nlSWE2D}, we let it satisfy the exact sinusoidal solutions of the form,
\begin{equation}\label{eqn:mms}
\begin{aligned}
h(x, y, t) &= H + 0.1 \cos(k_x (x - x_0)) \cos(k_y  (y - y_0))  \cos(\omega  t), \\
u(x, y, t) &= \cos(\omega  t)  \sin(k_y  (y - y_0))  \cos(k_x (x - x_0)), \\
v(x, y, t) &= \sin(\omega  t) \cos(k_y(y-y_0)) \sin(k_x(x-x_0)),
\end{aligned}
\end{equation}
where $x_0 = y_0 = 0.5$ and $(x, y) \in \Omega = [0, 1]^2$. Here $g = 9.81$, $H = 10$, $f=0$, with $k_x = k_y = 2\pi$ and $\omega = 2\pi$. We measure the discrete $l_2$ error and evaluate the convergence rate as, $q_{\alpha,i} = \log(\alpha_{i-1}/\alpha_{i})/\log({m_i/m_{i-1}})$, where $\alpha_i$ is the $l_2$-error at grid level $i$, and $m_i$ is the number of cells in $x$ or $y$ directions. 
\begin{figure}

     \begin{subfigure}{0.48
\textwidth}
    \centering
    \includegraphics[width = \linewidth]{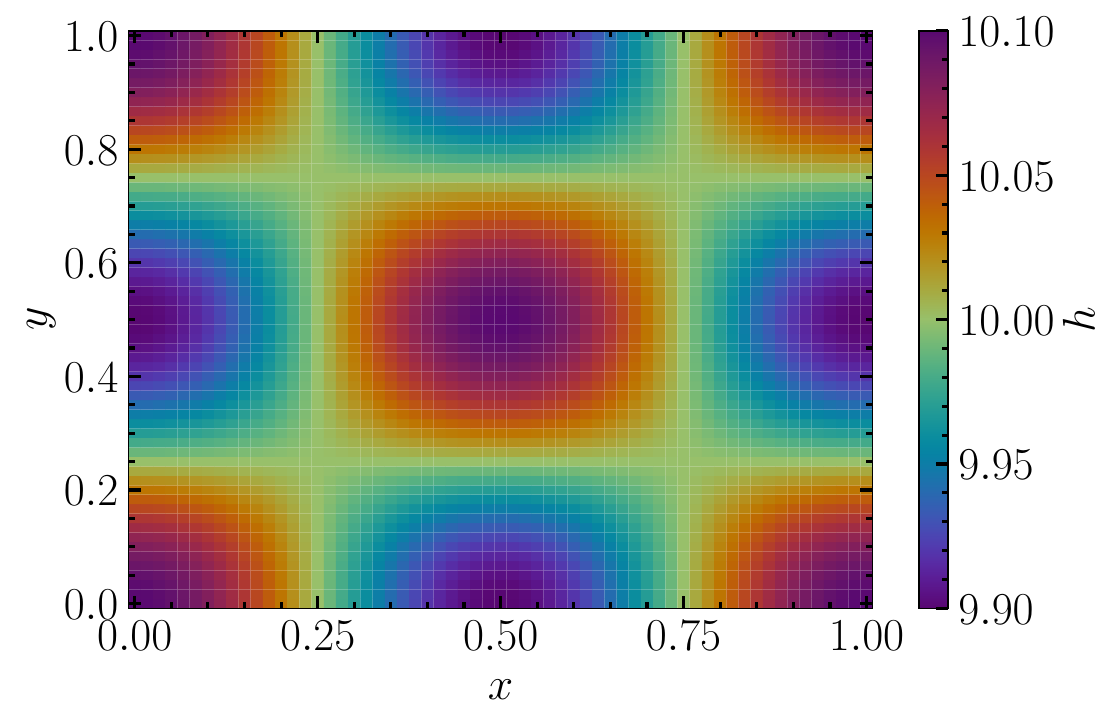}
    \caption{}
    \label{fig:mms}
    \end{subfigure}
      \begin{subfigure}{
0.48\textwidth}
    \centering
    \includegraphics[width = \linewidth]{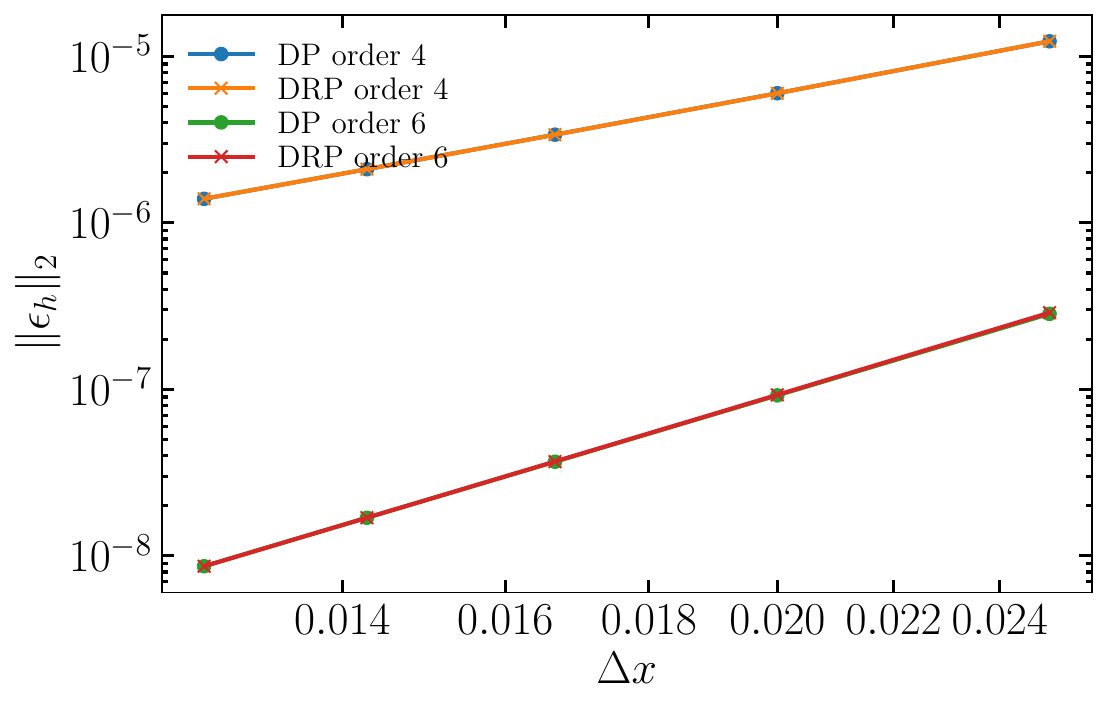}
    \caption{}
    \label{fig:mms_error}
    \end{subfigure}
    \caption{ (a) Contour of $h$ for the MMS solution \eqref{eqn:mms} at $t = 0$. (b) $l_2$ error of height as a function grid size, $\Delta x$. }
\end{figure}
\begin{table}
\centering
\begin{tabular}{|c|c|c|c|c|}
\hline 
$m$ & $q_{{h}_{\text{DP4}}}$ & $q_{{h}_{\text{DP6}}}$ & $q_{{h}_{\text{DRP4}}}$ &$ q_{{h}_{\text{DRP6}}} $\\ \hline \hline
31  &                &                &  &  \\ 
41  & 3.0340               & 5.0207               & 3.0322 & 5.0582 \\ 
51  & 3.0475               & 5.0253               & 3.0449  & 5.0905  \\ 
61  & 3.0342               & 5.0210               & 3.0323  & 5.0569 \\
71  & 3.0259               & 5.0179               & 3.0247 & 5.0389   \\
81  & 3.0204               & 5.0155               &  3.0196 &  5.0280 \\ \hline
\end{tabular}
\caption{\label{tab:1}Convergence rate, \revi{$q_h$} of height as a function of spatial resolution for DP and DRP operators of order 4 and 6, computed from the manufactured solution, with grid-dependent hyper-viscosity scaling as $\sim (\Delta x)^3$ for order 4 and $\sim (\Delta x)^5$ for order 6, and $\delta = 0.1$. Note that the errors converge at the expected rate.}
\end{table}
We consider the convergence of depth ($h$) at final time $t_{\text{end}} = 0.5$ in Table~\ref{tab:1}, which all indicate the expected convergence rates when incorporating grid-dependent hyper-viscosity, which scales as $(\Delta x)^3$ for order 4, and $(\Delta x)^5$ for order 6.

\subsubsection{2D merging vortex problem}
Next we consider the ``merging vortex problem" \cite{mcrae2014energy}. This problem is modelled by the 2D rotating SWE \eqref{eqn:nlSWE2D} 
in the spatial domain $\Omega=(0,2 \pi] \times(0,2 \pi]$ with periodic BCs. The initial conditions are a pair of Gaussian vortex  with the in-compressible stream function
\begin{equation}\label{eqn:psi}
\psi=e^{-5\left((y-\pi)^2+(x-2.6 \pi / 3)^2\right)}+e^{-5\left((y-\pi)^2+(x-3.5 \pi / 3)^2\right)}. 
\end{equation}
 The initial conditions for the velocity field are defined through $\mathbf{u}(x,y,0) = \nabla^{\perp} \psi $, $\nabla^{\perp} = (-\partial_y, \partial_x)$ and the initial condition for water height $h$ is obtained from linear geostrophic balance $f \times \mathbf{u}^{\perp}(x,y,0) + g\times\nabla h(x,y,0)=0 \implies h(x,y,0) = H + f/g\times\psi(x,y)$, with $f = g = H = 8 $. 

 The continuous 2D rotating SWEs \eqref{eqn:nlSWE2D} preserves infinitely many invariants. However, any discrete approximation can (approximately) preserve a finite number of invariants \cite{ricardo2023conservation}. For the numerical approximation we select  a subset of the invariants, namely: 1) total energy/entropy $E(t)$, 2) total enstrophy $E_s(t)$, 3) total vorticity $W(t)$ and  4) total mass $M(t)$,  defined by
\begin{equation}\label{eq:continuous_invariants}
\begin{aligned}
    E(t) = \int_{\Omega} e dx dy, \quad {E}_{\mathrm{s}}(t) = \int_{\Omega}   h^{-1} \omega^2  dx dy, \quad {W}(t) = \int_{\Omega} \omega dx dy,  \quad M(t) = \int_{\Omega} h dx dy,
\end{aligned}
\end{equation}
which the numerical method should accurately preserve.
Here, $e = \frac{1}{2}(g h^2 + h u^2 + h v^2  )$
is the elemental energy/entropy and $\omega$ is the absolute vorticity defined in \eqref{eqn:nlSWE2D}. The energy/entropy $E(t)$ is critical for nonlinear stability. The enstrophy ${E}_{\mathrm{s}}$ is a higher order moment and weakly bounds the derivatives of the solution. More importantly, the enstrophy is crucial for controlling grid-scale errors and ensuring that high frequency oscillations do not dominate the accuracy of the solution.

The numerical approximations  of the invariants are given by
\begin{equation}\label{eq:discrete_invariants}
\begin{aligned}
    \mathcal{E}(t) = \sum^{N}_{i,j=1} e_{ij} \Delta x \Delta y, \quad \mathcal{E}_{\mathrm{s}}(t) = \sum^{N}_{i,j=1} h_{ij}^{-1} \omega^2_{ij} \Delta x \Delta y, \quad \mathcal{W}(t) = \sum^{N}_{i,j=1} \omega_{ij} \Delta x \Delta y,  \quad 
    \mathcal{M}(t) = \sum^{N}_{i,j=1} h_{ij} \Delta x \Delta y,
\end{aligned}
\end{equation}
where $e_{ij}$ and $\boldsymbol{\omega}$ are defined in \eqref{eqn:nl2dSWE2D} and~\eqref{eqn:elemental_2d}, respectively, and  where $i,j$ are indices corresponding to grid points. 
We define the relative changes in the discrete invariants
\begin{equation}\label{eq:change_invariants}
\begin{aligned}
    \Delta_r{\mathcal{E}}(t) = \frac{\mathcal{E}(t)-\mathcal{E}(0)}{\mathcal{E}(0)}, \quad \Delta_r{\mathcal{E}}_s(t) = \frac{\mathcal{E}_s(t)-\mathcal{E}_s(0)}{\mathcal{E}_s(0)}, \quad \Delta_r{\mathcal{W}}(t) = \frac{\mathcal{W}(t)-\mathcal{W}(0)}{\mathcal{W}(0)}, \quad 
    \Delta_r{\mathcal{M}}(t) = \frac{\mathcal{M}(t)-\mathcal{M}(0)}{\mathcal{M}(0)}.
\end{aligned}
\end{equation}

The semi-discrete approximation of the 2D rotating SWEs \eqref{eqn:nlSWE2D} in periodic domains using the DP SBP framework is derived in \ref{sec:hyper2d}. In Theorem \ref{thm:swp2d_Diss_nonlinear}, we  prove that  the semi-discrete approximation conserves the discrete energy/entropy without hyper-viscosity $\delta =0$ and dissipates  energy/entropy with  hyper-viscosity $\delta >0$. Semi-discrete conservation of total mass and vorticity also follow, but the proofs are not given. 


As before, the semi-discrete approximations are evolved in time using  classical fourth order accurate explicit Runge-Kutta method with $\operatorname{CFL}=0.1$. In Fig.~\ref{fig:2d_experiment} we show snapshots of potential vorticity using DP-SBP operators of order 4 on a uniform grid of size $251^2$ grid points. 
The left panel shows the snapshots of the absolute vorticity with hyper-viscosity ($\delta = 0.5$) and right panel shows the snapshots of the absolute vorticity without hyper-viscosity ($\delta = 0$). It is significantly noteworthy that in both cases the numerical solutions remain bounded through out the simulation. This is consistent with the energy stability proof of Theorem \ref{thm:swp2d_Diss_nonlinear}. However, without hyper-viscosity ($\delta = 0$, right panel), the structure in the solution is completely destroyed by spurious numerical artefacts. The addition of hyper-viscosity to the scheme ($\delta = 0.5$, left panel) eliminates the spurious wave modes which can destroy the accuracy of the numerical solution. With hyper-viscosity  the numerical method preserves approximate geostrophic balance and the solution remains accurate throughout the simulation. The numerical solutions are comparable to the results in past works \cite{mcrae2014energy} using compatible finite element methods. We also note that in \cite{mcrae2014energy}, in order to avoid spurious numerical oscillations the authors utilised the so-called Artificial Potential Vorticity Method (APVM) in their finite element method. In general for nonlinear problems, although energy/entropy conservation may ensure the boundedness of numerical solutions, it does not guarantee convergence of numerical errors. Suitable amount of numerical dissipation is necessary to control high frequency errors.


We also compute the relative change in the discrete invariants \eqref{eq:discrete_invariants} as defined in \eqref{eq:change_invariants}. 
It is important to note that while our semi-discrete approximation (without hyper-viscosity) is energy/entropy conserving, the \revi{nature of the stable Runge-Kutta time-stepping scheme can result in energy decay through numerical dissipation}. Therefore, our fully-discrete scheme \revi{may not necessarily be fully} energy conserving. However, when the scheme is stable the conservation errors are expected to converge to zero with mesh refinement.

The time evolution of relative change in the discrete invariants are shown in Figs.~\ref{fig:error_enstrophy} and~\ref{fig:error_enstrophy_order_6} for DP and DRP SBP operators of order 4 and 6. Without hyper-viscosity, note that total mass and total absolute vorticity are conserved up to machine precision, and energy is dissipated by the Runge-Kutta time stepping scheme. However, without hyper-viscosity, the enstrophy grows linearly leading to the growth of high frequency errors which contaminate the solution. With the addition of hyper-viscosity, the magnitude of the enstrophy is relatively small and decays linearly with time. With hyper-viscosity, while total absolute vorticity is conserved up to machine precision, total energy and total  enstrophy are dissipated, and  total mass is conserved up to discretisation errors. In any case, the conservation errors converge with mesh refinement. The  convergence rates of energy and mass conservation errors are shown in Table~\ref{tab:convdisp_mergingvortex}. The $L_2$ errors are displayed in Fig.~\ref{fig:merging_vortex_conv_error}, for both DP and DRP operators of order 4 and 6, computed at end time $t = 4$. We find nearly third order convergence for both energy and mass in the fourth order cases, and near fourth order convergence for the sixth order ones. \revi{ For enstrophy, which is a higher order moment and does not converge in the usual way as does mass and energy, we find second order convergence for either cases.} We note that the sub-convergence observed is not due to the spatial order but due to time discretisation when running until $t = 4$, since we utilise the explicit fourth order Runge Kutta (RK4) method. Thus, these demonstrably indicate that when introducing artificial dissipation, we achieve satisfactory convergence rates of mass and energy due to the self-consistent properties of the dissipation operator we have derived.
%
%
\begin{figure}
    \begin{subfigure}{0.45
\textwidth}
    \centering
    \includegraphics[width = 1\linewidth]{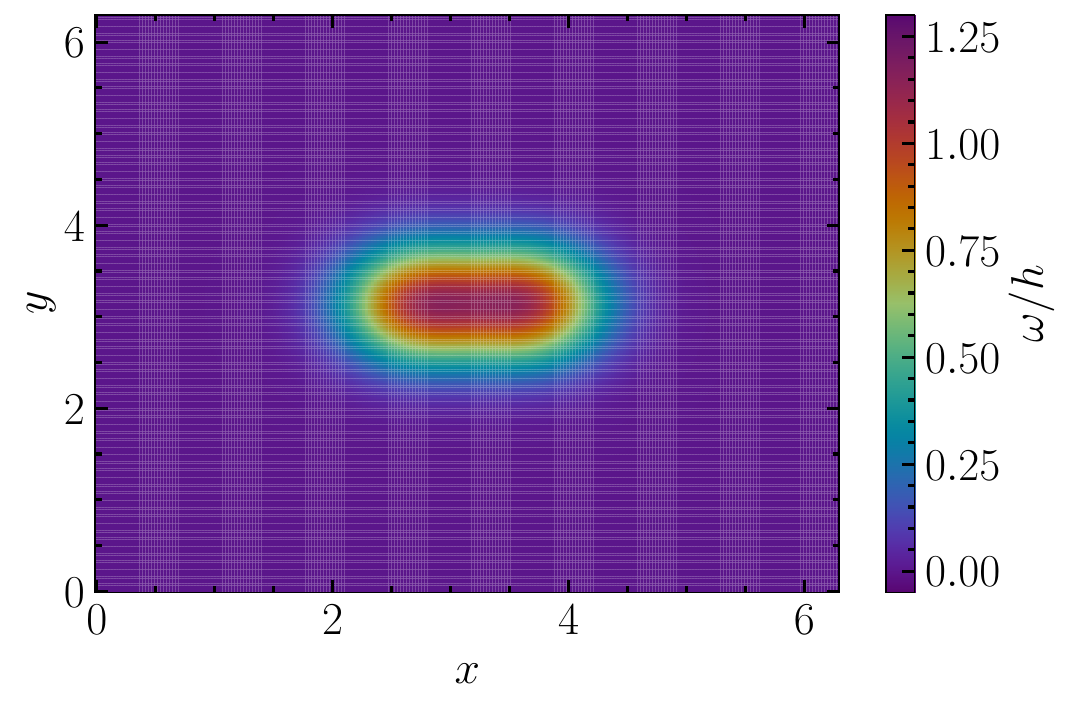}

    \end{subfigure}
    \begin{subfigure}{0.45\textwidth}
    \centering
    \includegraphics[ width = 1\linewidth]{figures/init_vort_new-2.pdf}
    \end{subfigure}

        \begin{subfigure}{0.45
\textwidth}
    \centering
    \includegraphics[width = 1\linewidth]{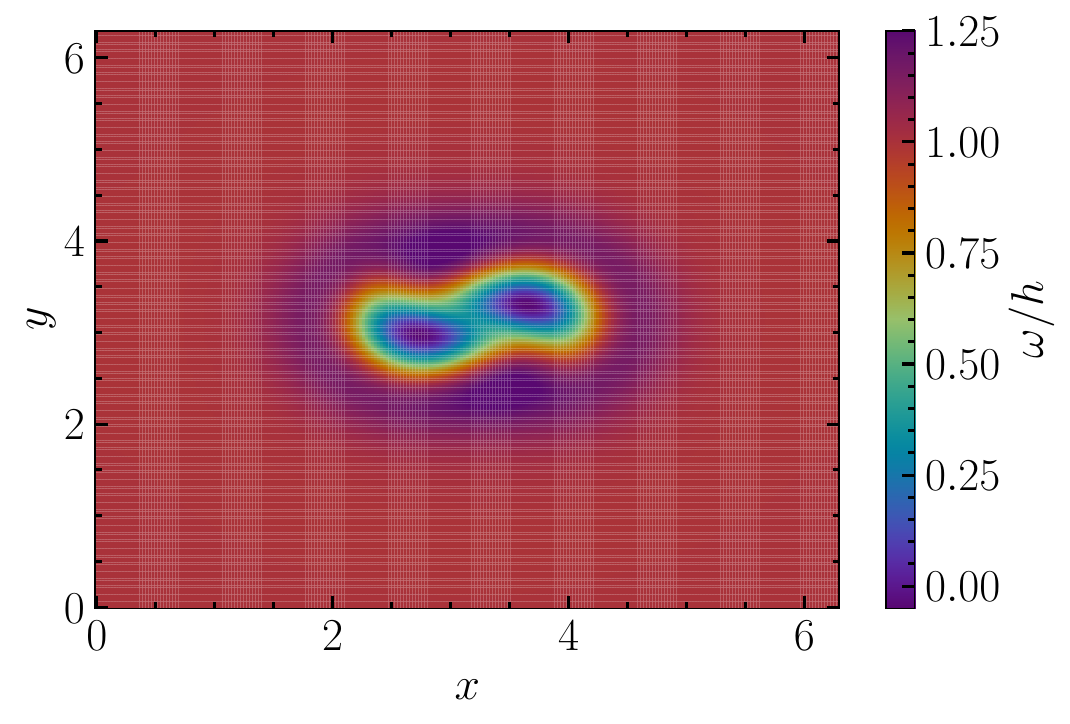}

    \end{subfigure}
        \begin{subfigure}{0.45\textwidth}
    \centering
    \includegraphics[ width = 1\linewidth]{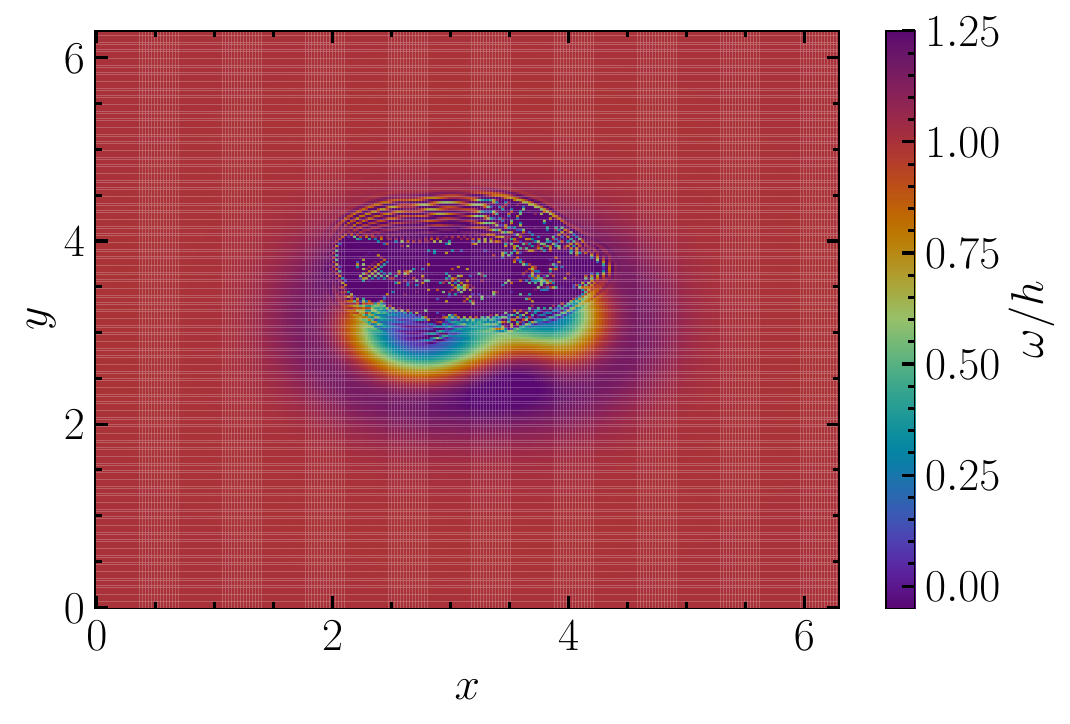}

    \end{subfigure}

                \begin{subfigure}{0.45
\textwidth}
    \centering
    \includegraphics[width = 1\linewidth]{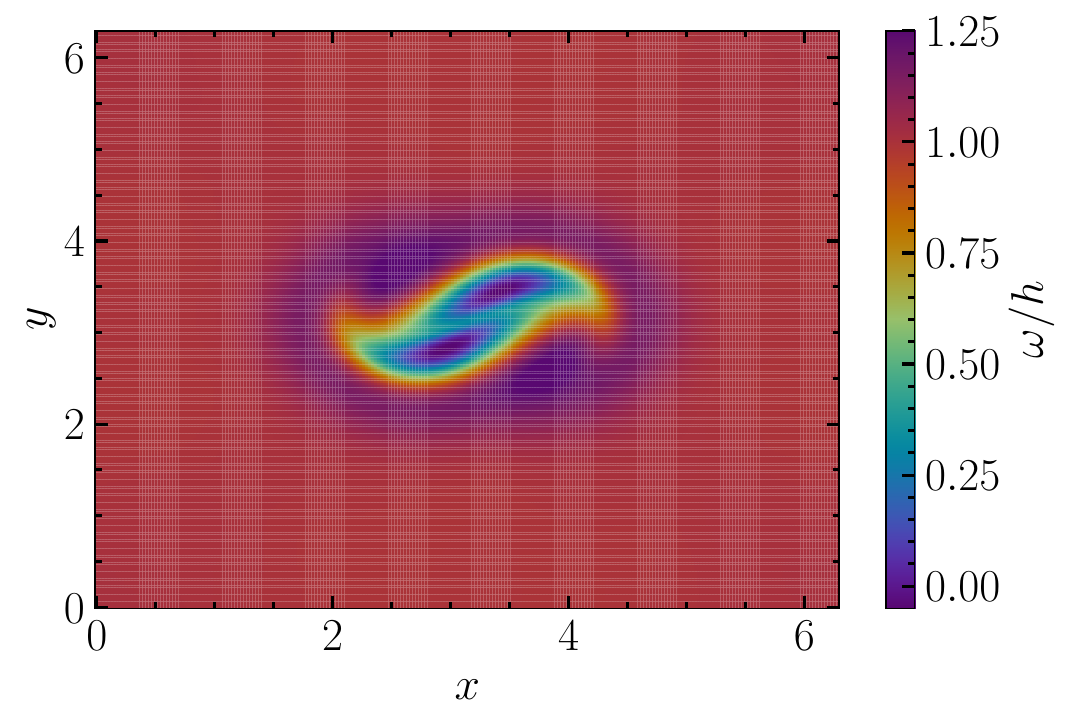}
  
    \end{subfigure}
        \begin{subfigure}{0.45
\textwidth}
    \centering
    \includegraphics[width = 1\linewidth]{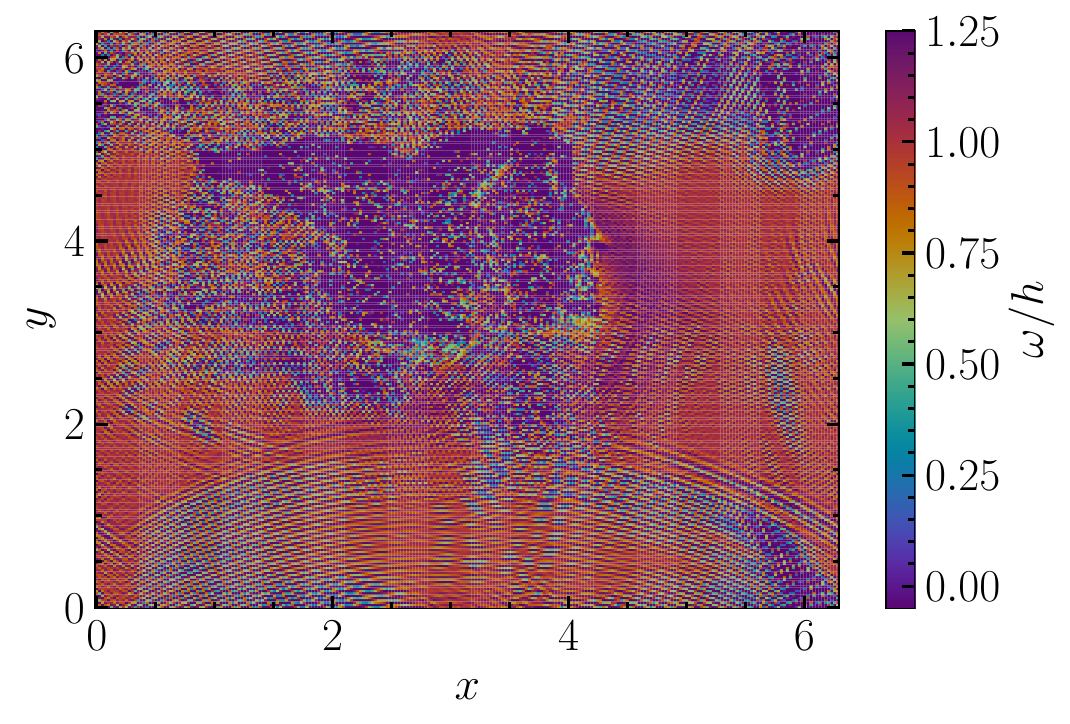}
  
    \end{subfigure}

                        \begin{subfigure}{0.45
\textwidth}
    \centering
    \includegraphics[width = 1\linewidth]{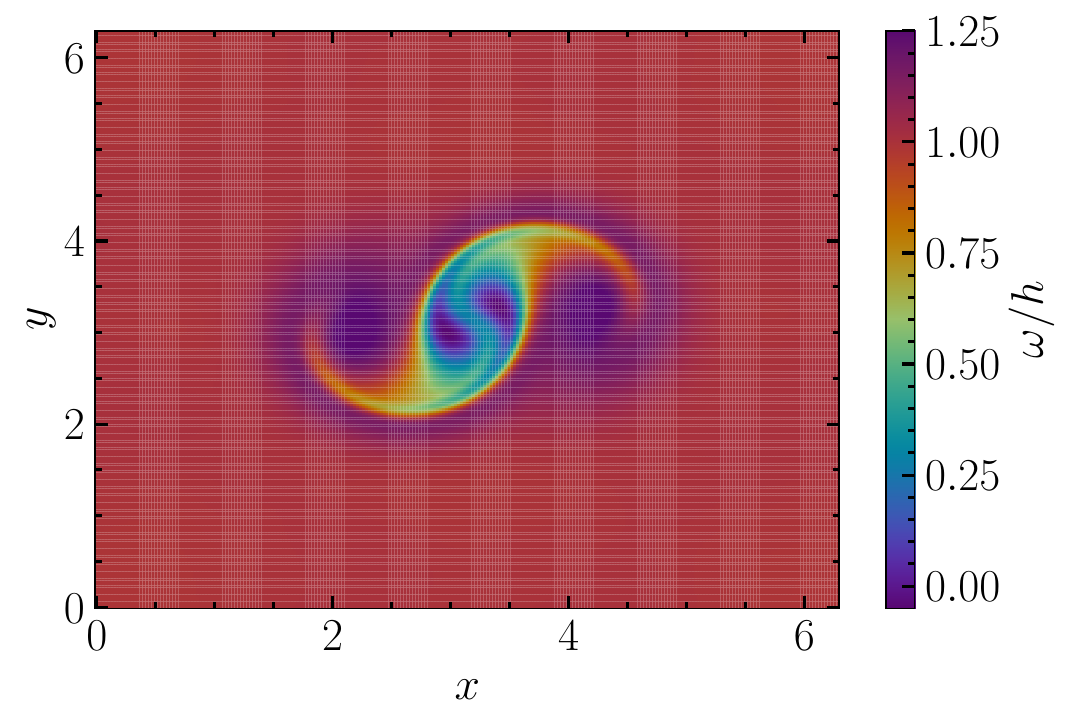}

    \end{subfigure}
                \begin{subfigure}{0.45
\textwidth}
    \centering
    \includegraphics[width = 1\linewidth]{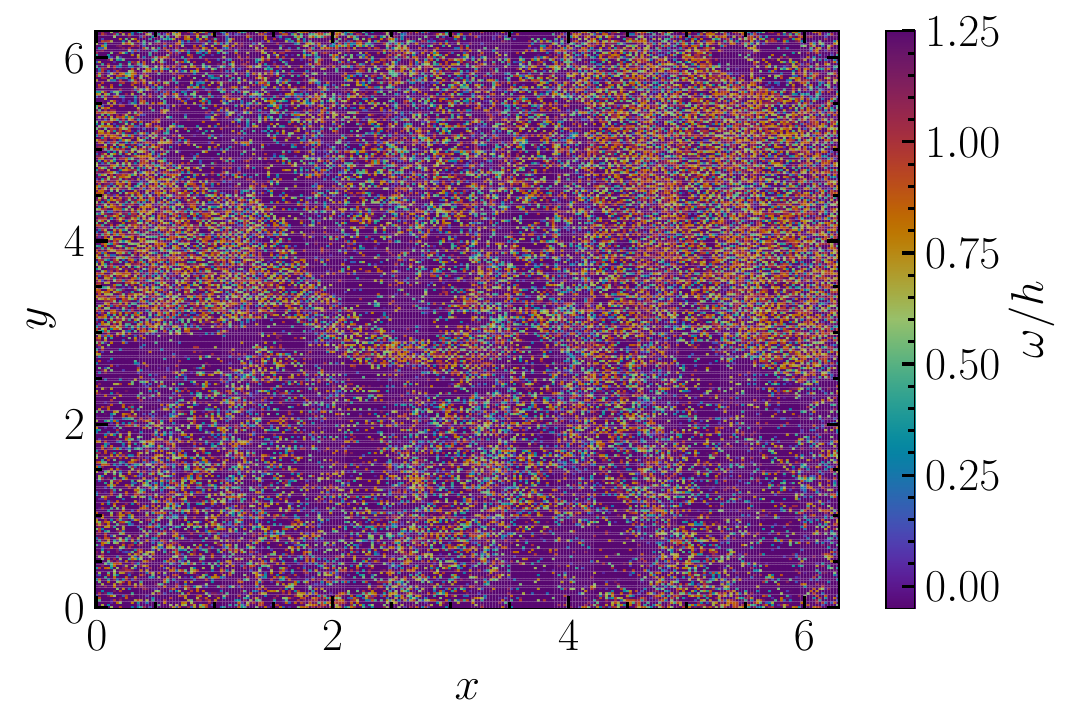}
 
    \end{subfigure}




  
  
 
\caption{\label{fig:2d_experiment} Potential vorticity, $\omega/h$ for the 2D merging vortex problem \eqref{eqn:psi} using DP operators of order 4. The left panels show the solutions with hyper-viscosity, while the right show those without hyper-viscosity, where $\delta = 0.5$. Snapshots are at $t = 0, 0.3, 0.7,1.5$.}
\end{figure}


\begin{figure}
    \centering
    \includegraphics[width = 0.99\linewidth]{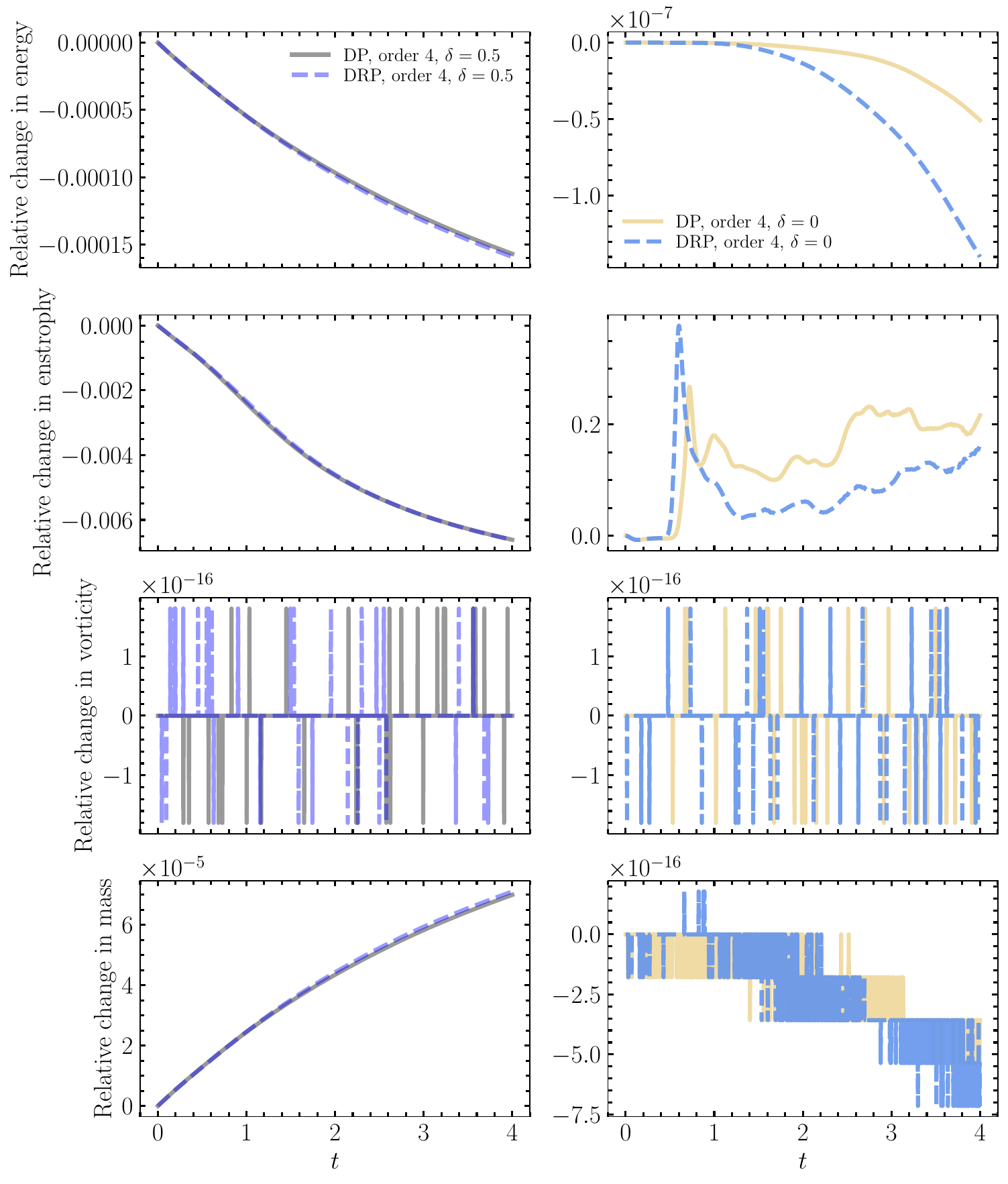}
    \caption{Relative change of the total energy, vorticity, enstrophy and mass for the 2D merging vortex problem using 4th order accurate DP and DRP SBP operators. \revi{The left panel corresponds to  hyperviscosity and the right panel corresponds  to no hyperviscosity.} Clearly with no hyperviscosity, mass and vorticity remains conserved up to machine precision, but enstrophy grows which eventually pollutes the solution everywhere. With  hyperviscosity, vorticity remains conserved up to machine precision, but enstrophy and energy are dissipated.}
    \label{fig:error_enstrophy}
\end{figure}
\begin{figure}
    \centering
    \includegraphics[width = 0.89\linewidth]{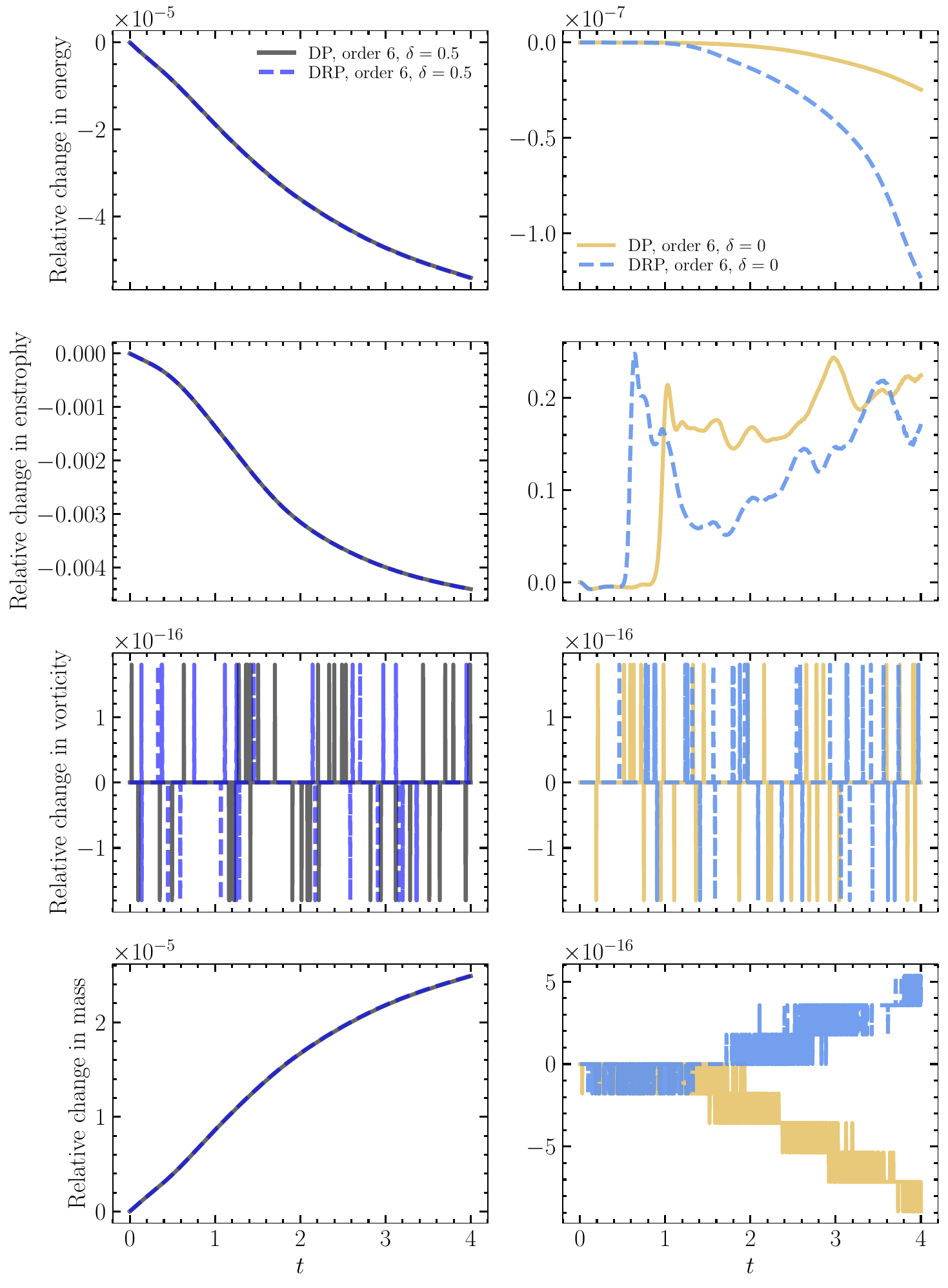}
    \caption{Relative change of the total energy, vorticity, enstrophy and mass for the 2D merging vortex problem using 6th order accurate DP and DRP SBP operators. \revi{The left panel corresponds to  hyperviscosity and the right panel corresponds  to no hyperviscosity.} Clearly with no hyperviscosity, mass and vorticity remains conserved up to machine precision, but enstrophy grows which eventually pollutes the solution everywhere. With  hyperviscosity, vorticity remains conserved up to machine precision, but enstrophy and energy are dissipated.}
    \label{fig:error_enstrophy_order_6}
\end{figure}
\begin{figure}
\centering
     \begin{subfigure}{0.48
\textwidth}
    \centering
    \includegraphics[width = \linewidth]{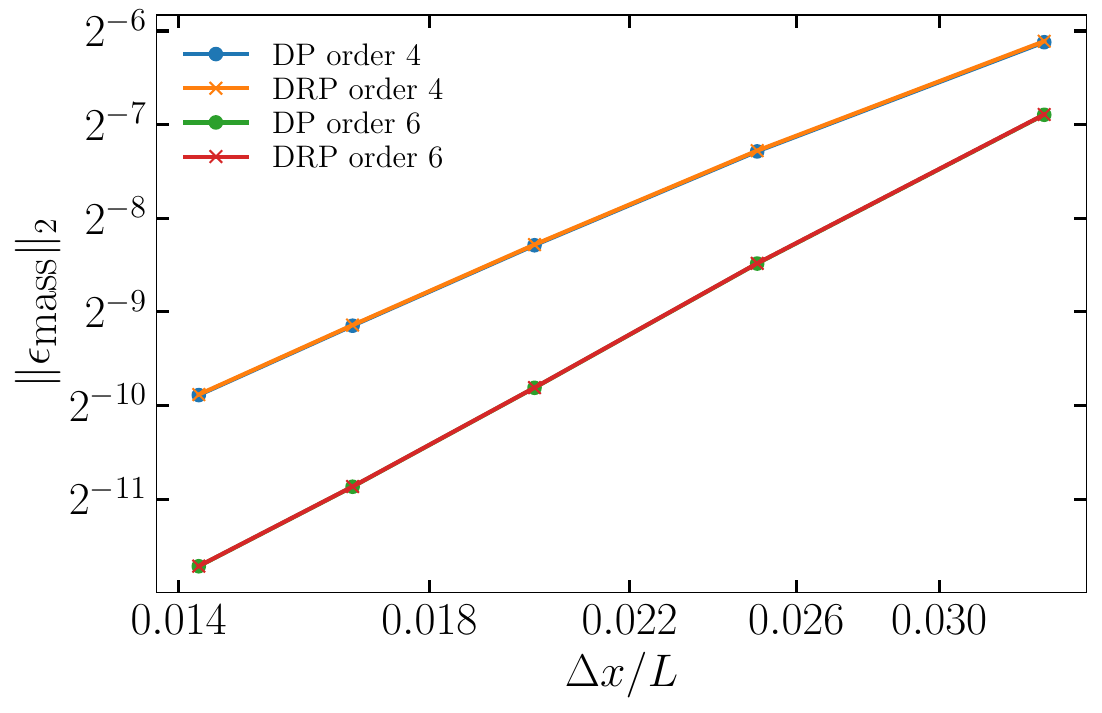}
    \caption{}
    \label{fig:mass_error_merging}
    \end{subfigure}
 \revi{     \begin{subfigure}{
0.48\textwidth}
    \centering
    \includegraphics[width = \linewidth]{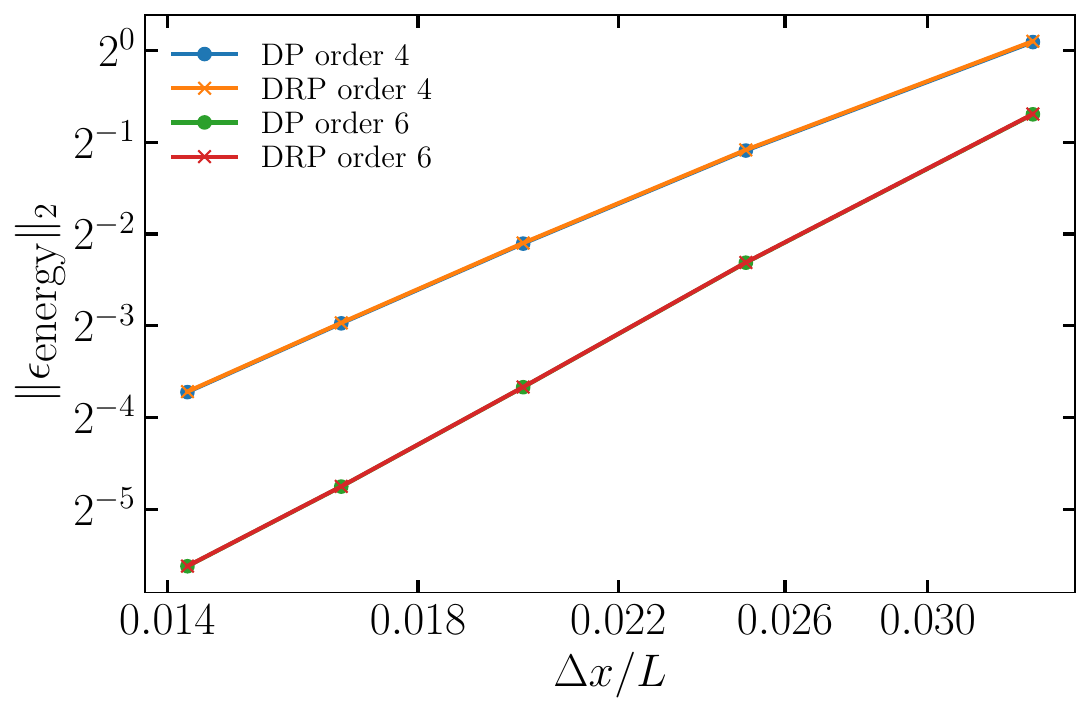}
    \caption{}
    \label{fig:energy_error_merging}
    \end{subfigure}
    \begin{subfigure}{0.5\textwidth}
\centering
    \includegraphics[width=\linewidth]{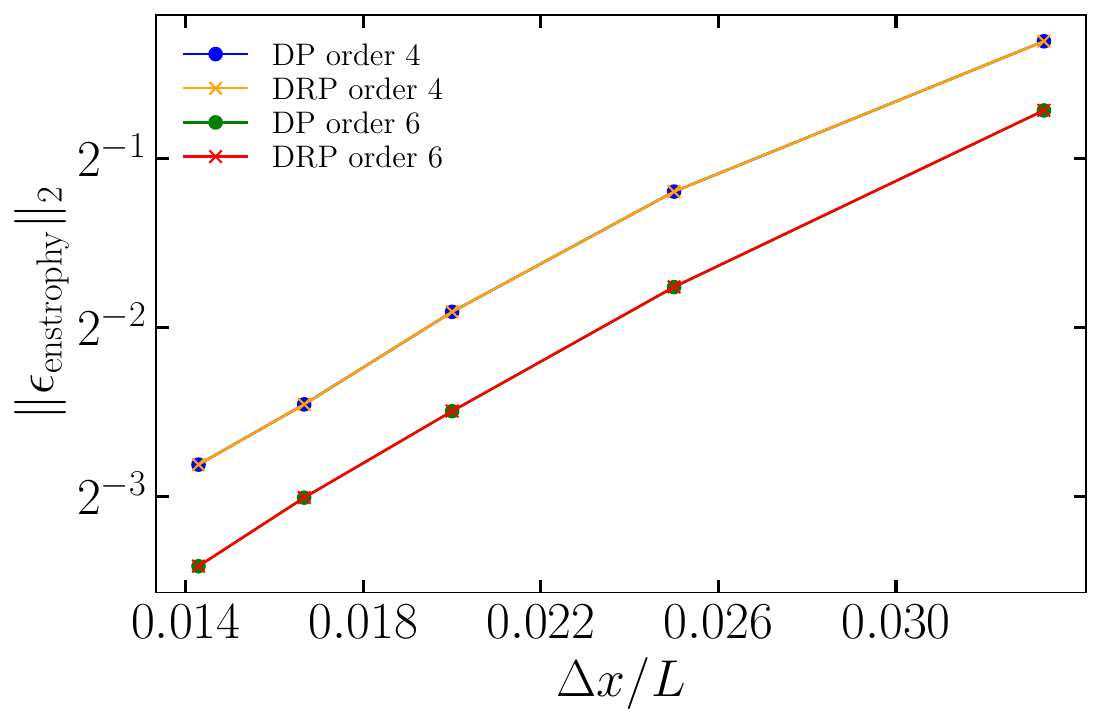}
        \label{fig:enstrophy_error_merging}
    \end{subfigure}

    }\caption{\label{fig:merging_vortex_conv_error}$L_2$ conservation errors of mass, energy \revi{and enstrophy} for the case with hyper-viscosity, $\delta = 0.5 $ as a function of uniform grid spacing $\Delta x/L$, for both DP and DRP operators of order 4 and 6 in the merging vortex problem. The conservation errors converge to zero as expected.}
\end{figure}
\begin{table}
\centering

\begin{tabular}{c|c|c|c|c|c|c|c|c}
\hline 
$m$  & \multicolumn{2}{l}{$q_{\text{DP4}}$} & \multicolumn{2}{l}{$q_{\text{DP6}}$} & \multicolumn{2}{l}{$q_{\text{DRP4}}$} & \multicolumn{2}{l}{$q_{\text{DRP6}}$}\\ \hline \hline
   & energy             & mass            & energy             & mass            & energy             & mass             & energy             & mass             \\
31 & -                  & -               & -                  & -               & -                  & -                & -                  & -                \\
41 &  2.9318            &  2.8919         & 4.0100             & 3.9404          & 2.9361             & 2.8979           & 4.0112             & 3.9412           \\
51 &  3.2222            &  3.1834         & 4.3111             & 4.2144          & 3.2241             & 3.1860           & 4.3133             & 4.2165                \\
61 &  3.3657            &  3.3269         & 4.1944             & 4.0909          & 3.3675             & 3.3291           & 4.1948             & 4.0911                \\
71 &  3.4232            &  3.3837         & 3.9652             & 3.8782          & 3.4252             &  3.3862          & 3.9654             &  3.8783                \\ \hline \hline
\end{tabular}
\caption{\label{tab:convdisp_mergingvortex}Convergence rate, \revi{$q_h$} of energy and mass as a function of spatial resolution for DP and DRP operators of order 4 and 6, computed from the merging vortex problem, with grid-dependent hyper-viscosity scaling as roughly $\sim (\Delta x)^3$ for order 4 and $\sim (\Delta x)^4$ for order 6, and $\delta = 0.5$ (the sub-convergence is due to fourth order time discretisation).}

\end{table}
\subsubsection{Barotropic shear instability}
Here we consider the classical Galewsky \cite{galewsky2004initial} test case, but instead on a 2D plane (see Peixoto et al.~\cite{peixoto2019semi} for a similar setup), which involves triggering a barotropic shear within zonal jets by initialising the flow with a thin fluid discontinuity, supplemented with small perturbations. Here we use the following initial conditions,
\begin{equation}
u(x, y)_{\pm} =\pm u_0\text{sech}(k_v(y - y_{0_{\pm}})), \quad v(x, y)=0
\end{equation}
where $u_0=50$ is the maximum speed and $y_{0_{+}} = 3\times 10^7$, $y_{0_{-}} = 1\times 10^7$. To ensure linearised geostrophic balance between the depth and velocity fields (i.e., that the initial conditions are analytically in a (pseudo) steady state), we define the depth perturbation $h$, as,
$$
h(x, y)=-\frac{f}{g} \int_0^y u(x, s) d s .
$$
whose solution can be found analytically. Small Gaussian perturbations, $\widetilde{h}$ are added to $h$ to trigger the barotropic instability,
$$
\left.\widetilde{h}(x, y)=0.02 H\left[\exp \left\{-k d_1(x, y)\right)\right\}+\exp \left\{-k d_2(x, y)\right\}\right]
$$
where $k=1$, and $d_i(x, y)=\frac{\left(x-x_i\right)^2}{L_x^2}+\frac{\left(y-y_i\right)^2}{L_y^2}, i=\{1,2\}$, are the square Euclidean distances of $(x, y)$ to the points $p_1=\left(x_1, y_1\right)=\left(0.75 L_x, 0.95 L_y\right), p_2=\left(x_2, y_2\right)=$ $\left(0.25 L_x, 0.2375 L_y\right)$, respectively, and we use $L_x = L_y = 4\times 10^7$, with $H = 1 \times 10^4$, $f = 2 \times 7.292 \times 10^{-5}$, $k_v = 10^{-7}$ and $g = 1$.

Figure~\ref{fig:KH} displays the time evoluton of the two-dimensional vorticity contour at the resolution of $401^2$ using DP operators of order 4 with hyper-viscosity ($\delta = 10$). The initial Gaussian perturbations trigger fast inertial gravity waves which subsequently interact with the zonal jet flow, resulting into large scale vortices from baroclinic generation of vorticity. These shearing motions result in the so-called Kelvin-Helmholtz instability typically observed in geostrophic flows. 

Our numerical scheme resolves the large and small-scale features of the vortices very well, which is not typically the case in certain variational numerical methods like CG or DG due to the development of spurious wave modes; or low-order dissipative schemes that do not resolve turbulent motions occurring on multiple scales.
\begin{figure}
        \begin{subfigure}{0.5
\textwidth}
    \centering
    \includegraphics[width = 1\linewidth]{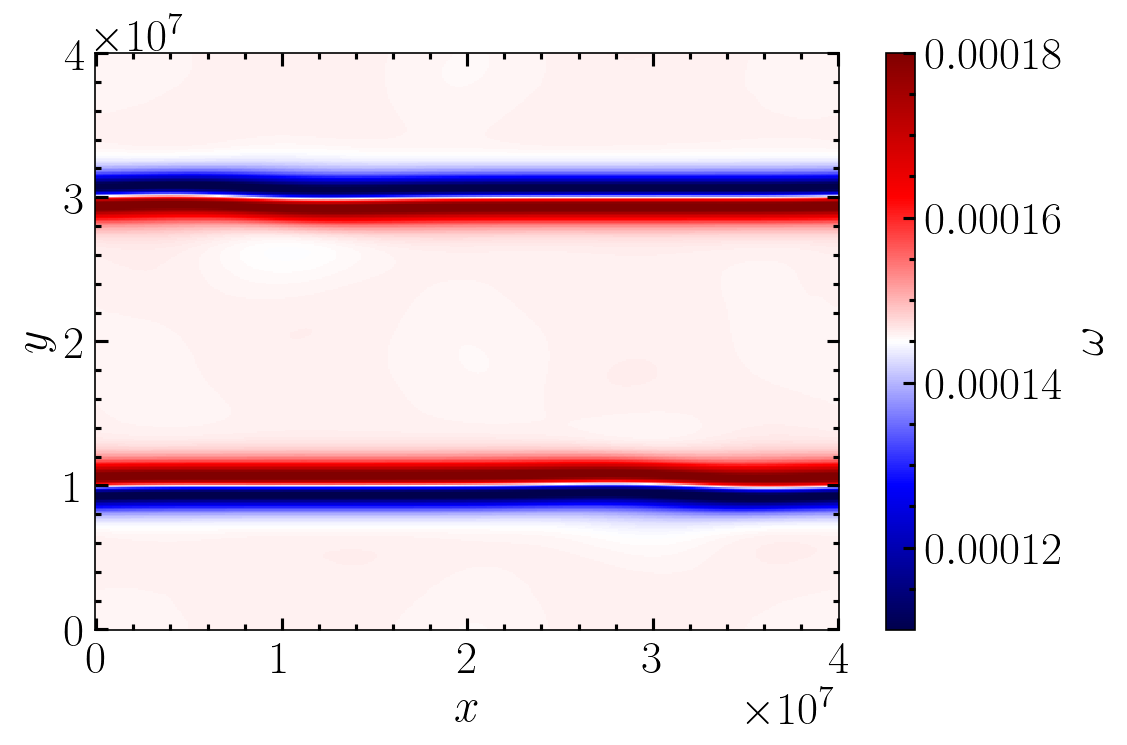}
    \caption{\label{fig:with_hyper_u}$t = 0$ days}
    \end{subfigure}
    \begin{subfigure}{0.5\textwidth}
    \centering
    \includegraphics[ width = 1\linewidth]{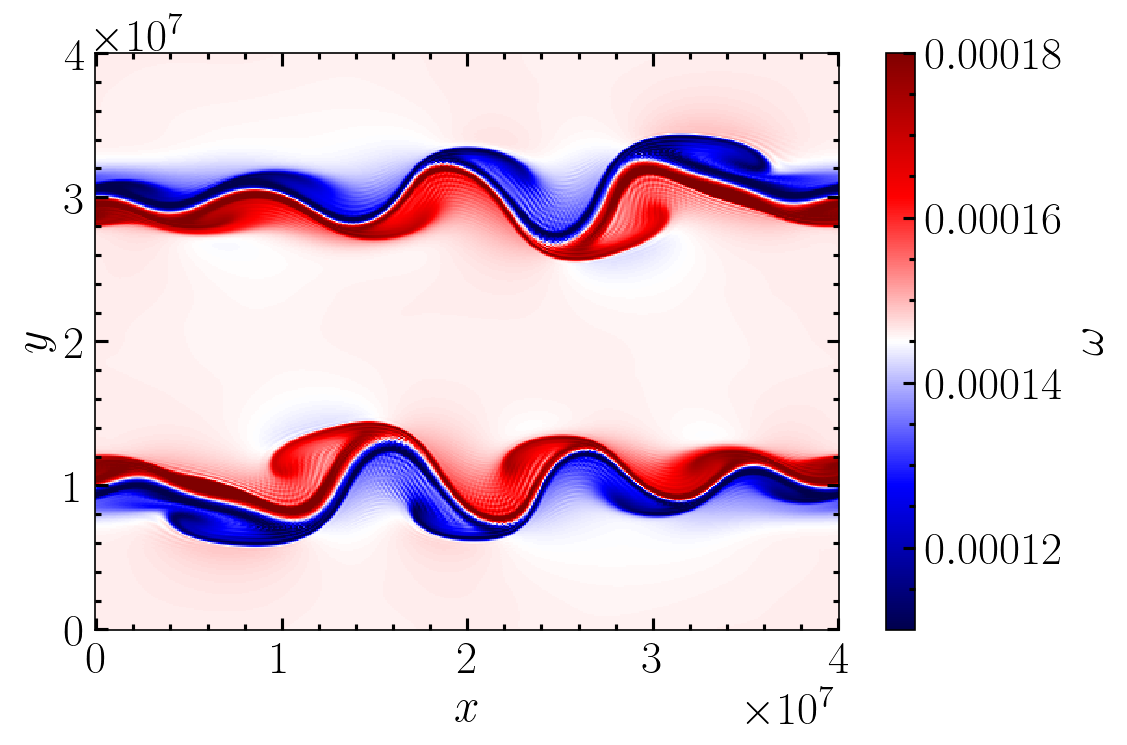}
    \caption{\label{fig:lake_small}$t = 14$ days}
    \end{subfigure}
    \begin{subfigure}{0.5
\textwidth}
    \centering
    \includegraphics[width = 1\linewidth]{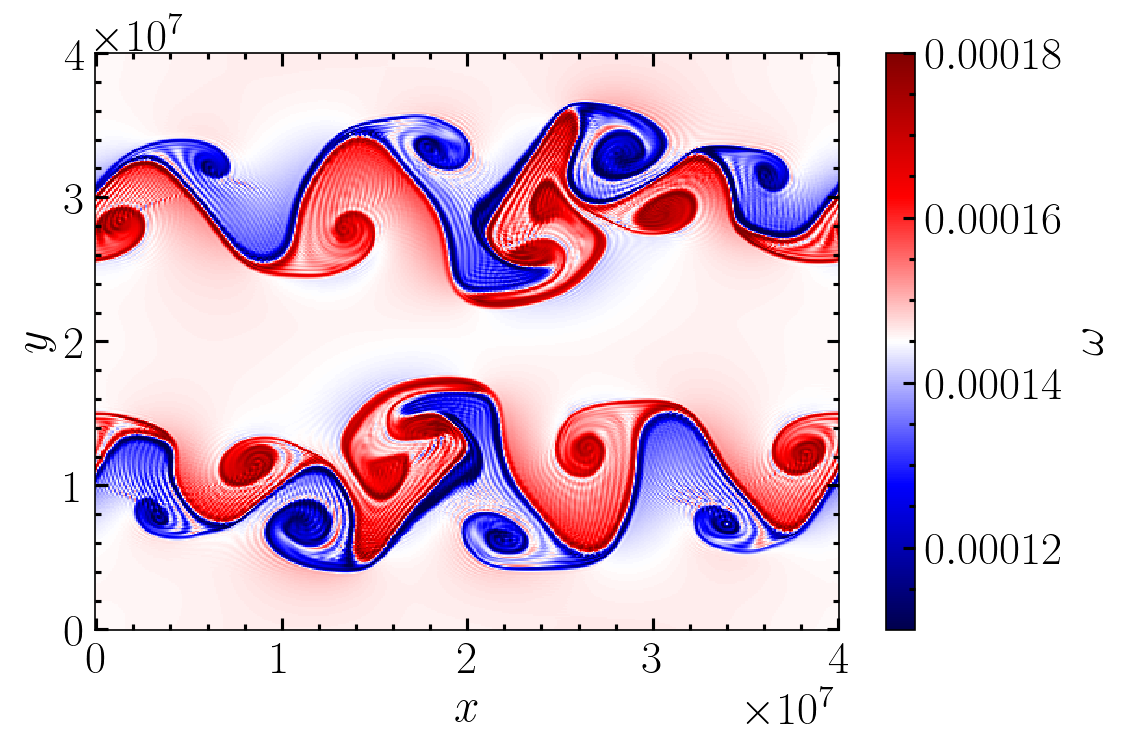}
    \caption{\label{fig:no_hyper_u}$t = 20$ days}
    \end{subfigure}
        \begin{subfigure}{0.5\textwidth}
    \centering
    \includegraphics[ width = 1\linewidth]{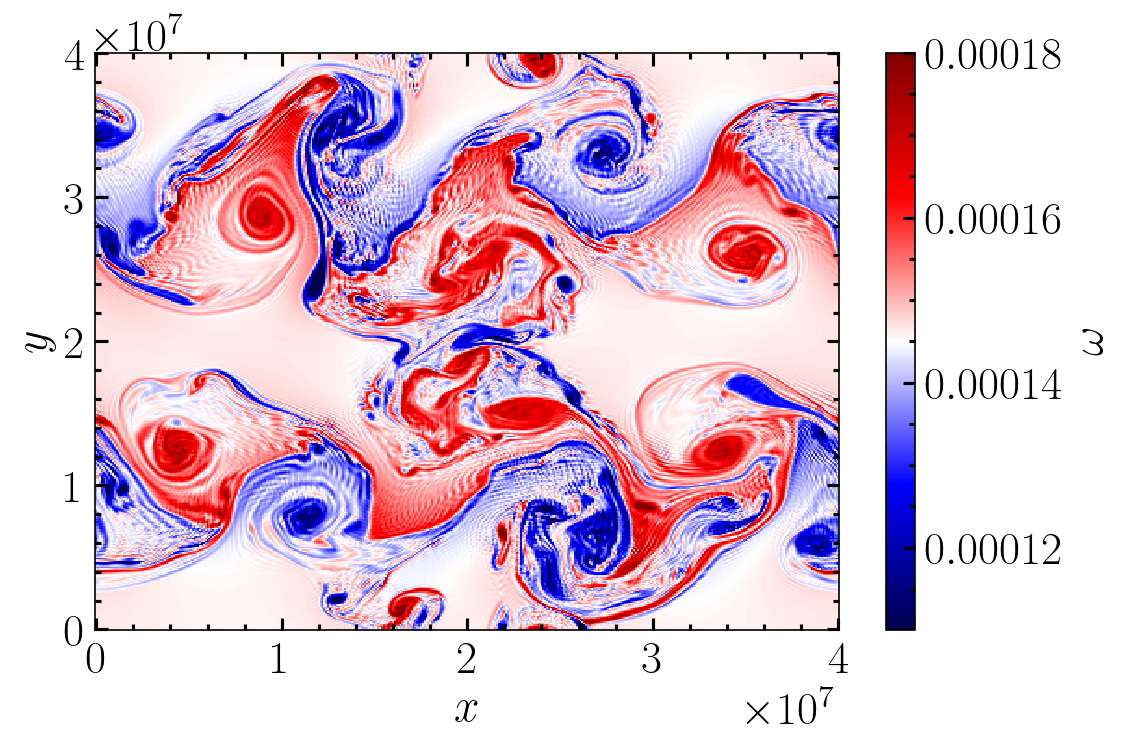}
    \caption{\label{fig:with_hyper_h}$t = 36$ days}
    \end{subfigure}
    \begin{subfigure}{0.5
\textwidth}
    \centering
    \includegraphics[width = 1\linewidth]{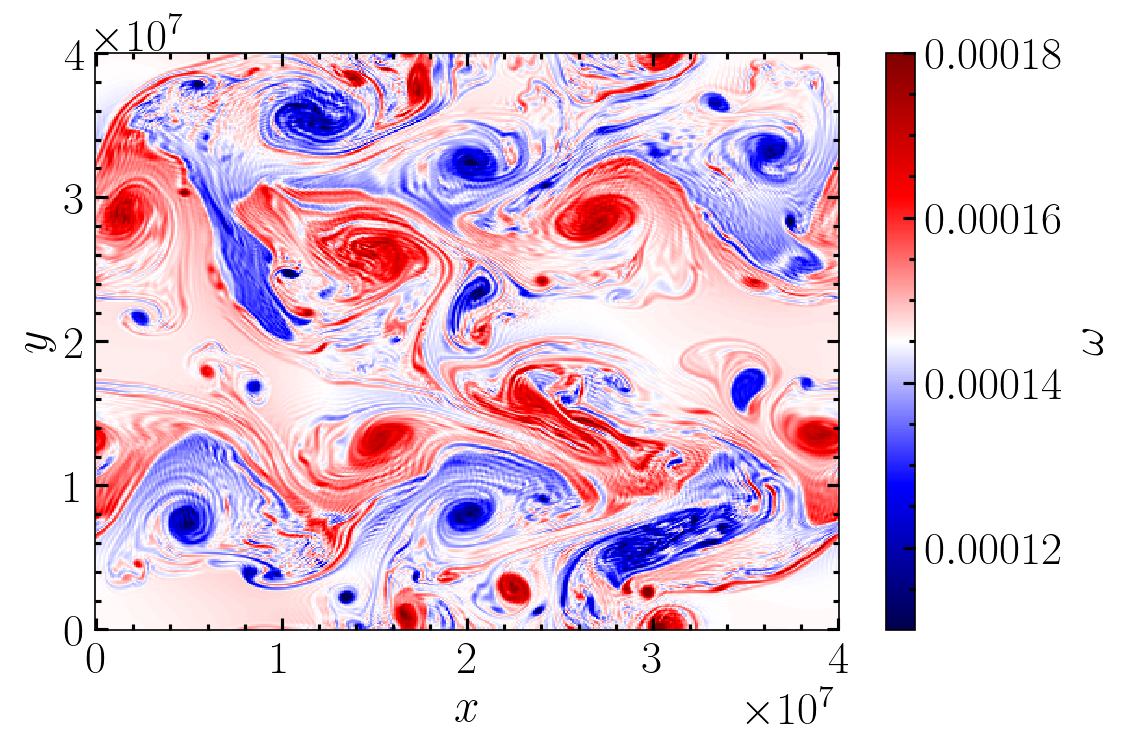}
    \caption{\label{fig:with_hyper_u}$t = 54$ days}
    \end{subfigure}
        \begin{subfigure}{0.5
\textwidth}
    \centering
    \includegraphics[width = 1\linewidth]{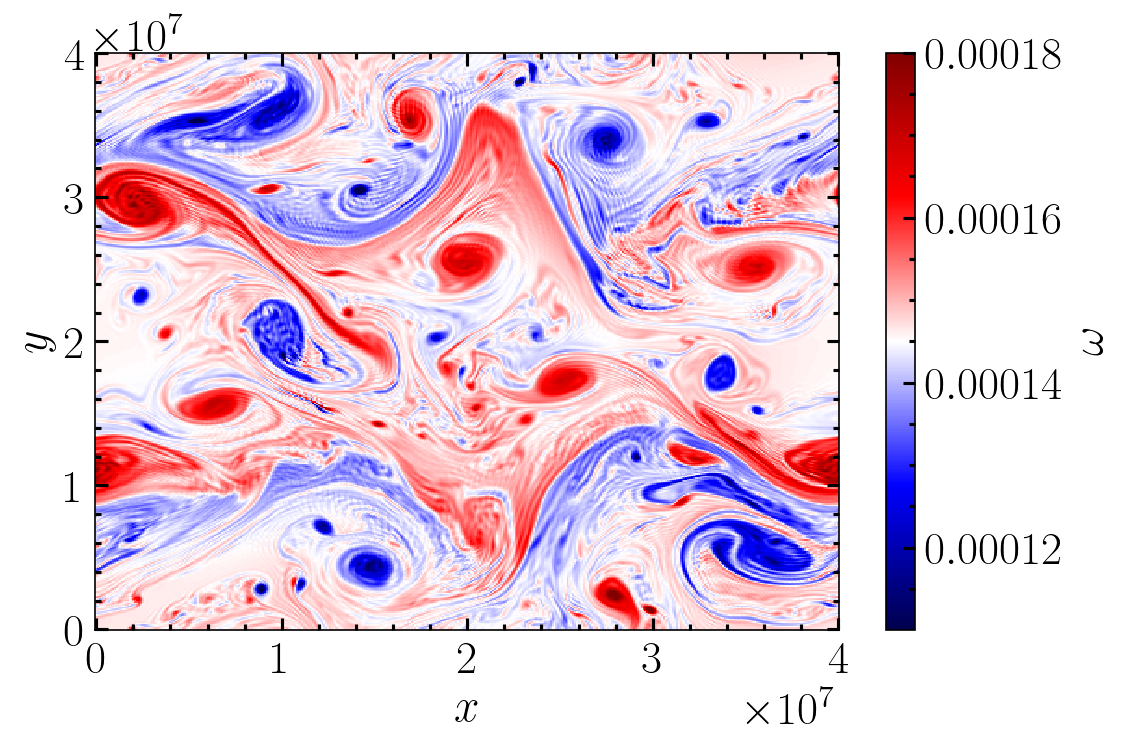}
    \caption{\label{fig:with_hyper_u}$t = 76$ days}
    \end{subfigure}

\caption{\label{fig:KH} Snapshots of 2D vorticity ($\omega$) contour at various times in the shear instability simulation. A movie of this simulation can be found at \url{https://youtu.be/jJayYT-eMTI?si=5dIkQ2S1_ILSY6vn}. }
\end{figure}
Moreover, in order to evaluate the fidelity of our numerical scheme in resolving turbulent scales, we show the shell-integrated kinetic energy and enstrophy spectra in our simulations (Fig.~\ref{fig:spectra}), where the kinetic energy spectra is defined using the Fourier transformed velocities $(\hat{u}(\vec{k}), \hat{v}(\vec{k}))$, as
$$
E_n=\sum_{n \leq\|\vec{k}\|<n+1} E_{\vec{k}}, \quad \text { where } \quad E_{\vec{k}}=\frac{1}{2}\left(\hat{u}(\vec{k}) \hat{u}^*(\vec{k})+\hat{v}(\vec{k}) \hat{v}^*(\vec{k})\right)
$$
where $\vec{k}=\left(k_1, k_2\right)$ represents a longitudinal mode, ${ }^*$,  the complex conjugate, $\|\vec{k}\|=\sqrt{k_1^2+k_2^2}$, and $E_n$ represents the spectral density with respect to the wavenumber $n$ and wavelength $L / n$ with $L$ the size of the domain. The enstrophy spectra, assuming isotropy is defined as $E_\omega = k^2 E_k $, along with the expected power law cascades from two-dimensional hydrodynamic turbulence \cite{boffetta2012two}. The kinetic energy spectra in Fig.~\ref{fig:ekin_spec}, shows a reasonable Kraichnan spectra ($k^{-3}$) within the inertial range, up to small-scale intermittency or finite resolution effects. Therefore, we resolve the direct enstrophy cascade reasonably well despite somewhat lower resolution than most canonical two-dimensional turbulence studies (e.g., \cite{boffetta2010evidence}). Simlarly, the enstrophy spectra (Fig.~\ref{fig:ens_spec}) also displays a reasonable $k^{-1}$ law within a narrow region of the inertial range, while we note that the enstrophy cascade, which is a second moment of velocity, typically takes significant spatial resolution to achieve good convergent behaviour in the inertial range.

\begin{figure}
 
    \centering

    \begin{subfigure}{0.48\textwidth}
    \centering
    \includegraphics[ width = 1\linewidth]{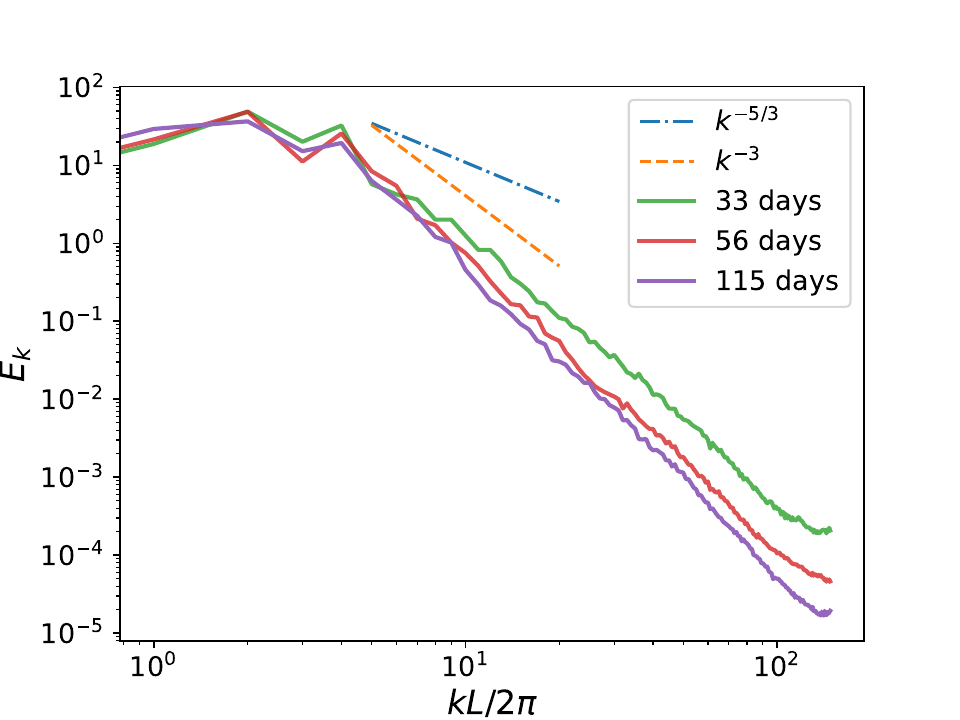}
    \caption{\label{fig:ekin_spec} }
    \end{subfigure}
    \begin{subfigure}{0.48
\textwidth}
    \centering
    \includegraphics[width = 1\linewidth]{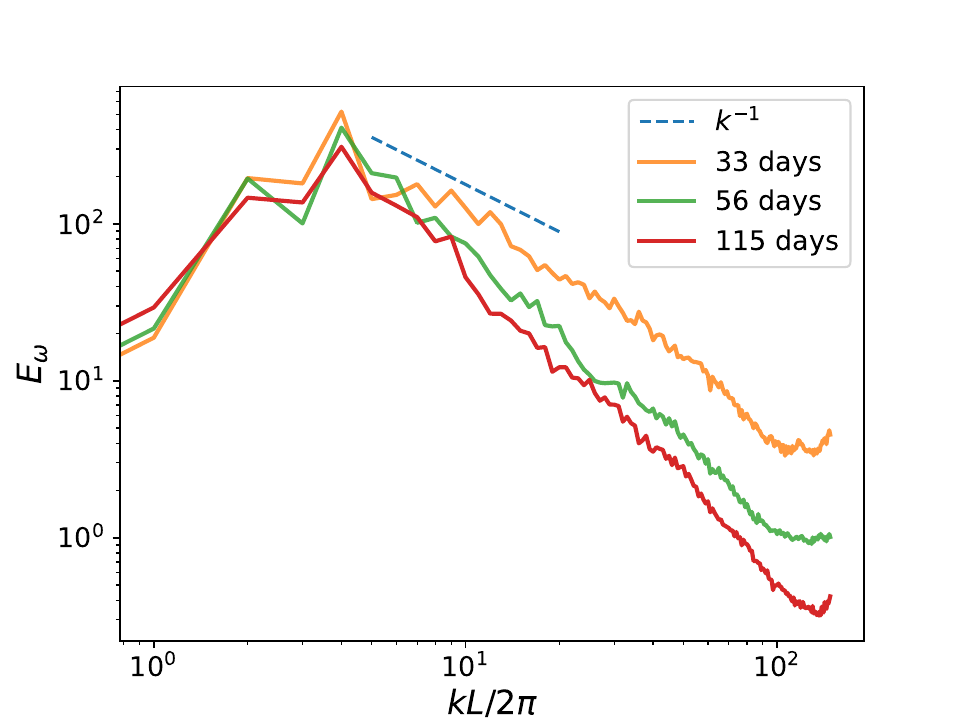}
    \caption{\label{fig:ens_spec} }
    \end{subfigure}
    \caption{\label{fig:spectra} Spectra from the baroclinic instability test case, (a) Kinetic energy spectra. (b) Enstrophy spectra. Power law lines expected from two-dimensional inverse energy cascade and direct enstrophy cascade are placed for the eye.}
\end{figure}
\section{Conclusion}\label{sec:conclusion}
A novel, energy/entropy stable, and high order accurate numerical scheme has been developed and analysed to solve the 1D and 2D nonlinear and linear SWEs in vector invariant form. We prove both nonlinear and local linear stability  results, as well as a priori error estimates for the semi-discrete approximations of the linear and nonlinear SWEs. The analysis is verified by computing the numerical eigenvalues of the evolution operator, for both linear and nonlinear problems. 
The scheme is suitable for the simulations of  subcritical flow regimes typically prevalent in  atmospheric and geostrophic flow problems. The numerical scheme  is based on the recently derived high order accurate dual-pairing SBP operators \cite{mattsson2017diagonal,williams2021provably,CWilliams2021}. Convergence, accuracy, and well-balanced properties are validate via the method of manufactured solutions and canonical test problems such as the dam break and lake at rest. The method has been validated with a number of canonical test cases as well as with the method of manufactured smooth solutions, including 2D merging vortex problems  and barotropic shear instability, with fully developed turbulence.

For nonlinear problems, energy/entropy stability although ensures the boundedness of numerical solutions but it does not guarantee convergence of numerical errors with mesh refinement. Suitable amount of numerical dissipation is necessary to control high frequency errors which could dominate the numerical simulations. Using the dual-pairing SBP framework, we designed high order accurate and nonlinear hyper-viscosity operator which dissipates entropy. The hyper-viscosity operator is effective in controlling the oscillations associated with the  shocks and discontinuities and high frequency grid-scale errors. 

  Another important results here are new well-posed BCs for the 1D vector invariant form of the SWE formulated in terms of fluxes of the primitive variables, and applicable to linear and nonlinear problems. We prove that the BCs are well-posed and energy/entropy stable. Numerical experiments are performed to verify the theoretical results, and in particular the nonlinear transmissive BCs. 

Future work will extend the method to spherical and complex geometries using  curvilinear grids. We will also extend the 1D nonlinear BCs to 2D which will enable efficient numerical simulations in non-periodic bounded domains.


\section{Acknowledgements}
We thank Kenny Wiratama and David Lee for insightful discussions about higher order discretisation methods for hyperbolic systems. We are also grateful to Christopher Williams, Alberto Mart{\'i}n, Rudi Prihandoko, James R. Beattie, Neco Kriel, Markus Hegland and Eric W. Hester for general discussions regarding this work. J.K.J.H. acknowledges funding via the ANU Chancellor's International Scholarship, the Space Plasma, Astronomy and Astrophysics Research Award, the Boswell Technologies Endowment Fund and Higher Degree Research (HDR) Award for Space Exploration. We also acknowledge computational resources provided by the National Computational Infrastructure (NCI) under the grant \texttt{xx52}, which is supported by the Australian Government.
\section{Data Availability Statement}
The datasets produced/analysed in the present study are available upon reasonable request to the corresponding author (J.K.J.H.).

\appendix

\section{Additional convergence studies for traditional SBP and odd-order DP and DRP operators}\label{appen:addconv}
Here in Tables~\ref{table:sine} and  \ref{table:sine2}, we show additional convergence tables with the Gaussian MMS conducted using the traditional SBP operators (Table~\ref{table:sine}), as well as odd-order DP and DRP operators (Table~\ref{table:sine2}). As confirmed earlier, we observe superconvergence for 4th and 6th order SBP operators. However, traditional SBP operators achieve globally $3$rd order convergence rate for order 4 interior operator, consistent with the $p+1$ global accuracy expected from Laplace transform analysis, and $5$th order global convergence rate for order 6 interior operator. 

However, as shown in Table~\ref{table:sine2}, the DP and DRP odd order operators exhibit super-convergent behaviour, which has also been shown in a number of recent works \cite{mattsson2017diagonal,mattsson2018compatible,lundgren2020efficient}, due to the nature of these types of upwind SBP operators, while we note that the precise reason for this is still unknown.
\begin{table}[H]
 \resizebox{\textwidth}{!}	{\begin{tabular}{c c c c c c c c c}
		\hline  \hline
   \multicolumn{1}{c}{\texttt{SBP, order 4, MMS} } \\
& \multicolumn{4}{c}{Linear} 
 & 
 \multicolumn{4}{c}{Nonlinear}  \\ 
 		\multicolumn{1}{c}{$N$} & $\log_2 \| \textrm{error}\|_{2,u}$ & $\log_2\|\textrm{error}\|_{2,h}$ & $q_u$    & $q_h$     & $\log_2 \| \textrm{error}\|_{2,u}$ & $\log_2\|\textrm{error}\|_{2,h}$ & $q_u$    & $q_h$    \\
    \multicolumn{1}{c}{(1)} & (2) & (3) & (4) & (5) & (6) & (7) & (8) & (9)\\ \hline
		41  & -5.8130      &  -5.7980                          & -        & -       &  -8.7047 &  -8.6876 & -  & -   \\ 
		81  & -8.9000                            & -8.8860                         & 3.0870 & 3.0880 & -5.0059 & -5.1572 & 3.6988 & 3.5304 \\ 
		161  & -11.9541                            & -11.9402                        & 3.0541 & 3.0542  & -11.9884 & -11.9613 & 3.2837 & 3.2737\\ 
 		321 & -14.9679                            & -14.9541                          & 3.0139 & 3.0139 & -15.0822 & -15.0326 & 3.0937 & 3.0714  \\ 
 		641 & -17.9714                           & -17.9576                        & 3.0035 & 3.0035 & -18.1097 & -18.0501 & 3.0276 & 3.0175\\ \hline \hline
 	\end{tabular}}
  	\resizebox{\textwidth}{!}{\begin{tabular}{c c c c c c c c c}

  \multicolumn{1}{c}{\texttt{SBP, order 6, MMS} }
  \\
	&	\multicolumn{4}{c}{Linear} 
&
\multicolumn{4}{c}{Nonlinear}  \\ 
		\multicolumn{1}{c}{$N$} & $\log_2 \| \textrm{error}\|_{2,u}$ & $\log_2\|\textrm{error}\|_{2,h}$ & $q_u$    & $q_h$     & $\log_2 \| \textrm{error}\|_{2,u}$ & $\log_2\|\textrm{error}\|_{2,h}$ & $q_u$    & $q_h$     \\
    \multicolumn{1}{c}{(1)} & (2) & (3) & (4) & (5) & (6) & (7) & (8) & (9)\\  \hline
		41  & -5.5447      &  -5.4670                         & -        & - & -4.1892 & -4.1338 & - & -          \\ 
		81  & -9.5738                            & -9.5487                         & 4.0291 & 4.0817  & -9.1780 & -9.1673  & 4.9888 & 5.0335\\ 
		161  & -15.1548                            & -15.1388                       & 5.5810 & 5.5901 & -14.7453 & -14.7211 &  5.5673 & 5.5538 \\ 
		321 & -20.3594                           & -20.3471                         & 5.2047  & 5.2083 & -19.8871 & -19.8401 & 5.1417 & 5.1190\\
  		641 & -25.3107                           & -25.3672                         & 4.9513  & 5.0201 & -24.8160 & -24.8400 & 4.9289 & 4.9999 \\ 
	 \hline \hline
	\end{tabular}}
 \caption{\label{table:sine}Log$_2 L^2$ norm errors and convergence rate of SBP operators of even order (order 4 and 6) for primitive variables $u,h$.}
 \end{table}
 \begin{table*}
\centering
	  \resizebox{\textwidth}{!}{\begin{tabular}{c c c c c c c c c}
		\hline  \hline
  \multicolumn{1}{c}{\texttt{DP, order 5, MMS} } \\
	&	\multicolumn{4}{c}{Linear} 
&
\multicolumn{4}{c}{Nonlinear}  \\ 
		\multicolumn{1}{c}{$N$} & $\log_2 \| \textrm{error}\|_{2,u}$ & $\log_2\|\textrm{error}\|_{2,h}$ & $q_u$    & $q_h$     & $\log_2 \| \textrm{error}\|_{2,u}$ & $\log_2\|\textrm{error}\|_{2,h}$ & $q_u$    & $q_h$     \\
    \multicolumn{1}{c}{(1)} & (2) & (3) & (4) & (5) & (6) & (7) & (8) & (9)\\ \hline 
		41  & -7.2370                              &  -7.1095                           & -        & - & -7.0995 & -6.7450 & - & -         \\ 
		81  & -11.3378                            & -10.9423                         & 4.1008 & 3.8328 & -11.6030 & -10.6515 & 4.5035 & 3.9064\\ 
		161  & -15.2832                            & -15.0097                         & 3.9454 & 4.0674 & -15.4952 & -14.7624 & 3.8921 & 4.1109  \\ 
		321 & -18.9627                            & -18.8146                          & 3.6795 & 3.8049 & -19.1202 & -18.6201 & 3.6250 & 3.8577 \\ 
		641 & -22.5701                           & -22.5139                        & 3.6074 & 3.6993 & -22.7087 & -22.3616 & 3.5885 & 3.7415 \\ 
	\end{tabular}}
  \resizebox{\textwidth}{!}{\begin{tabular}{c c c c c c c c c}
		\hline  \hline
  \multicolumn{1}{c}{\texttt{DRP, order 5, MMS} } \\
	&	\multicolumn{4}{c}{Linear} 
&
\multicolumn{4}{c}{Nonlinear}  \\ 
		\multicolumn{1}{c}{$N$} & $\log_2 \| \textrm{error}\|_{2,u}$ & $\log_2\|\textrm{error}\|_{2,h}$ & $q_u$    & $q_h$     & $\log_2 \| \textrm{error}\|_{2,u}$ & $\log_2\|\textrm{error}\|_{2,h}$ & $q_u$    & $q_h$     \\
    \multicolumn{1}{c}{(1)} & (2) & (3) & (4) & (5) & (6) & (7) & (8) & (9)\\ \hline 
		41  & -6.7009        &  -5.7814                           & -        & - & -6.5358 & -5.1352 & - & -         \\ 
		81  & -10.2166      & -9.7475                         & 3.5158 & 3.9661 & -10.1635 & -9.2776 & 3.6276 & 4.1424\\ 
		161  & -14.2013                            & -13.9006                      & 3.9846 & 4.1531 & -14.1764 & -13.5178 & 4.0129 & 4.2402  \\ 
		321 & -17.9694                            & -17.8282                          & 3.7681 & 3.9276 & -17.9550 & -17.5155 & 3.7786 & 3.9977 \\ 
		641 & -21.6402                           & -21.6193                        & 3.6708 & 3.7911 & -21.6441 & -21.3592 & 3.6891 & 3.8437 \\ \hline \hline
	\end{tabular}}
\caption{\label{table:sine2}Log$_2 L^2$ norm errors and convergence rate of DP and DRP operators of odd order (order 5) for primitive variables $u,h$.}
\end{table*}
\section{2D energy/entropy stable semi-discrete approximation}\label{sec:hyper2d}
To extend the 1D semi-discrete approximation \eqref{eq:semi-discrete-SWE}  to the 2D rotating SWE \eqref{eqn:nlSWE2D} we discretise the 2D doubly-periodic domain $(x, y)\in \Omega = [0, 2\pi] \times [0, 2\pi]$ into $(N+1)\times (N+1)$ grid points  $(x_i, y_j)= (i\Delta{x}, j\Delta{y})$ with the uniform spatial step $\Delta{x}=\Delta{y} = 2\pi/N$ for $N\in \mathbb{N}$. A 2D scalar field $v(x, y)$   on the grid $(x_i, y_j)$ is re-arranged row-wise into a vector  $\mathbf{v}\in \mathbb{R}^{(N+1)\times (N+1)}$.
The 2D rotating SWE in semi-discrete form can be written as:
    \begin{equation}\label{eqn:nl2dSWE2D}
    \begin{aligned}
   \frac{d \mathbf{h}}{d t } + \nabla_{\mathbf{D}_+} \cdot (
    \mathbf{U}\mathbf{h}) = G_h, \quad \frac{d \mathbf{U}}{d t} + \boldsymbol{\omega} \mathbf{U}^\perp + \nabla_{\mathbf{D}_-} \biggl(\frac{|\mathbf{U}|^2}{2} + g\mathbf{h}\biggl) = \mathbf{G}_u\textcolor{blue}{,}  \\ \quad \boldsymbol{\omega} = \mathbf{D}_{-x}\mathbf{v} - \mathbf{D}_{-y} \mathbf{u} + f_c,   \quad \mathbf{q}(0) =  \mathbf{f},
    \end{aligned}
\end{equation}
where $\mathbf{U} = (\mathbf{u}, \mathbf{v})^T$, $\mathbf{U}^\perp = (-\mathbf{v}, \mathbf{u})^T$, $|\mathbf{U}|^2 = \mathbf{u}^2+\mathbf{v}^2$, $\nabla_{\mathbf{D}_{\eta}}  = (\mathbf{D}_{\eta x}, \mathbf{D}_{\eta y})^T$, with $\eta \in \{+,-\}$. is the discrete approximation of  gradient operator $\nabla = (\partial/\partial x, \partial/\partial y)^T$. As before, note that the derivatives for continuity equation are approximated with $D_+$ and the derivatives for the momentum equation are approximated with the dual operator $D_-$.  Using the Kronecker product $\otimes$, the 2D DP derivative operators are defined as
$\mathbf{D}_{\pm x} = D_\pm \otimes \mathbf{I}_y$, $\mathbf{D}_{\pm y} = \mathbf{I}_x\otimes D_\pm $ where $D_\pm $ are 1D DP operators closed with periodic BCs, and satisfy $D_{+} + D_{-}^T =0$.  The 2D SBP norm is defined as $\mathbf{P} = \Delta{x}\Delta{y} \mathbf{I}$ with $\mathbf{I}=\left(\mathbf{I}_x\otimes \mathbf{I}_y\right)$ and $\mathbf{I}_x, \mathbf{I}_y \in \mathbb{R}^{(N+1)\times (N+1)}$ being identity matrices. 

We define the 2D energy functional  
\begin{equation}\label{eqn:energy_2d}
\begin{aligned}
\|\mathbf{q}\|_{\mathbf{P}\mathbf{W}^{\prime}}^2:=\mathbf{q}^T \left(\mathbf{P}\mathbf{W}^{\prime}\right)\mathbf{q}
=
\Delta{x}\Delta{y}\sum_{i,j=1}^{N}\left(g h_{ij}^2 + h_{ij} u_{ij}^2 + h_{ij} v_{ij}^2  \right),  
\quad \mathbf{W}^{\prime} = \begin{bmatrix}
g\mathbf{I} & \mathbf{0}  & \mathbf{0} \\
\mathbf{0} & \text{diag}(\mathbf{h})  & \mathbf{0} \\
\mathbf{0} & \mathbf{0} & \text{diag}(\mathbf{h}) 
\end{bmatrix},
\end{aligned}
\end{equation} 
and the elemental energy
\begin{equation}\label{eqn:elemental_2d}
e_{ij}: = \frac{1}{2}\left(g h_{ij}^2 + h_{ij} u_{ij}^2 + h_{ij} v_{ij}^2  \right),
\end{equation}
defines an entropy functional for subcritical flows $u_{ij}^2 + v_{ij}^2 < gh_{ij}$. 
That is the elemental energy $e_{ij}$ is a convex function in terms of the prognostic variables $h_{ij}, u_{ij},  v_{ij}$. It is also noteworthy that,
\begin{equation}\label{eqn:W2d}
\begin{aligned}
\mathbf{W}\mathbf{q} 
= 
\begin{bmatrix}
\frac{|\mathbf{U}|^2}{2} + g\mathbf{h}\\
    \mathbf{u}\mathbf{h}\\
    \mathbf{v}\mathbf{h}
\end{bmatrix},
\quad 
  \mathbf{q}^T\left(\mathbf{P}\mathbf{W}\right)\frac{d\mathbf{q} }{dt} = \frac{1}{2}\frac{d}{dt}\biggl(\mathbf{q}^T\left(\mathbf{P}\mathbf{W}^{\prime}\right)\mathbf{q}\biggl), 
  \quad 
 \mathbf{W} = \begin{bmatrix}
g\mathbf{I} & \frac{1}{2}\text{diag}(\mathbf{u}) & \frac{1}{2} \text{diag}(\mathbf{v})\\
\frac{1}{2}\text{diag}(\mathbf{u}) & \frac{1}{2}\text{diag}(\mathbf{h})
& 0 \\
\frac{1}{2}\text{diag}(\mathbf{v}) & 0 & \frac{1}{2}\text{diag}(\mathbf{h})\end{bmatrix}, \quad \mathbf{W} = \mathbf{W}^T.
\end{aligned}
\end{equation}

To make the 2D analysis amenable to the 1D theory we will reformulate the semi-discrete approximation \eqref{eqn:nl2dSWE2D}.
As before  the 2D hyper-viscosity operator is constructed such that when appended to the semi-discrete formulation \eqref{eqn:nl2dSWE2D}, it becomes,
\begin{equation}\label{eq:semi-discrete-SWE_Diss2D}
   \frac{d\mathbf{q}}{dt} + \nabla_{\mathbf{D}} \cdot \mathbf{F} + \mathbf{\Omega}_f = \mathbf{\mathcal{K}} \mathbf{q} + \mathbf{G}, 
   \quad \mathbf{\Omega}_f= 
   \begin{bmatrix}
   \mathbf{0}\\
    -\boldsymbol{\omega}\mathbf{v}\\
    \boldsymbol{\omega}\mathbf{u}
\end{bmatrix},
\quad
\boldsymbol{\omega} = \mathbf{D}_{-x}\mathbf{v} - \mathbf{D}_{-y} \mathbf{u} + f_c,
   \quad \mathbf{q}(0) =  \mathbf{f},
\end{equation}
where
$\nabla_{\mathbf{D}} \cdot \mathbf{F} = \mathbf{D}_x\mathbf{F}_x + \mathbf{D}_y\mathbf{F}_y $
and the flux and its gradients are given by
$$
\begin{aligned}
\mathbf{F}_x := 
\begin{bmatrix}
    \mathbf{u}\mathbf{h}\\
    \frac{|\mathbf{U}|^2}{2} + g\mathbf{h}\\
     \mathbf{0}
\end{bmatrix},
\quad
\mathbf{F}_y := 
\begin{bmatrix}
    \mathbf{v}\mathbf{h}\\
    \mathbf{0}\\
    \frac{|\mathbf{U}|^2}{2} + g\mathbf{h}
\end{bmatrix},
\quad
\mathbf{D}_x\mathbf{F}_x := 
\begin{bmatrix}
    \mathbf{D}_{+x}\left(\mathbf{u}\mathbf{h}\right)\\
    \mathbf{D}_{-x}\left(\frac{|\mathbf{U}|^2}{2} + g\mathbf{h}\right)\\
    \mathbf{0}
\end{bmatrix},
\quad
\mathbf{D}_y\mathbf{F}_y := 
\begin{bmatrix}
    \mathbf{D}_{+y}\left(\mathbf{v}\mathbf{h}\right)\\
    \mathbf{0}\\
    \mathbf{D}_{-y}\left(\frac{|\mathbf{U}|^2}{2} + g\mathbf{h}\right)
\end{bmatrix}.
\end{aligned}
$$
Here the 2D nonlinear hyper-viscosity operator $\mathbf{\mathcal{K}}$ is given by 
\begin{equation}\label{eq:SWE_Diss_2D}
\begin{aligned}
\mathbf{\mathcal{K}}=\mathbf{W}^{-1}\left(I_3\otimes \mathbf{P}^{-1}\right)\left(\left(I_3\otimes \mathcal{A}\otimes \mathbf{I}_y\right) + \left(I_3\otimes \mathbf{I}_x\otimes \mathcal{A}\right)\right), \quad \mathcal{A} = \mathcal{A}^T, \quad \quad \mathbf{u}^T\mathcal{A}\mathbf{u} \le 0, \quad \forall \mathbf{u}\in \mathbb{R}^{N+1},
\end{aligned}
\end{equation}
where $\mathcal{A}$ is the 1D hyper-viscosity operator derived in Section \ref{sec:hyperviscosity} and $I_3 \in \mathbb{R}^{3\times 3}$ is the  identity matrix.

 We are now ready to prove that the numerical  method \eqref{eqn:nl2dSWE2D} or \eqref{eq:semi-discrete-SWE_Diss2D} is strongly energy stable.  Analogous to the 1D result Theorem~\ref{thm:swp_Diss_nonlinear} we also have
\begin{theorem}\label{thm:swp2d_Diss_nonlinear}
    Consider the semi-discrete approximation \eqref{eqn:nl2dSWE2D} with periodic BCs and the nonlinear hyper-viscosity  operator $\mathcal{K}$ defined by \eqref{eq:SWE_Diss_2D}. Let the symmetric matrix $\mathbf{W}$ and the diagonal matrix $\mathbf{W}^{\prime}$ be defined as in \eqref{eqn:energy_2d} and \eqref{eqn:W2d}.
    For subcritical flows with $\textrm{Fr}^2 = (u_{ij}^2 + v_{ij}^2)/ (gh_{ij}) < 1$, the semi-discrete approximation \eqref{eqn:nl2dSWE2D} or \eqref{eq:semi-discrete-SWE_Diss2D}  is strongly stable. That is, the numerical solution $\mathbf{q}$ satisfies the estimate,
    $$
    \| \mathbf{q}\|_{W^{\prime}P} \leq K e^{\mu t}\biggl(\| \mathbf{f}\|_{W^{\prime}P} +  \max\limits_{\tau \in [0, t]}\|\mathbf{G(\tau})\|_{W^{\prime}P}\biggl),
    $$
 where
    $$
    \mu = \max_{\tau \in[0,t]}{\frac{\mathbf{q}^T\left(\left(I_3\otimes \mathcal{A}\otimes \mathbf{I}_y\right) + \left(I_3\otimes \mathbf{I}_x\otimes \mathcal{A}\right)\right)\mathbf{q}}{\| \mathbf{q}\|_{W^{\prime}P}^2}} \leq 0, \quad {K = \max\left(1, {(1-e^{-\mu t})}/{\mu}\right)}.$$
\end{theorem}
\begin{proof}
  Multiplying \eqref{eqn:nl2dSWE2D} with $\mathbf{q}^T\mathbf{P}\mathbf{W}$
from the left yields, 
\begin{align}
\begin{split}  
    \left(\mathbf{q},\frac{d \mathbf{q}}{d t}\right)_{\mathbf{W}\mathbf{P}} &=-\left(\mathbf{q},\nabla _{\mathbf{D}} \cdot \mathbf{F}\right)_{\mathbf{W}\mathbf{P}} +(\mathbf{q}, \mathcal{K}\mathbf{q})_{\mathbf{W}\mathbf{P}} +   {\left(\mathbf{q},\mathbf{G} \right)_{\mathbf{W}\mathbf{P}}}\\
    &= -\Delta x \Delta y(\mathbf{W}\mathbf{q})^T(\mathbf{D}_x \mathbf{F}_x + \mathbf{D}_y \mathbf{F}_y )+ {\left(\mathbf{q},\mathbf{G} \right)_{\mathbf{W}\mathbf{P}}} +(\mathbf{q}, \mathcal{K}\mathbf{q})_{\mathbf{W}\mathbf{P}} + (\mathbf{q}, \omega k \times \mathbf{u})_{\mathbf{W} \mathbf{P}} \\
    &= - \Delta x \Delta y\biggl[\biggl(\frac{|\mathbf{U}|^2}{2} + g\mathbf{h}\biggl)^T \mathbf{D}_{+x}\biggl(\mathbf{u}\mathbf{h}\biggl) + \biggl(\mathbf{u}\mathbf{h}\biggl)^T\mathbf{D}_{-x}\biggl(\frac{|\mathbf{U}|
    ^2}{2} + g\mathbf{h}\biggl) + \\& \biggl(\frac{|\mathbf{U}|^2}{2}   + g\mathbf{h}\biggl)^T \mathbf{D}_{+y}\biggl(\mathbf{v}\mathbf{h}\biggl)  + \biggl(\mathbf{v}\mathbf{h}\biggl)^T\mathbf{D}_{-y}\biggl(\frac{|\mathbf{U}|^2}{2} + g\mathbf{h}\biggl)\biggl] + {\left(\mathbf{q},\mathbf{G} \right)_{\mathbf{W}\mathbf{P}}} \\& +(\mathbf{q}, \mathcal{K}\mathbf{q})_{\mathbf{W}\mathbf{P}} + (\mathbf{q}, \omega k \times \mathbf{u})_{\mathbf{W} \mathbf{P}}.
\end{split}
    \end{align}
Where we note that vorticity does not contribute to the energy conservation equation since, 
\begin{align}
\begin{split}
(\mathbf{q}, \omega k \times \mathbf{u})_{\mathbf{W} \mathbf{P}} 
= \Delta x \Delta y \sum_{i=1}^N \sum_{j=1}^N  \omega_{ij}h_{ij}   \mathbf{u}_{ij} \cdot (\mathbf{k} \times \mathbf{u}_{ij}) 
= 0, 
\end{split}
\end{align}
which is analogous to the continuous setting. Noting the negative duality of the periodic DP operators, $\mathbf{D}_+ = -\mathbf{D}_-^T$, it is easy to see that the first bracketed term all cancel each other out. Finally, we have,
    \begin{align}
\begin{split}  
      \frac{1}{2}\frac{d }{dt} \|\mathbf{q}\|_{\mathbf{W}^{\prime}\mathbf{P}}^2  
    =  (\mathbf{q}, \mathcal{K}\mathbf{q})_{\mathbf{W}\mathbf{P}}+ {\left(\mathbf{q},\mathbf{G} \right)_{\mathbf{W}\mathbf{P}}} \le \mu \|\mathbf{q}\|_{\mathbf{W}^{\prime}\mathbf{P}}^2 + {\left(\mathbf{q},\mathbf{G} \right)_{\mathbf{W}\mathbf{P}}} = \mu \|\mathbf{q}\|_{\mathbf{W}^{\prime}\mathbf{P}}^2 + (\mathbf{q},\widetilde{\mathbf{G}})_{\mathbf{W}^{\prime}\mathbf{P}}.
\end{split}
    \end{align}
    where 
    $
   \mu = \max_{\tau \in[0,t]}{\frac{\mathbf{q}^T\left(\left(I_3\otimes \mathcal{A}\otimes \mathbf{I}_y\right) + \left(I_3\otimes \mathbf{I}_x\otimes \mathcal{A}\right)\right)\mathbf{q}}{\| \mathbf{q}\|_{W^{\prime}P}^2}} \leq 0, \quad {K = \max\left(1, {(1-e^{-\mu t})}/{\mu}\right)}. \leq 0,
    $ $\widetilde{\mathbf{G}} = {(\mathbf{W}^{\prime})}^{-1}\mathbf{W}\mathbf{G}$ and at subcritical flow regime we have $|\widetilde{\mathbf{G}}| \le  |\mathbf{G}|$.
    Cauchy-Schwartz inequality yields
\begin{align}
\begin{split}  
    \frac{1}{2} \frac{d}{dt} \| \mathbf{q}\|_{\mathbf{W}^{\prime}\mathbf{P}}^2   \leq \mu \|\mathbf{q}\|_{\mathbf{W}^{\prime}\mathbf{P}}^2 + \|\mathbf{q}\|_{\mathbf{W}^{\prime}} \|\mathbf{G}\|_{\mathbf{W}^{\prime}\mathbf{P}} \iff
    \frac{d}{dt} \| \mathbf{q} \|_{\mathbf{W}^{\prime}\mathbf{P}} \leq \mu \|\mathbf{q}\|_{\mathbf{W}^{\prime}\mathbf{P}} +  \|\mathbf{G}\|_{\mathbf{W}^{\prime}\mathbf{P}}. \\
\end{split}
    \end{align}
    Combining Gr\"onwall's Lemma and Duhamel's principle gives,
$$
\| \mathbf{q}\|_{\mathbf{W}^{\prime}\mathbf{P}} \leq Ke^{\mu t}\biggl(\| \mathbf{f}\|_{\mathbf{W}^{\prime}\mathbf{P}} +  \max\limits_{\tau \in [0, t]}\|\mathbf{G(\tau})\|_{\mathbf{W}^{\prime}\mathbf{P}}\biggl).
$$
\end{proof}
\revi{ We can also prove that the  semi-discrete total absolute vorticity is globally conserved in a periodic domain.
\begin{theorem}\label{theo:conservation_of_vorticity}
 Consider the semi-discrete DP SBP FD approximation \eqref{eqn:nl2dSWE2D} with periodic BCs.  At time $t \ge 0$, let the semi-discrete total absolute vorticity be denoted by $\mathcal{W}(t)= \sum^{N}_{i,j=1} \omega_{ij} \Delta x \Delta y$ where $\boldsymbol{\omega}=\mathbf{D}_{-x}{\mathbf{v}} - \mathbf{D}_{-y}{\mathbf{u}} + {f_c}$ is the semi-discrete absolute vorticity. If the Coriolis frequency $f_c$ is time-invariant, {that is $df_c/dt =0$}, then we have
\begin{align}
    \frac{d \mathcal{W}(t)}{dt} = 0,  \quad \forall t\ge 0.
\end{align}
\end{theorem}
\begin{proof}
Consider
\begin{align}
    \frac{d \boldsymbol{\omega}}{dt} = \mathbf{D}_{-x}\frac{d\mathbf{v}}{dt} - \mathbf{D}_{-y} \frac{d\mathbf{u}}{dt} + \frac{d f_c}{dt} = \mathbf{D}_{-x}\frac{d\mathbf{v}}{dt} - \mathbf{D}_{-y} \frac{d\mathbf{u}}{dt} ,
\end{align}
and
\begin{align}
    \frac{d \mathcal{W}(t)}{dt} = \left(\mathbf{1}, \frac{d \boldsymbol{\omega}}{dt}\right)_P = \left(\mathbf{1},\mathbf{D}_{-x}\frac{d\mathbf{v}}{dt}\right)_P - \left(\mathbf{1},\mathbf{D}_{-y} \frac{d\mathbf{u}}{dt} \right)_P.
\end{align}
Using the DP SBP property $\mathbf{D}_{-x}+\mathbf{D}_{+x}^T = \mathbf{0}$, $\mathbf{D}_{-y}+\mathbf{D}_{+y}^T = \mathbf{0}$ then we have
\begin{align}
    \frac{d \mathcal{W}(t)}{dt}  = -\left(\mathbf{D}_{+x}\mathbf{1},\frac{d\mathbf{v}}{dt}\right)_P + \left(\mathbf{D}_{+y}\mathbf{1}, \frac{d\mathbf{u}}{dt} \right)_P =0, \quad \forall\, t\ge 0,
\end{align}
where $\mathbf{1} = \left(1, 1, \ldots, 1\right)\in \mathbb{R}^{(N+1)^2}$ with $\mathbf{D}_{+x}\mathbf{1} =\vb{0}$,  $\mathbf{D}_{+y}\mathbf{1} =\vb{0}$. 
The proof is complete. 
\end{proof}
}

\begin{remark}
    While we are able to prove mass and potential enstrophy conservation in the continuous setting, the semi-discrete approximation of these quantities require local conservation. We will not prove it, but we have already shown numerically that their conservation at the discrete level is maintained up to discretisation errors.
\end{remark}

\bibliographystyle{elsarticle-num} 
\bibliography{refs}





\end{document}